\documentclass[12pt,a4paper,leqno]{article}

\usepackage{latexsym}
\usepackage{amssymb,exscale}
\usepackage[centertags]{amsmath}
\usepackage{amsthm}
\usepackage[all]{xy}
\usepackage{graphicx}
\usepackage{multicol}
\usepackage{makeidx}

\numberwithin{equation}{section}

\swapnumbers
\theoremstyle{definition}
\newtheorem{theorem}[equation]{Theorem}
\newtheorem{lemma}[equation]{Lemma}
\newtheorem{corollary}[equation]{Corollary}
\newtheorem{definition}[equation]{Definition}
\newtheorem{remark}[equation]{Remark}

\newcommand{\cpHs}{CP HS}
\newcommand{\cpHss}{CPs HS}
\newcommand{\hcpHs}{HCP HS}
\newcommand{\hcpHss}{HCPs HS}
\newcommand{\cpmc}{CP MC}
\newcommand{\hcpmc}{HCP MC}
\newcommand{\cpmcs}{CPs MC}
\newcommand{\hcpmcs}{HCPs MC}
\newcommand{\cpps}{CP PS}
\newcommand{\cppss}{CPs PS}
\newcommand{\hcpps}{HCP PS}
\newcommand{\hcppss}{HCPs PS}
\renewcommand{\Re}{\operatorname{Re}}
\renewcommand{\Im}{\operatorname{Im}}
\DeclareMathOperator{\codim}{codim}
\DeclareMathOperator{\tr}{tr}
\DeclareMathOperator{\mes}{mes}

\DeclareMathOperator{\dist}{dist}
\DeclareMathOperator{\const}{const}
\DeclareMathOperator{\sgn}{sgn}
\DeclareMathOperator{\supp}{supp}
\renewcommand{\(}{\bigl(}
\renewcommand{\)}{\bigr)\vphantom{)}}
\newcommand{\al}{\alpha}
\newcommand{\be}{\beta}
\newcommand{\de}{\delta}
\newcommand{\ga}{\gamma}
\newcommand{\la}{\lambda}
\newcommand{\Om}{\Omega}
\newcommand{\om}{\omega}
\newcommand{\R}{\mathbb R}
\newcommand{\C}{\mathbb C}
\newcommand{\Z}{\mathbb Z}
\newcommand{\N}{\mathbb N}
\newcommand{\A}{\mathcal A}
\newcommand{\B}{\mathcal B}

\newcommand{\M}{\mathcal M}
\newcommand{\F}{\mathcal F}
\newcommand{\cE}{\mathcal E}
\newcommand{\cN}{\mathcal N}
\newcommand{\ccP}{\mathcal P}
\newcommand{\cC}{\mathcal C}
\newcommand{\One}{\mathbf1}
\newcommand{\ti}{\tilde}
\newcommand{\sif}{$\sigma$-field}
\newcommand{\ip}[2]{\langle#1,#2\rangle}
\newcommand{\Ex}{\mathbb E}
\renewcommand{\Pr}[1]{\,\mathbb P\,\(\,#1\,\)\,}
\newcommand{\cP}[2]{\mathbb{P}\,\(\,#1\,\big|\,#2\,\)\,}
\newcommand{\pd}{\partial}
\newcommand{\ot}{\otimes}
\renewcommand{\phi}{\varphi}
\newcommand{\eps}{\varepsilon}
\def\emailwww#1#2{\par\quad {\tt #1}\par\quad {\tt #2}\medskip}
\newenvironment{abbrevlist}{\begin{list}{}
 {\setlength{\parsep}{0ex}
  \setlength{\itemsep}{0ex}
  \setlength{\leftmargin}{1.1cm}
  \setlength{\labelwidth}{1.1cm}
  \setlength{\labelsep}{0cm}
 }}
 {\end{list}}
\makeatletter
  {\end{multicols}}%
\makeatother
\makeatletter
\renewcommand*\l@section[2]{%
  \ifnum \c@tocdepth >\z@
    \addpenalty\@secpenalty
    \addvspace{0.25em \@plus\p@}%
    \setlength\@tempdima{2.5em}%
    \begingroup
      \parindent \z@ \rightskip \@pnumwidth
      \parfillskip -\@pnumwidth
      \leavevmode \bfseries
      \advance\leftskip\@tempdima
      \hskip -\leftskip
      #1\nobreak\hfil \nobreak\hb@xt@\@pnumwidth{\hss #2}\par
    \endgroup
  \fi}
\renewcommand*\numberline[1]{\hb@xt@\@tempdima{\hfil#1\hskip1em}}
\makeatother

\begin{document}

\title{Non-Isomorphic Product Systems}

\author{Boris Tsirelson}

\date{}

\maketitle

\stepcounter{footnote}
\footnotetext{%
 Supported by the \textsc{israel science foundation} founded by the
 Israel Academy of Sciences and Humanities.}

\begin{abstract}
Uncountably many mutually non-isomorphic product systems (that is,
continuous tensor products of Hilbert spaces) of types $ II_0 $ and $
III $ are constructed by probabilistic means (random sets and
off-white noises), answering four questions of W.~Arveson.
\end{abstract}

\setcounter{tocdepth}{1}
\tableofcontents

\vfill

\section*{Introduction}
\addcontentsline{toc}{section}{Introduction}
\smallskip

\hfill\parbox{9cm}{%
It is very difficult to construct examples of types $ II $ and $ III
$.\\
\mbox{}\hfill W.~Arveson \cite[p.~166]{Ar99}
}

\medskip

Product systems appeared in Arveson \cite{Ar89} as a tool for
investigating quantum dynamics (namely, $ E_0 $-semigroups). Fruitful
interplay between continuous (tensor) products of various spaces,
started in Tsirelson and Vershik \cite{TV}, is crucial for the
present work. For Hilbert spaces, the idea goes back to Araki and
Woods \cite{AW}; for probability spaces --- to Feldman
\cite{Fe}; measure classes\footnote{%
 Called `measure-type spaces'; the term `measure class' is borrowed
 from Arveson.}
appeared in \cite{TV}. Also, the word `continuous' in `continuous
product' may refer to various bases, not just $ [0,\infty)
$. Especially, complete Boolean algebras are used in \cite{AW} and
\cite{Fe}. Relations between product systems and continuous products
are clarified in Sect.~1.

Main results: Existence of a continuum of mutually non-isomorphic
product systems of type $ II_0 $ (Theorem \ref{6.6}), both symmetric
and asymmetric (Theorems \ref{7.1}, \ref{7.2}), and of type $ III $
(Theorem \ref{13.6}). Four questions asked by Arveson \cite{Ar89},
\cite{Ar96}, \cite{Ar97}, \cite{Ar99} are thus answered.

Main ingredients: Anatoly Vershik's idea of a continuous product of
measure classes; an idea of Jonathan Warren (private communication,
Nov.~1999) of constructing a continuous product of measure classes out
of a given random set; and the author's ideas about appropriate
invariants (Sect.~2), random sets (Sect.~5--7), FHS spaces and
off-white noises (Sect.~8--13).

\begin{center}\textsc{Abbreviations}\end{center}

\noindent
\begin{minipage}{8cm}
\begin{abbrevlist}
\item[CP\hfill] Continuous product\dots
\item[CPs\hfill] Continuous products\dots
\item[HCP\hfill] Homogeneous continuous product\dots
\end{abbrevlist}
\end{minipage}
\hfill
\begin{minipage}{5.5cm}
\begin{abbrevlist}
\item[HS\hfill] \dots of Hilbert spaces
\item[MC\hfill] \dots of measure classes
\item[PS\hfill] \dots of probability spaces
\end{abbrevlist}
\end{minipage}

\medskip\noindent
Thus, `\hcpHss' means: homogeneous continuous products of Hilbert
spaces.

\section{Basic notions}

\begin{sloppypar}
An oversimplified caricature of our objects is the Hilbert space $ H =
L_2 \( \{-1,+1\}^\Z \) $ of all square integrable measurable functions
on the probability space $ \{-1,+1\}^\Z $ of all two-sided sequences
of $ \pm1 $, the product of $ \Z $ copies of $ \{-1,+1\} $, both $ -1
$ and $ +1 $ being of probability $ 1/2 $. Relevant structures on $ H
$ are the shift operation $ \theta : H \to H $ and decompositions $ H
= H_{-\infty,\infty} = H_{-\infty,t-1} \otimes H_{t,\infty} $, where $
H_{s,t} = L_2 \( \{-1,+1\}^{\Z\cap[s,t]} \) $. Shifts $ \theta^u $ tie
together these decompositions by sending $ H_{-\infty,t-1} $ to $
H_{-\infty,t+u-1} $ and $ H_{t,\infty} $ to $ H_{t+u,\infty} $. The
caricature uses discrete time ($ \Z $), while our true objects use
continuous time ($ \R $), which leads to a deeper theory.
\end{sloppypar}

\begin{quote}
Throughout, either by assumption or by construction, all Hilbert
spaces will be separable.\footnote{%
 For now, they are complex (that is, over $ \C $), but in Sect.~2 we
 will return to the point.}%
\index{real versus complex}\index{complex versus real}
\index{Hilbert space: always separable}
\end{quote}

We consider a one-parameter group of unitary operators on the tensor
product of two Hilbert spaces $ H_- $ and $ H_+ $ (denoted also $
H_{-\infty,0} $ and $ H_{0,\infty} $) and corresponding automorphisms
$ \al_t $ of the operator algebra $ \B ( H_- \ot H_+ ) $,
\begin{equation}\label{1}
\al_t (A) = e^{itX} A e^{-itX} \, , \quad A \in \B ( H_- \ot H_+ ) \,
, \quad t \in \R \, ,
\end{equation}
satisfying two equivalent conditions:
\begin{gather}\label{2}
\al_t ( \B(H_-) \ot \One ) \subset \B(H_-) \ot \One \quad \text{for }
 t \le 0 \, ; \\
\al_t ( \One \ot \B(H_+) ) \subset \One \ot \B(H_+) \quad \text{for }
 t \ge 0 \, ;
\end{gather}
here $ \B(H_-) \ot \One $ is the algebra (in fact, type $ I_\infty $
factor) of all operators $ A \ot \One $ on $ H_- \ot H_+ $, where $ A
\in \B(H_-) $ and $ \One $ is the identity operator on $ H_+
$.\index{zz@$ \One $, identity operator on a Hilbert space}

For each $ t $,
\begin{equation}\label{3}
\A_{-\infty,\infty} = \A_{-\infty,t} \ot \A_{t,\infty} \, ,
\end{equation} 
where
\begin{gather*}
\A_{-\infty,t} = \al_t ( \B(H_-) \ot \One ) \, , \quad
 \A_{t,\infty} = \al_t ( \One \ot \B(H_-) ) \, ; \\
\A_{-\infty,s} \subset \A_{-\infty,t} \text{ and } \A_{s,\infty}
\supset \A_{t,\infty} \text{ for } s \le t \, ;
\end{gather*}
by \eqref{3} we mean that $ \A_{-\infty,t} $ and $ \A_{t,\infty} $
mutually commute and generate $ \A_{-\infty,\infty} $ as a von Neumann
algebra. Moreover, introducing
\[
\A_{s,t} = \A_{-\infty,t} \cap \A_{s,\infty}
\]
we get
\[
\A_{r,t} = \A_{r,s} \ot \A_{s,t} \quad \text{for } -\infty < r < s < t
< \infty \, ,
\]
see \cite[(1.2)]{Ar}. If, in addition, both intersections $ \cap_t
\A_{-\infty,t} $ and $ \cap_t \A_{t,\infty} $ are trivial (that is,
equal $ \C \cdot \One $), then the pair $ \( (e^{itX})_{t\in\R},
\B(H_-)\ot\One \) $ is a \emph{history}\index{history}
as defined by Arveson \cite[Sect.~1.2]{Ar}.

Endomorphisms $ \al_t^- : \B(H_-) \to \B(H_-) $, $ \al_t^+ : \B(H_+)
\to \B(H_+) $ defined for $ t \ge 0 $ by
\begin{gather*}
\al_{-t} ( A \ot \One ) = \al_t^- (A) \ot \One \quad \text{for } A \in
 \B(H_-) \, , \\
\al_t ( \One \ot A ) = \One \ot \al_t^+ (A) \quad \text{for } A \in
 \B(H_+) \, ,
\end{gather*}
form so-called $ E_0 $-semigroups\index{E0@$ E_0 $-semigroup} $ \al^-,
\al^+ $,
see \cite[Sect.~1.2]{Ar}, irrespective of the additional condition of
tail triviality (the latter means that the $ E_0 $-semigroups are
pure).

It is well-known (see \cite[Sect.~3.1]{Ar}) that each $ E_0
$-semigroup (pure or not) leads to a product
system\index{product system}
(of Hilbert
spaces). Namely, for $ t > 0 $, the set $ \cE(t) = \{ T \in \B(H_+) :
\forall A \in \B(H_+) \; \al_t^+ (A) T = TA \} $ is actually a Hilbert
space, and multiplication acts like tensoring in the sense that $ ST =
W_{s,t} ( S \ot T ) $ for some (unique) unitary operator $ W_{s,t} :
\cE(s) \ot \cE(t) \to \cE(s+t) $.

On the other hand, our objects are a bit more structured than just
histories; the decomposition of our Hilbert space into $ H_- \ot H_+ $
is only an example. The general decomposition is, roughly speaking,
\[
H_{-\infty,\infty} = H_{-\infty,t_1} \ot H_{t_1,t_2} \ot \dots \ot
H_{t_{n-1},t_n} \ot H_{t_n,\infty}
\]
for $ -\infty < t_1 < \dots < t_n < \infty $. Here is a rigorous
definition.

\begin{definition}\label{1.1}
(a) 
A \emph{continuous product of Hilbert spaces} (`\cpHs', for
short)\index{CP HS@\cpHs\ = continuous product of Hilbert spaces}
is $
\( (H_{s,t})_{-\infty\le s<t\le\infty},
\linebreak[0]
(W_{r,s,t})_{-\infty\le r<s<t\le\infty} \) $ where each $ H_{s,t} $
is a Hilbert space, $ \dim H_{s,t} > 0 $, and each $ W_{r,s,t} :
H_{r,s} \otimes H_{s,t} \to H_{r,t} $ is a unitary operator,
satisfying the condition

(a1) [associativity] the diagram 
\[
\xymatrix{
 & H_{r,s} \otimes H_{s,t} \otimes H_{t,u}
  \ar[dl]_{W_{r,s,t}\otimes\One}
  \ar[dr]^{\One\otimes W_{s,t,u}} \\
H_{r,t} \otimes H_{t,u}
  \ar[dr]_{W_{r,t,u}} & &
 H_{r,s} \otimes H_{s,u}
  \ar[dl]^{W_{r,s,u}} \\
 & H_{r,u}
}
\]
is commutative whenever $ -\infty\le r<s<t<u\le\infty $.

The \cpHs\ is \emph{tail trivial,} if it satisfies the condition

\begin{sloppypar}
(a2) [tail triviality]:
\[
\bigcap_{t\in\R} \A_{-\infty,t} = \C \cdot \One = \bigcap_{t\in\R}
\A_{t,\infty} \, ,
\]
where $ \A_{-\infty,t} = \{ W_{-\infty,t,\infty} ( A \ot \One )
W^{-1}_{-\infty,t,\infty} : A \in \B(H_{-\infty,t}) \} $, $
\A_{t,\infty} = \{ W_{-\infty,t,\infty} ( \One \ot A ) 
W^{-1}_{-\infty,t,\infty} : A \in \B(H_{t,\infty}) \} $.
\end{sloppypar}

(b) An isomorphism\index{isomorphism!of \cpHss}
of one \cpHs\ $ \( (H_{s,t})_{s<t},
(W_{r,s,t})_{r<s<t} \) $ to another \cpHs\ $ \( (H'_{s,t})_{s<t},
(W'_{r,s,t})_{r<s<t} \) $ is a family $ (U_{s,t})_{-\infty\le
s<t\le\infty} $ of unitary operators $ U_{s,t} : H_{s,t} \to
H'_{s,t} $ such that the diagram
\[
\xymatrix@C+3cm{
H_{r,s} \otimes H_{s,t}
  \ar[r]^{W_{r,s,t}}
  \ar[d]_{U_{r,s} \otimes U_{s,t}} &
 H_{r,t}
  \ar[d]^{U_{r,t}} \\
H'_{r,s} \otimes H'_{s,t}
  \ar[r]^{W'_{r,s,t}} &
H'_{r,t}
}
\]
is commutative whenever $ -\infty \le r < s < t \le \infty $.
\end{definition}

\begin{sloppypar}
Given a \cpHs\ $ \( (H_{s,t})_{s<t}, (W_{r,s,t})_{r<s<t} \) $ and $ u
\in \R $, the object $ \( (H_{s+u,t+u})_{s<t},
\linebreak[0]
(W_{r+u,s+u,t+u})_{r<s<t} \) $ is also a \cpHs; an isomorphism $
\theta^u $ from the former to the latter may be called a shift (by $ u
$) of the given \cpHs. A one-parameter group of shifts turns a \cpHs\
into a homogeneous \cpHs, as defined below.
\end{sloppypar}

\begin{definition}\label{1.2}
\begin{sloppypar}
(a) A \emph{homogeneous continuous product of Hilbert 
spaces}\index{HCP HS@\hcpHs\ = homogeneous continuous product of Hilbert 
spaces}\index{homogeneous!CP HS@\cpHs}
(`\hcpHs', for short) is $ \( (H_{s,t})_{-\infty\le s<t\le\infty},
\linebreak[0]
(W_{r,s,t})_{-\infty\le r<s<t\le\infty}, \linebreak[0]
(\theta^u_{s,t})_{u\in\R,-\infty\le s<t\le\infty} \) $
where $ \(
(H_{s,t})_{s<t}, (W_{r,s,t})_{r<s<t} \) $ is a \cpHs, and $
(\theta^u_{s,t})_{u;s<t} $ is a family of unitary operators
\[
\theta^u_{s,t} : H_{s,t} \to H_{s+u,t+u}
\]
such that $ \theta^0_{s,t} = \One_{s,t} $, and
\end{sloppypar}

(a1) the diagram
\[
\xymatrix@C+3cm{
H_{r,s} \otimes H_{s,t}
  \ar[r]^{W_{r,s,t}}
  \ar[d]_{\theta^u_{r,s} \otimes \theta^u_{s,t}} &
 H_{r,t}
  \ar[d]^{\theta^u_{r,t}} \\
H_{r+u,s+u} \otimes H_{s+u,t+u}
  \ar[r]^{W_{r+u,s+u,t+u}} &
H_{r+u,t+u}
}
\]
is commutative whenever $ u \in \R $, $ -\infty \le r < s < t \le
\infty $;

(a2) the diagram
\[
\xymatrix{
H_{s,t}
  \ar[rr]^{\theta^{u+v}_{s,t}}
  \ar[rd]_{\theta^u_{s,t}} & &
 H_{s+u+v,t+u+v} \\
& H_{s+u,t+u}
  \ar[ru]_{\theta^v_{s+u,t+u}}
}
\]
is commutative whenever $ u,v \in \R $, $ -\infty \le s < t \le \infty
$;

(a3) operators $ \theta^u_{-\infty,\infty} $ are a strongly continuous
unitary group; that is, $ \theta^u_{-\infty,\infty} = e^{iuX} $ for
some self-adjoint operator $ X $ on $ H_{-\infty,\infty} $.

The \hcpHs\ is called tail trivial, if the underlying \cpHs\ is tail
trivial.

\begin{sloppypar}
(b) An isomorphism of one \hcpHs\ $ \( (H_{s,t})_{s<t},
(W_{r,s,t})_{r<s<t}, (\theta^u_{s,t})_{u;s<t} \) $ to another \hcpHs\ $
\( (H'_{s,t})_{s<t}, (W'_{r,s,t})_{r<s<t}, (\theta^{\prime u}_{s,t})_{u;s<t}
\) $ is an isomorphism $ (U_{s,t})_{s<t} $ of the \cpHs\ $ \(
(H_{s,t})_{s<t}, (W_{r,s,t})_{r<s<t} \) $ to the \cpHs\ $ \(
(H'_{s,t})_{s<t}, (W'_{r,s,t})_{r<s<t} \) $ (as defined by
\ref{1.1}(b)) such that the diagram
\[
\xymatrix@C+3cm{
H_{s,t}
  \ar[r]^{\theta^u_{s,t}}
  \ar[d]_{U_{s,t}} &
 H_{s+u,t+u}
  \ar[d]^{U_{s+u,t+u}} \\
H'_{s,t}
  \ar[r]^{\theta^{\prime u}_{s,t}} &
H'_{s+u,t+u}
}
\]
is commutative whenever $ u \in \R $, $ -\infty \le s < t \le \infty
$.
\end{sloppypar}
\end{definition}

Given an \hcpHs, we may take $ H_- = H_{-\infty,0} $, $ H_+ =
H_{0,\infty} $, and consider the unitary group $ (
W^{-1}_{-\infty,0,\infty} e^{iuX} W_{-\infty,0,\infty} )_{u\in\R} $ on
$ H_- \ot H_+ $, then \eqref{1} and \eqref{2} are satisfied; indeed,
\eqref{2} follows from the lemma below.

\[
\begin{gathered}\includegraphics{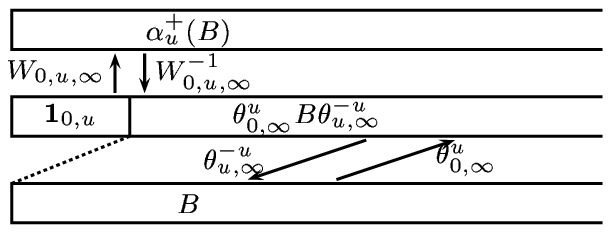}\end{gathered}
\]

\begin{lemma}\label{1.3}
For every \hcpHs\ $ \( (H_{s,t})_{s<t}, (W_{r,s,t})_{r<s<t},
(\theta^u_{s,t})_{u;s<t} \) $, every $ u \ge 0 $, $ A \in
\B(H_{-\infty,0}) $, $ B \in \B(H_{0,\infty}) $,
\[
\al_{-u} ( A \ot \One ) = \al_u^- (A) \ot \One \, , \quad
\al_u ( \One \ot B ) = \One \ot \al_u^+ (B) \, ,
\]
where
\begin{gather*}
\al_u^- (A) = W_{-\infty,u,0} ( \theta^u_{-\infty,0} A
 \theta^{-u}_{-\infty,u} \ot \One_{u,0} ) W^{-1}_{-\infty,u,0} \, , \\
\al_u^+ (B) = W_{0,u,\infty} ( \One_{0,u} \ot \theta^u_{0,\infty} B
 \theta^{-u}_{u,\infty} ) W^{-1}_{0,u,\infty} \, .
\end{gather*}
\end{lemma}

\begin{proof}
Only the formula for $ \al^+ $ will be proven; the formula for $ \al^-
$ is similar. We have
\begin{gather*}
\al_u ( \One \ot B ) = W^{-1}_{-\infty,0,\infty} e^{iuX}
 W_{-\infty,0,\infty} ( \One \ot B ) W^{-1}_{-\infty,0,\infty}
 e^{-iuX} W_{-\infty,0,\infty} \, , \\
e^{iuX} W_{-\infty,0,\infty} = \theta^u_{-\infty,\infty}
 W_{-\infty,0,\infty} = W_{-\infty,0,\infty} ( \theta^u_{-\infty,0}
 \ot \theta^u_{0,\infty} )
\end{gather*}
by \ref{1.2}(c,a);
\[
W_{-\infty,0,\infty} ( \One_{-\infty,0} \ot W_{0,u,\infty} ) =
W_{-\infty,u,\infty} ( W_{-\infty,0,u} \ot \One_{u,\infty} )
\]
by \ref{1.1}(a1); thus,
\begin{multline*}
W^{-1}_{-\infty,0,\infty} e^{iuX} W_{-\infty,0,\infty} =
 W^{-1}_{-\infty,0,\infty} W_{-\infty,u,\infty} ( \theta^u_{-\infty,0}
 \ot \theta^u_{0,\infty} ) = \\
= ( \One_{-\infty,0} \ot W_{0,u,\infty} ) (
 W_{-\infty,0,u} \ot \One_{u,\infty} )^{-1} ( \theta^u_{-\infty,0}
 \ot \theta^u_{0,\infty} ) \, .
\end{multline*}
On the other hand,
\[
( \theta^u_{-\infty,0} \ot \theta^u_{0,\infty} ) ( \One_{-\infty,0} \ot
B ) ( \theta^u_{-\infty,0} \ot \theta^u_{0,\infty} )^{-1} =
\One_{-\infty,u} \ot ( \theta^u_{0,\infty} B (\theta^u_{0,\infty})^{-1}
) \, ;
\]
thus
\begin{multline*}
( W_{-\infty,0,u} \ot \One_{u,\infty} )^{-1} ( \theta^u_{-\infty,0}
 \ot \theta^u_{0,\infty} ) ( \One_{-\infty,0} \ot B ) (
 \theta^u_{-\infty,0} \ot \theta^u_{0,\infty} )^{-1} ( W_{-\infty,0,u}
 \ot \One_{u,\infty} ) \\
= ( W_{-\infty,0,u} \ot \One_{u,\infty} )^{-1} ( \One_{-\infty,u} \ot
 \theta^u_{0,\infty} B (\theta^u_{0,\infty})^{-1} ) ( W_{-\infty,0,u}
 \ot \One_{u,\infty} ) = \\
= \One_{-\infty,0} \ot \One_{0,u} \ot \theta^u_{0,\infty} B
 (\theta^u_{0,\infty})^{-1} \, ,
\end{multline*}
and so,
\begin{multline*}
\al_u ( \One \ot B ) = ( \One_{-\infty,0} \ot W_{0,u,\infty} ) (
 \One_{-\infty,0} \ot \One_{0,u} \ot \theta^u_{0,\infty} B
 (\theta^u_{0,\infty})^{-1} ) ( \One_{-\infty,0} \ot W_{0,u,\infty}
 )^{-1} \\
= \One_{-\infty,0} \ot W_{0,u,\infty} ( \One_{0,u} \ot
 \theta^u_{0,\infty} B (\theta^u_{0,\infty})^{-1} )
 W^{-1}_{0,u,\infty} = \One_{-\infty,0} \ot \al_u^+ (B) \, .
\end{multline*}
\end{proof}

\begin{corollary}\label{1.4}
Let $ \( (H_{s,t})_{s<t}, (W_{r,s,t})_{r<s<t},
(\theta^u_{s,t})_{u;s<t} \) $ be a tail trivial \hcpHs. Then the
unitary group $ \( W^{-1}_{-\infty,0,\infty} \theta^t_{-\infty,\infty}
W_{-\infty,0,\infty} )_{t\in\R} $ on $ H_{-\infty,0} \ot H_{0,\infty}
$ and the algebra $ \B(H_{-\infty,0}) \ot \One $ form a 
history\index{history!out of \hcpHs}
\textup{(as defined in \cite[Sect.~1.2]{Ar}).}
\end{corollary}

\begin{proof}
Property \eqref{2} follows from Lemma \ref{1.3}; tail triviality
follows from Definition \ref{1.1}(a2).
\end{proof}

Every \hcpHs\ (tail trivial or not) leads (via its $ E_0 $-semigroup $
\al^+ $) to a product system\index{product system!out of \hcpHs}
$ (\cE(t))_{t>0} $. As was said, an element of $ \cE(t) $ is an
operator $ T \in \B(H_+) $ such that $ \al_t (A) T = TA $ for all $ A
\in \B(H_+) $. According to Lemma \ref{1.3} we rewrite it as
\[
W_{0,t,\infty} \underbrace{ ( \One_{0,t} \ot \theta^t_{0,\infty} ) (
\One_{0,t} \ot A ) ( \One_{0,t} \ot \theta^t_{0,\infty} )^{-1} }_{ (
\One_{0,t} \ot \theta^t_{0,\infty} A \theta^{-t}_{t,\infty} ) }
W^{-1}_{0,t,\infty} T = T A \, ,
\]
or just
\[
( \One_{0,t} \ot A ) \ti T = \ti T A
\]
where
\[
\ti T = ( \One_{0,t} \ot \theta^{-t}_{t,\infty} ) W^{-1}_{0,t,\infty}
T, \quad \ti T : H_{0,\infty} \to H_{0,t} \ot H_{0,\infty} \, .
\]
Applying it to a one-dimensional operator $ A $ we see that $ \ti T x
\in H_{0,t} \ot x $ for each $ x \in H_{0,\infty} $; therefore $ \ti T
$ must be of the form
\[
\ti T x = y \ot x
\]
for some $ y \in H_{0,t} $. Also, for every $ y $, such $ \ti T $
belongs to $ \cE(t) $. We have a natural unitary map from $ H_{0,t} $
to $ \cE(t) $ (the latter consists of operators, but is isometric to a
Hilbert space); namely, $ y \mapsto Ty $,
\[
T_y x = W_{0,t,\infty} ( \One_{0,t} \ot \theta^t_{0,\infty} ) ( y \ot
x ) = W_{0,t,\infty} ( y \ot \theta^t_{0,\infty} x ) \, .
\]

Let $ \ti T_1 \in \cE(t_1) $, $ \ti T_2 \in \cE(t_2) $, $ \ti T_1 x =
y_1 \ot x $, $ \ti T_2 x = y_2 \ot x $, $ y_1 \in H_{0,t_1} $, $ y_2
\in H_{0,t_2} $, and $ T = T_2 T_1 $. Then
\begin{multline*}
Tx = T_2 T_1 x = W_{0,t_2,\infty} ( y_2 \ot \theta^{t_2}_{0,\infty}
 T_1 x ) = \\
= W_{0,t_2,\infty} \( y_2 \ot \underbrace{ W_{t_2,t_1+t_2,\infty} (
 \theta^{t_2}_{0,t_1} \ot \theta^{t_2}_{t_1,\infty} ) }_{
 \theta^{t_2}_{0,\infty} W_{0,t_1,\infty} }
 ( y_1 \ot \theta^{t_1}_{0,\infty} x ) \) = \\
= W_{0,t_2,\infty} \( y_2 \ot W_{t_2,t_1+t_2,\infty} (
 \theta^{t_2}_{0,t_1} y_1 \ot \theta^{t_1+t_2}_{0,\infty} x ) \) = \\
= W_{0,t_1+t_2,\infty} \( W_{0,t_2,t_1+t_2} ( y_2 \ot
 \theta^{t_2}_{0,t_1} y_1 ) \ot \theta^{t_1+t_2}_{0,\infty} x \) \, ,
\end{multline*}
that is,
\[
\ti T x = y \ot x \, , \quad \text{where } y = W_{0,t_2,t_1+t_2} ( y_2
\ot \theta^{t_2}_{0,t_1} y_1 ) \, .
\]
The latter is the multiplication of the product system, transplanted
from $ (\cE(t))_{t>0} $ to $ (H_{0,t})_{t>0} $. We may transplant also
the Borel structure, which leads to the following conclusion.

\begin{theorem}\label{1.5}
Let $ \( (H_{s,t})_{s<t}, (W_{r,s,t})_{r<s<t},
(\theta^u_{s,t})_{u;s<t} \) $ be a \hcpHs. Then the disjoint union
\[
E = \biguplus_{t>0} H_{0,t} = \{ (t,x) : 0<t<\infty, x \in H_{0,t} \}
\]
is a product system\index{product system!out of \hcpHs}
\textup{(as defined in \cite[Def.~3.2.1]{Ar}),}
being equipped with the binary operation
\[
\( (s,x), (t,y) \) \mapsto \( s+t, W_{0,s,s+t} ( x \ot
\theta^{s}_{0,t} y ) \)
\]
and the Borel structure induced by such an embedding into operators on
$ H_{0,\infty} $:
\[
(s,x) \mapsto \( y \mapsto W_{0,s,\infty} ( x \ot
\theta^{s}_{0,\infty} y ) \) \, .
\]
\end{theorem}

Our definition of \cpHs\ stipulates Hilbert spaces $ H_{s,t} $ not only
for bounded intervals $ (s,t) $, but also for unbounded intervals $
(-\infty,t) $, $ (s,\infty) $ and $ (-\infty,\infty) $. These
unbounded intervals are more important than it would seem! We may
define a \emph{local \cpHs}\index{local \cpHs}
by excluding unbounded intervals from
Definition \ref{1.1}. Even under tail triviality (item (a2) of
Def.~\ref{1.1}), a global \cpHs\ is far from being uniquely determined
by its local \cpHs. Indeed, asymptotic behavior at infinity is
evidently missing in a local \cpHs. Strangely enough, something else,
seemingly local, is also missing. A global \cpHs\ induces a natural
Borel structure on its product system $ \uplus H_{0,t} $ (see Theorem
\ref{1.5}), but a local \cpHs\ does not, for the following reason.

Let $ \phi : \R \to \R $ be a function, additive in the sense that $
\phi(s+t) = \phi(s) + \phi(t) $ for all $ s,t $. It is well-known that
such $ \phi $ is either linear, $ \phi(t) = ct $, or quite pathologic
(non-measurable); the latter case is possible (which follows easily
from existence of a Hamel basis). Let $ \( (H_{s,t})_{s<t},
(W_{r,s,t})_{r<s<t}, (\theta^u_{s,t})_{u;s<t} \) $ be a local
\hcpHs.\footnote{%
 It is not clear, how to `localize' item (a3) of
 Def.~\ref{1.2}. Anyway, we may restrict ourselves to local \hcpHss\
 that can be extended to global \hcpHss.}
Define $ U_{s,t} \in \B(H_{s,t}) $ as the scalar operator $
U_{s,t} = e^{i\phi(t-s)} \One_{s,t} $; then $ (U_{s,t})_{s<t} $ is an
automorphism of the local \hcpHs\ (as defined by Def.~\ref{1.2}(b), but
excluding unbounded intervals). If the local \hcpHs\ is induced by a
global \hcpHs, then $ \uplus H_{0,t} $ has its natural Borel structure,
but the transformation of $ \uplus H_{0,t} $ that multiplies each $
H_{0,t} $ by $ e^{i\phi(t)} $ is not Borel, unless $ \phi $ is Borel
(therefore, linear). Such a pathologic automorphism of the local
\hcpHs\
cannot be extended to an automorphism of the global \hcpHs. Here is a
simple explanation: a family of scalar operators $ U_{s,t} = u_{s,t}
\One $, $ u_{s,t} \in \C $, $ |u_{s,t}| = 1 $, is an automorphism if
and only if $ u_{r,s} u_{s,t} = u_{r,t} $ and $ u_{r,s} = u_{r+t,s+t}
$. Thus, $ u_{r,s} = \overline{ f(r) } f(s) $ for some $ f $ such that
$ \overline{ f(r) } f(s) = \overline{ f(r+t) } f(s+t) $. If the case $
s = \infty $ is admitted, we get $ \overline{ f(r) } f(\infty) =
\overline{ f(r+t) } f(\infty) $, and $ f $ must be constant!

Or rather, $ f $ must be constant on $ (-\infty,\infty) $, not $
[-\infty,\infty] $. We may take two numbers $ u_-, u_+ \in \C $, $
|u_-| = 1 $, $ |u_+| = 1 $, and define $ u_{-\infty,t} = u_- $, $
u_{t,\infty} = u_+ $ for $ t \in \R $; of course, $ u_{-\infty,\infty}
= u_- u_+ $ and $ u_{r,s} = 1 $ for $ r,s \in \R $. Then $
(u_{s,t})_{s<t} $ is an automorphism of the global \hcpHs, whose
restriction to the local \hcpHs\ is trivial. We see that restriction (of
automorphisms, from global to local \hcpHs) is neither surjective nor
injective.

\begin{lemma}\label{1.6}
\begin{sloppypar}
Let $ \( (H_{s,t})_{-\infty\le s<t\le\infty},
(W_{r,s,t})_{-\infty\le r<s<t\le\infty},
(\theta^u_{s,t})_{u\in\R,-\infty\le s<t\le\infty} \) $ and
$ \( (H'_{s,t})_{-\infty\le s<t\le\infty},
(W'_{r,s,t})_{-\infty\le r<s<t\le\infty},
(\theta^{\prime u}_{s,t})_{u\in\R,-\infty\le s<t\le\infty} \) $
be two \hcpHss. If the corresponding product systems $ \uplus H_{0,t}
$, $ \uplus H'_{0,t} $ are 
isomorphic\index{isomorphism!of product systems and local \hcpHss}
then the corresponding local \hcpHss\ 
$ \( (H_{s,t})_{-\infty<s<t<\infty},
(W_{r,s,t})_{-\infty<r<s<t<\infty},
(\theta^u_{s,t})_{u\in\R,-\infty<s<t<\infty} \) $ and
$ \( (H'_{s,t})_{-\infty<s<t<\infty},
(W'_{r,s,t})_{-\infty<r<s<t<\infty},
(\theta^{\prime u}_{s,t})_{u\in\R,-\infty<s<t<\infty} \) $
are isomorphic.
\end{sloppypar}
\end{lemma}

\begin{proof}
Isomorphism between $ \uplus H_{0,t} $ and $ \uplus H'_{0,t} $ is a
family $ (U_{0,t})_{0<t<\infty} $ of unitary operators $ U_{0,t} :
H_{0,t} \to H'_{0,t} $ such that the diagram
\[
\xymatrix@C+3cm{
H_{0,s} \otimes H_{0,t}
  \ar[r]^{ W_{0,s,s+t} ( \One \ot \theta^s_{0,t} ) }
  \ar[d]_{U_{0,s} \otimes U_{0,t}} &
 H_{0,s+t}
  \ar[d]^{U_{0,s+t}} \\
H'_{0,s} \otimes H'_{0,t}
  \ar[r]^{ W'_{0,s,s+t} ( \One \ot \theta^{\prime s}_{0,t} ) } &
H'_{0,s+t}
}
\]
is commutative for all $ s,t \in (0,\infty) $. In addition, $
(U_{0,t})_{0<t<\infty} $ satisfies a Borel measurability condition,
but this one will not be used. We define $ U_{s,t} : H_{s,t} \to
H'_{s,t} $ for $ -\infty < s < t < \infty $ by
\[
U_{s,t} = \theta^{\prime s}_{0,t-s} U_{0,t-s} \theta^{-s}_{s,t} \, .
\]
We have
\begin{multline*}
U_{r,t} W_{r,s,t} = \theta^{\prime r}_{0,t-r} U_{0,t-r}
 \theta^{-r}_{r,t} W_{r,s,t} = \\
= \theta^{\prime r}_{0,t-r} U_{0,t-r} W_{0,s-r,t-r} (
 \theta^{-r}_{r,s} \ot \theta^{-r}_{s,t} ) = \\
= \theta^{\prime r}_{0,t-r} U_{0,t-r} W_{0,s-r,t-r} ( \One \ot
 \theta^{s-r}_{0,t-s} ) ( \theta^{-r}_{r,s} \ot \theta^{-s}_{s,t} ) =
 \\
= \theta^{\prime r}_{0,t-r} W'_{0,s-r,t-r} ( \One \ot \theta^{\prime
 s-r}_{0,t-s} ) ( U_{0,s-r} \ot U_{0,t-s} ) ( \theta^{-r}_{r,s} \ot
 \theta^{-s}_{s,t} ) = \\
= W'_{r,s,t} ( \theta^{\prime r}_{0,s-r} \ot \theta^{\prime
 r}_{s-r,t-r} ) ( \One \ot \theta^{\prime s-r}_{0,t-s} ) ( U_{0,s-r}
 \ot U_{0,t-s} ) ( \theta^{-r}_{r,s} \ot \theta^{-s}_{s,t} ) = \\
= W'_{r,s,t} ( \theta^{\prime r}_{0,s-r} U_{0,s-r} \theta^{-r}_{r,s}
 \ot \theta^{\prime r}_{s-r,t-r} \theta^{\prime s-r}_{0,t-s} U_{0,t-s}
 \theta^{-s}_{s,t} ) = \\
= W'_{r,s,t} ( U_{r,s} \ot U_{s,t} ) \, ,
\end{multline*}
which means commutativity of the diagram of Def.~\ref{1.1}(b). Also,
\begin{multline*}
U_{s+u,t+u} \theta^u_{s,t} = \theta^{\prime s+u}_{0,t-s} U_{0,t-s}
 \theta^{-s-u}_{s+u,t+u} \theta^u_{s,t} = \\
= \theta^{\prime s+u}_{0,t-s} U_{0,t-s} \theta^{-s}_{s,t} =
 \theta^{\prime u}_{s,t} \theta^{\prime s}_{0,t-s} U_{0,t-s}
 \theta^{-s}_{s,t} = \theta^{\prime u}_{s,t} U_{s,t} \, ,
\end{multline*}
which means commutativity of the diagram of Def.~\ref{1.2}(b).
\end{proof}

In principle, it would be more natural to define a \cpHs\ as local
(rather than global) and equipped with a Borel structure (rather than
the purely algebraic structure) in the spirit of
\cite[Def.~3.2.1]{Ar}. However, given the aim to construct
non-isomorphic product systems, we consider
\begin{quote}
global \hcpHss\ up to local
isomorphisms;\index{isomorphism!local, of global \hcpHss}
\end{quote}
that is, our objects are global \hcpHss, but our isomorphisms are
isomorphisms of corresponding local \cpHss\ (no Borel structure
stipulated). It means that

(a) we deal with a subclass of the class of all product systems,

(b) on that subclass we introduce an equivalence relation coarser than
isomorphism of product systems.

Therefore, existence of many nonequivalent objects within our
framework implies existence of (at least) equally many non-isomorphic
product systems.

It was noted (after Theorem \ref{1.5}) that the Borel structure on a
local \cpHs\ emerges via globality and homogeneity. It is natural to
ask what happens to the Borel structure in a local inhomogeneous
setup. The next result presents an answer which will not be used, and
I omit the proof; the reader may skip to Sect.~2.

\begin{theorem}
For every local \cpHs\ $ \( (H_{s,t})_{-\infty<s<t<\infty},
(W_{r,s,t})_{-\infty<r<s<t<\infty} \) $ there exists a Borel structure on
the disjoint union $ \uplus H_{s,t} $ of all these Hilbert spaces such
that

\textup{(a)} $ \uplus H_{s,t} $ is a standard Borel space;

\textup{(b)} the projection $ \uplus H_{s,t} \to \{ (s,t) : -\infty <
s < t < \infty \} $ \textup{(sending each $ H_{s,t} $ to $ (s,t) $)}
is a Borel map;

\textup{(c)} the map $ (x,y) \mapsto U_{r,s,t} (x \otimes y) $ from the Borel
subset $ \cup_{r,s,t} H_{r,s} \times H_{s,t} $ of $ \( \uplus H_{s,t}
\) \times \( \uplus H_{s,t} \) $ to $ \uplus H_{s,t} $ is a Borel map;

\textup{(d)} there exists a sequence $ (e_n) $ of Borel maps $ e_n : \{ (s,t) :
-\infty < s < t < \infty \} \to \uplus H_{s,t} $ such that the sequence $
\( e_1(s,t), e_2(s,t), \dots \) $ is an orthonormal basis of $ H_{s,t}
$ for each $ (s,t) $ satisfying $ \dim H_{s,t} = \infty $. Otherwise,
the first $ n = \dim H_{s,t} $ elements of the sequence are such a
basis, and other elements vanish.
\end{theorem}

Do not think that the Borel structure is unique. Rotating each $
H_{s,t} $ by a phase factor $ u_{s,t} = \overline{ f(s) } f(t) $, $ f
: \R \to \{ z \in \C : |z| = 1 \} $, we get an automorphism of the
local \cpHs, even if $ f $ is not a Borel function.

\section{Some invariants}

\smallskip

\hfill\parbox{9cm}{%
The fundamental problem in the theory of $ E_0 $-semigroups is to find
a complete set of computable invariants for cocycle conjugacy. There
is some optimism that such a classification is possible, but we are
far from achieving that goal.\\
\mbox{}\hfill W.~Arveson \cite[Sect.~3]{Ar}
}

\medskip

In the absence of a complete set of computable invariants, some useful
specific invariants are needed, if we want to present non-isomorphic
product systems (or $ E_0 $-semigroups). However, I consider
continuous products of Hilbert spaces (\cpHss) instead of product
systems. Moreover, I consider \cpHss\ over $ [0,1] $ (rather than $
(-\infty,\infty) $) and waive homogeneity. Ultimately, my objects are
`global \hcpHss\ up to local isomorphisms', as explained in
Sect.~1. But for now it does not harm to define invariants on a wider
class of objects. A local isomorphism of global \hcpHss\ evidently
induces an isomorphism of corresponding \cpHss\ over $ [0,1] 
$.\index{CP HS o@\cpHs\ over $ [0,1] $}
Such \cpHss\ (and their
isomorphisms)\index{isomorphism!of \cpHss!over $ [0,1] $}
are defined similarly to
Def.~\ref{1.1}; the only modification is that $ r,s,t \in [0,1] $
rather than $ [-\infty,\infty] $.

So, an arbitrary \cpHs\ $ (H_{s,t})_{0\le s<t\le1} $ is
considered. Instead of $ W_{r,s,t} : H_{r,s} \ot H_{s,t} \to H_{r,t} $
I'll write simply $ H_{r,s} \ot H_{s,t} = H_{r,t} $, glueing these
spaces together and treating $ W_{r,s,t} $ as identity maps. (The
reader could restore all needed $ W_{r,s,t} $, thus returning to the
rigorous but cumbersome style of Sect.~1. Associativity stipulated by
Def.~\ref{1.1}(a1) eliminates ambiguity.)

A subset of $ [0,1] $ will be called
\emph{elementary,}\index{elementary set}
if it is the
union of a finite number of intervals. We treat elementary sets modulo
finite sets; say, $ (0.3,0.4) \cup (0.6,0.7) \cup (0.7,0.8) $ and $
[0.3,0.4) \cup \{0.5\} \cup [0.6,0.8] $ are equivalent elementary
sets. The general form of an (equivalence class of) elementary set is
\begin{equation}\label{12}
E = (s_1,t_1) \cup \dots \cup (s_n,t_n) \, , \quad 0 \le s_1 < t_1 <
\dots < s_n < t_n \le 1 \, ,
\end{equation}
$ n = 0,1,2,\dots $ ($ n=0 $ means $ E = \emptyset $). These
equivalence classes are a Boolean
algebra\index{Boolean algebra of elementary sets}
$ \cE $.

We define for each $ E \in \cE $ a Hilbert space
\[
H_E = H_{s_1,t_1} \ot \dots \ot H_{s_n,t_n}
\]
and observe that
\[
H_{E\cup F} = H_E \ot H_F \quad \text{whenever } E,F \in \cE, E\cap F
= \emptyset \, .
\]
In particular,
\[
H_{0,1} = H_{[0,1]} = H_E \ot H_{[0,1]\setminus E} \, ;
\]
accordingly, we define
\[
\A_E = \B (H_E) \ot \One_{[0,1]\setminus E} \subset \A_{0,1} = \B
(H_{0,1})
\]
(here $ \One_{[0,1]\setminus E} $ is the identity operator on $
H_{[0,1]\setminus E} $), getting
\begin{gather*}
\A_E \ot \A_F = \A_{E\cup F} \quad \text{for } E \cap F = \emptyset \,
 , \\
\A'_E = \A_{[0,1]\setminus E} \, .
\end{gather*}

\begin{definition}\label{2.1}
A sequence $ (E_k) $ of (equivalence classes of) elementary sets $
E_1, E_2, \dots \in \cE $ is 
\emph{infinitesimal}\index{infinitesimal sequence of elementary sets}
for the \cpHs\ $ (H_{s,t})_{0\le s<t\le1} $, if
\[
\sup_{A\in\A_k,\|A\|\le1} | \tr (AR) | \to 0 \quad \text{for } k \to
\infty
\]
for each trace-class operator $ R \in \B(H_{0,1}) $ satisfying $ \tr
(R) = 0 $; here $ \A_k = \A_{E_k} \subset \A_{0,1} $.
\end{definition}

Clearly, such a definition can be used in a much more general
situation, for an arbitrary sequence of von Neumann subalgebras of $
\B(H) $. A decreasing sequence ($ E_1 \supset E_2 \supset \dots $,
thus $ \A_1 \supset \A_2 \supset \dots $) is infinitesimal if and only
if the intersection of all $ \A_k $ is trivial (scalar operators
only); this fact will not be used, and I omit the proof. Rather, a
non-monotone sequence of such a form will be used:
\[
E_n = \bigcup_{k=0}^{n-1} \bigg( \frac1n \Big( k + \frac{1-\eps_n}2
\Big), \, \frac1n \Big( k + \frac{1+\eps_n}2 \Big) \bigg) \, ,
\]
the union of $ n $ equidistant intervals, of total length $ \mes E_n =
\eps_n $. The class of sequences $ (\eps_n) $ such that $ (E_n) $ is
infinitesimal, is an invariant of a \cpHs. Continuum of type III
product systems will be presented by means of that invariant.

\medskip

Let $ (H_{s,t})_{s<t} $ be a \cpHs\ over $ [0,1] $, and subspaces $
H'_{s,t} \subset H_{s,t} $, $ \dim H'_{s,t} > 0 $, satisfy $ H'_{r,s}
\ot H'_{s,t} = H'_{r,t} $ whenever $ 0 \le r < s < t \le 1 $; then $
(H'_{s,t})_{s<t} $ is another \cpHs\ on $ [0,1] $, embedded into $
(H_{s,t})_{s<t} $.\index{pair of embedded \cpHss}\index{embedded pair
of \cpHss}
An invariant for such a pair $ \( (H_{s,t})_{s<t},
(H'_{s,t})_{s<t} \) $ is proposed below. Of course, it is an invariant
of a pair (embedding), not of a single \cpHs. However, in some
situations (especially relevant to type $ II $ product systems) it
leads to some invariants of a single \cpHs. It happens, when subspaces $
H'_{s,t} $ can be unambiguously defined in terms of the given \cpHs\ $
(H_{s,t})_{s<t} $. Especially, we may ask about \emph{one-dimensional}
$ H'_{s,t} $.

\begin{definition}\label{2.2d}
Let $ (H_{s,t})_{s<t} $ be a \cpHs\ over $ [0,1] $.

(a) A \emph{decomposable vector}\index{decomposable!vector}
(of the \cpHs) is a vector $ \psi \in
H_{0,1} $, $ \psi \ne 0 $, such that for every $ t \in (0,1) $, $ \psi
$ is of the form $ \psi' \ot \psi'' $ for some $ \psi' \in H_{0,t} $,
$ \psi'' \in H_{t,1} $.

(b) The \cpHs\ is \emph{of type $ I $,}\index{type of \cpHs}
if decomposable vectors span the
whole $ H_{0,1} $ (as a closed linear space).

(c) The \cpHs\ is \emph{of type $ III $,}\index{type of \cpHs!3@$ III $}
if it has no decomposable vectors.

(d) The \cpHs\ is \emph{of type $ II $,} if it does not belong to
types $ I $, $ III $.

(e) The \cpHs\ is \emph{of type $ II_0 $,}\index{type of \cpHs!2@$ II_0 $}
if it is of type $ II $, and all decomposable vectors belong to a single
one-dimensional subspace of $ H_{0,1} $.
\end{definition}

For a \cpHs\ of type $ II_0 $ (examples will be given in Sect.~6),
one-dimensional subspaces $ H'_{s,t} \subset H_{s,t} $ satisfying $
H'_{r,s} \ot H'_{s,t} = H'_{r,t} $ are unique; in such a case, an
invariant of a pair leads to an invariant of a \cpHs\ itself. The
invariant is described below for an arbitrary pair, but the reader may
restrict himself to the case of $ \dim H'_{s,t} = 1 $ and type $ II_0
$.

Given an elementary set \eqref{12}, we introduce the subspace
\[
(H/H')_E = H_E \ot H'_{[0,1]\setminus E} = H'_{0,s_1} \ot H_{s_1,t_1}
\ot H'_{t_1,s_2} \ot \dots \ot H_{s_n,t_n} \ot H'_{t_n,1} \subset
H_{0,1} \, ,\index{He@$ (H/H')_E $, subspace}
\]
and the corresponding 
projection\index{Qe@$ Q_E $, projection}
\[
Q_E = \One_E \otimes P'_{[0,1]\setminus E} = P'_{0,s_1} \ot
\One_{s_1,t_1} \ot P'_{t_1,s_2} \ot \dots \ot \One_{s_n,t_n} \ot
P'_{t_n,1} \in \A_{0,1} \, ,
\]
where $ P'_{s,t} \in \A_{s,t} $ is the projection onto $ H'_{s,t}
$. Clearly,
\begin{gather*}
(H/H')_{E_1\cap E_2} = (H/H')_{E_1} \cap (H/H')_{E_2} \, , \\
Q_{E_1\cap E_2} = Q_{E_1} Q_{E_2} = Q_{E_2} Q_{E_1}
\end{gather*}
for all $ E_1, E_2 \in \cE $. (Note however that in general $
Q_{E_1\cup E_2} $ is strictly larger than $ Q_{E_1} + Q_{E_2} $ for $
E_1 \cap E_2 = \emptyset $.) Thus, $ \{ Q_E : E \in \cE \} $ is a
commuting set of projections, and $ H_{0,1} $ decomposes into the
corresponding direct integral. We'll describe the relevant measure
space explicitly, using the following quite general result of measure
theory.

\begin{lemma}\label{2.2}
Let $ X $ be a compact topological space, $ \A $ an algebra\footnote{%
 That is, $ \A $ is closed under complement and finite union.}
of subsets of $ X $, and $ \mu : \A \to [0,\infty) $ an additive
function satisfying the condition

\textup{(a) [regularity]}\index{regular!additive function}
for every $ A \in \A $ and $ \eps > 0 $ there exists
$ B \in \A $ such that $ \overline B \subset A $ \textup{(here $
\overline B $ is the closure of $ B $)} and $ \mu(B) \ge \mu(A) - \eps
$.

Then $ \mu $ has a unique extension to a measure on the \sif\
generated by $ \A $.
\end{lemma}

\begin{proof}
Due to a well-known theorem, it is enough to prove that $ \mu $ is $
\sigma $-additive on $ \A $. Let $ A_1 \supset A_2 \supset \dots $, $
A_1, A_2, \dots \in \A $, $ \cap A_k = \emptyset $; we have to prove
that $ \mu (A_k) \to 0 $. Given $ \eps > 0 $, we can choose $ B_k \in
\A $ such that $ \overline B_k \subset A_k $ and $ \mu (B_k) \ge \mu
(A_k) -
2^{-k} \eps $. Due to compactness, the relation $ \cap \overline B_k
\subset \cap A_k = \emptyset $ implies $ \overline B_1 \cap \dots \cap
\overline B_n = \emptyset $ for some $ n $. Thus, $ \mu (A_n) = \mu (
A_1 \cap \dots \cap A_n ) \le \mu ( B_1 \cap \dots \cap B_n ) + \mu (
A_1 \setminus B_1 ) + \dots + \mu ( A_n \setminus B_n ) < \eps $.
\end{proof}

\begin{remark}\label{2.3}
All $ A \in \A $ such that $ A $ and $ X \setminus A $ both satisfy
the regularity condition, are a subalgebra of $ \A $. (The proof is
left to the reader.) Therefore it is enough to check the condition for
$ A $ and $ X \setminus A $ where $ A $ runs over a set that generates
the algebra $ \A $.
\end{remark}

\begin{lemma}\label{2.4}
Let $ X $ be a compact topological space, $ \A $ an algebra of subsets
of $ X $, $ H $ a Hilbert space, and $ Q : \A \to \B(H) $ be such that
operators $ Q(A) $ are mutually commuting Hermitian projections, $ Q (
A \cup B ) = Q(A) + Q(B) $ whenever $ A \cap B = \emptyset $, $ A,B
\in \A $. Assume that for each $ x \in H $ the function $ \A \ni A
\mapsto \ip{Q(A)x}{x} \in [0,\infty) $ satisfies the regularity
condition \textup{\ref{2.2}(a).} Then $ Q $ has a unique extension to a
projection-valued
measure\index{projection-valued measure, extension to}
on the \sif\ generated by $ \A $.
\end{lemma}

\begin{proof}
We choose $ x \in H $ such that $ \ip{Px}{x} > 0 $ for every non-zero
projection $ P $ belonging to the set $ \ccP $ of all Hermitian
projections of the commutative von Neumann algebra generated by $ \{
Q(A) : A \in \A \} $. We apply Lemma \ref{2.2} to the function $ \mu :
\A \to [0,\infty) $, $ \mu(A) = \ip{Q(A)x}{x} $. Let $ \A $ be
equipped with the metric $ \dist (A,B) = \mu ( A \bigtriangleup B ) $,
and $ \ccP $ --- with the metric $ \dist (P,Q) = \ip{|P-Q|x}{x} $, then
$ Q : \A \to \ccP $ is isometric, and can be extended by continuity to
the \sif\ generated by $ \A $.
\end{proof}

We return to $ (H_{s,t})_{s<t} $ and $ Q_E $. The increasing
projection-valued function $ t \mapsto Q_{(0,t)} $ has at most a
countable set of discontinuity points (since $ H_{0,1} $ is
separable), as well as the decreasing $ t \mapsto Q_{(t,1)} $. We
denote by
$ D $\index{D@$ D $, countable set of discontinuity points}
the union of these two sets, and restrict ourselves to
elementary sets \eqref{12} such that $ s_1, t_1, \dots, s_n, t_n
\notin D $ --- these will be called
\emph{regular}\index{regular!elementary set}
elementary sets.\footnote{%
 In the homogeneous setup, the set $ D $ is necessarily empty; all
 elementary sets are regular. The reader may restrict himself to such
 a case.}
Note that continuity (at $ \eps = 0+ $) of $ Q_{(0,t-\eps)} = (
\One_{0,t-\eps} \ot P'_{t-\eps,t} ) \ot P'_{t,1} $ implies continuity
of $ \One_{0,t-\eps} \ot P'_{t-\eps,t} $; similarly,  continuity of $
Q_{(t+\eps,1)} = P'_{0,t} \ot ( P'_{t,t+\eps} \ot \One_{t+\eps,1} ) $
implies continuity of $ P'_{t,t+\eps} \ot \One_{t+\eps,1} $. Together,
they imply continuity of $ ( \One_{0,t-\eps} \ot P'_{t-\eps,t} ) \ot (
P'_{t,t+\eps} \ot \One_{t+\eps,1} ) = Q_{(0,t-\eps)\cup(t+\eps,1)} $
provided that $ t \in (0,1) \setminus D $. It follows easily that $
Q_E $ treated as a function of $ 2n $ variables $ s_1, t_1, \dots,
s_n, t_n $ (according to \eqref{12}) is continuous in these variables
at every point such that $ s_1, t_1, \dots, s_n, t_n \in (0,1)
\setminus D $.

We introduce topological space 
$ \cC $\index{C@$ \cC $, the space of closed sets}
of all closed sets $ C \subset
[0,1] $, equipped with the Vietoris 
topology,\index{Vietoris topology}
that is, the Hausdorff metric\index{Hausdorff metric}
\[
\dist ( C_1, C_2 ) = \max_{t\in[0,1]} \bigg| \min_{s\in C_1} |s-t| -
\min_{s\in C_2} |s-t| \bigg|
\]
under the convention that $ \min_{s\in\emptyset} |s-t| = 1 $ for all $
t $. It is well-known that $ \cC $ is compact.

We denote by $ \A_0 $ the algebra of subsets of $ \cC $ generated by
subsets of the form $ \{ C \in \cC : C \cap I \ne \emptyset \} $,
where $ I $ runs over intervals $ I \subset [0,1] $ of all kinds (say,
$ (s,t) $, and $ [s,t) $, and $ [t,t] = \{ t \} $, and $ \emptyset $,
etc). Every partition $ [0,1] = I_1 \cup \dots \cup I_n $ of $ [0,1] $
into a finite number $ n $ of intervals $ I_1, \dots, I_n $ leads to a
partition of $ \cC $ into $ 2^n $ sets of $ \A_0 $ (just choose for
each $ k = 1, \dots, n $ one of two possibilities, either $ C \cap I_k
= \emptyset $ or $ C \cap I_k \ne \emptyset $), that is, a finite
subalgebra of $ \A_0 $ containing $ 2^{2^n} $ sets.

A finer partition of $ [0,1] $ leads to a finer partition of $ \cC $
and a larger finite subalgebra of $ \A_0 $. The whole $ \A_0 $ is the
union of these finite subalgebras (over all partitions of $ [0,1] $).

However, we restrict ourselves to `regular' intervals, whose endpoints
belong to $ (0,1) \setminus D $ (though, intervals $ [0,t) $, $ [0,t]
$, $ (t,1] $ and $ [t,1] $ are regular for $ t \in (0,1) \setminus D
$), and the corresponding `regular' subalgebra $ \A \subset \A_0 $,
the union of $ 2^{2^n} $-subalgebras over all `regular' partitions.

We may define a projection-valued additive function $ \mu : \A \to
\B(H_{0,1}) $ by
\begin{equation}\label{13}
\mu ( \{ C : C \cap I = \emptyset \} ) = Q_{[0,1]\setminus I}
\end{equation}
for any `regular' interval $ I \subset [0,1] $. Indeed, each of the $
2^n $ elements of the partition of $ \cC $ that correspond to $ [0,1]
= I_1 \cup \dots \cup I_n $ is the intersection of $ n $
sets of the form $ \{ C : C \cap I_k = \emptyset \} $ or $ \{ C : C
\cap I_k \ne \emptyset \} $, and the corresponding projection is the
product of $ n $ projections of the form $ Q_{[0,1]\setminus I_k} $ or
$ \One - Q_{[0,1]\setminus I_k} $.

It follows from \eqref{13} (using multiplicativity) that
\begin{equation}\label{14}
\mu ( \{ C : C \subset E \} ) = Q_E
\end{equation}
for any `regular' elementary set $ E $.

\begin{lemma}
The additive function $ \mu $ has a unique extension to a projection
measure\index{projection-valued measure, extension to}
on the Borel \sif\ of $ \cC $.
\end{lemma}

\begin{proof}
The $ 2^n $ elements of the partition of $ \cC $ that correspond to $
[0,1] = I_1 \cup \dots \cup I_n $ are sets of diameter $ \le \max (
|I_1|, \dots, |I_n| ) $. It follows easily that the \sif\ generated by
$ \A $ is the whole Borel \sif\ of $ \cC $. Due to Lemma \ref{2.4}, it
is enough to check the regularity condition for $ \ip{\mu(A)x}{x}
$. Due to Remark \ref{2.3} we may do it only for $ A = \{ C : C \cap
I = \emptyset \} $ and its complement. It remains to note that $
Q_{[0,s]\cup[t,1]} $ is continuous in $ s,t $ whenever $ s,t \in (0,1)
\setminus D $; the same for $ Q_{[0,s]} $ and $ Q_{[t,1]} $.
\end{proof}

The projection-valued Borel measure $ \mu $ on $ \cC $ that satisfies
\eqref{13}, \eqref{14}, will be called the \emph{spectral
measure}\index{spectral measure on $ \cC $, projection-valued!of an
 embedded pair}\index{embedded pair of \cpHss!spectral measure on $ \cC $}
of the pair $ (H_{s,t})_{s<t} $, $ (H'_{s,t})_{s<t} $.

For every $ x \in H_{0,1} $, the positive Borel measure $ A \mapsto
\ip{ \mu(A)x }{ x } $ will be called the spectral measure of $ x $.

The $ \sigma $-ideal of Borel sets $ A_0 \subset \cC $ such that $
\mu(A_0) = 0 $, or equivalently, the $ \sigma $-filter of Borel sets $
A_1 \subset \cC $ such that $ \mu(\cC \setminus A_1) = 0 $, is an
invariant of the pair $ (H_{s,t})_{s<t} $, $ (H'_{s,t})_{s<t} $. It is
convenient to express the relation $ \mu(\cC \setminus A_1) = 0 $ by
saying that `almost all spectral sets belong to $ A_1 $'.

Note that
\begin{equation}\label{15}
\mu \( \{ C : C \ni t \} \) = 0 \quad \text{for each } t \in (0,1)
\setminus D \, ;
\end{equation}
in other words, almost all spectral sets do not contain a given point
$ t $ (unless $ t $ belongs to at most countable set $ D $). It
follows (via Fubini's theorem) that almost all spectral sets are of zero
Lebesgue measure, and therefore, nowhere dense.

\medskip

Till now, Hilbert spaces over $ \C $ were 
considered.\index{real versus complex}\index{complex versus real}
However, some
constructions are more natural over $ \R $. Definition \ref{1.1} works
equally well over $ \R $, giving `real \cpHs'. Every Hilbert space $ H
$ over $ \R $ (`real Hilbert space') has its
complexification\index{complexification!of real Hilbert space}
$ H^\C
$, a Hilbert space over $ \C $ (`complex Hilbert space') containing $
H $ as an $ \R $-linear subspace. Every $ A \in \B(H) $ has its
complexification $ A^\C \in \B(H^\C) $; we have $ \| A^\C \| = \| A \|
$, $ (A^*)^\C = (A^\C)^* $, $ (AB)^\C = A^\C B^\C $. On the other
hand, the general form of an element $ A \in \B(H^\C) $ is $ A = (\Re
A)^\C + i (\Im A)^\C $; $ \Re A, \Im A \in \B(H) $.\footnote{%
 Here $ \Re A $ is not $ (A+A^*)/2 $. Rather, it may be thought of as
 the \emph{component-wise} real part of a matrix.}
Note that $ \max (
\| \Re A \|, \| \Im A \| ) \le \| A \| \le \| \Re A \| + \| \Im A \|
$. Every real unitary (that is,
invertible linear isometric) operator $ U : H_1 \to H_2 $ between real
Hilbert spaces has its complexification $ U^\C : H_1^\C \to H_2^\C $,
a complex unitary operator. Compositions are respected: $ (U_2 U_1)^\C
= U_2^\C U_1^\C $. Tensor products are also respected: $ (H_1 \ot
H_2)^\C = H_1^\C \ot H_2^\C $ (up to a canonical isomorphism), and $
(A\ot B)^\C = A^\C \ot B^\C $. Thus, every real \cpHs\ $ \(
(H_{s,t})_{s<t}, (W_{r,s,t})_{r<s<t} \) $ has its
complexification\index{complexification!of real \cpHs}
$ \( (H^\C_{s,t})_{s<t}, (W^\C_{r,s,t})_{r<s<t} \) $, a complex
\cpHs. Isomorphic real \cpHss\ have 
isomorphic\index{isomorphism!and complexification}
complexifications.

Can it happen that two nonisomorphic real \cpHss\ have isomorphic
complexifications?\index{isomorphism!and complexification}
I do not know.

However, invariants introduced in this section are robust to
complexification, as stated below.

Definition \ref{2.1} works equally well over $ \R $.

\begin{lemma}\label{2.7}
A sequence of elementary sets is infinitesimal for a real \cpHs\ if and
only if it is 
infinitesimal\index{infinitesimal sequence of elementary sets!and
 complexification}
for the complexified \cpHs. \textup{(The
time set is $ [0,1] $, as before.)}
\end{lemma}

\begin{proof}
\begin{sloppypar}
We have $ H_{0,1} = H_{E_k} \ot H_{(0,1)\setminus E_k} $, $ \A_k =
\B(H_{E_k}) \ot \One_{(0,1)\setminus E_k} $, and $ H^\C_{0,1} =
H^\C_{E_k} \ot H^\C_{(0,1)\setminus E_k} $, $ \A^\C_k = \B(H^\C_{E_k})
\ot \One_{(0,1)\setminus E_k} $. The general form of an element $ A
\in \A^\C_k $ is $ A = (\Re A)^\C + i (\Im A)^\C $; $ \Re A, \Im A \in
\A_k $; $ \max ( \| \Re A \|, \| \Im A \| ) \le \| A \| \le \| \Re A
\| + \| \Im A \| $. Similarly, the general form of a trace-class
operator $ R \in \B(H^\C_{0,1}) $ is $ (\Re R)^\C + i (\Im R)^\C $,
where $ \Re R, \Im R \in \B(H_{0,1}) $ are trace-class operators, and $
\tr (R) = \tr(\Re R) + i\tr(\Im R) $. Thus, $ AR = \( (\Re A)^\C + i
(\Im A)^\C \) \( (\Re R)^\C + i (\Im R)^\C \) = \( (\Re A) (\Re R)
\)^\C - \( (\Im A) (\Im R) \)^\C + i \( (\Re A) (\Im R) \)^\C + i \(
(\Im A) (\Re R) \)^\C  $; $ \tr \( \( (\Re A) (\Re R) \)^\C \) = \tr
\( (\Re A) (\Re R) \) $, etc.; $ | \tr(AR) | \le | \tr \( (\Re A) (\Re
R) \) | + \text{\small (three more terms)} $;
\[
\sup_{A\in\A^\C_k, \|A\|\le1} | \tr(AR) | \le 2 \sup_{A\in\A_k,
\|A\|\le1} | \tr(A \Re R) | + 2 \sup_{A\in\A_k, \|A\|\le1} | \tr(A \Im
R) | \, .
\]
If the sequence is infinitesimal for the real \cpHs, then the
right-hand side tends to $ 0 $, and therefore the left-hand side tends
to $ 0 $, which means that the sequence is infinitesimal for the
complexified \cpHs.
\end{sloppypar}

On the other hand, if $ R \in \B(H_{0,1}) $ is a  trace-class
operator, then $ \tr(AR) = \tr( A^\C R^\C ) $ for $ A \in \A_k $, and
\[
\sup_{A\in\A_k, \|A\|\le1} | \tr(AR) | \le \sup_{A\in\A^\C_k,
\|A\|\le1} | \tr(A R^\C) | \, ,
\]
which gives us the second implication.
\end{proof}

We turn to a pair, a \cpHs\ $ (H_{s,t})_{s<t} $ on $ [0,1] $ and an
embedded \cpHs\ $ (H'_{s,t})_{s<t} $, $ H'_{s,t} \subset H_{s,t} $,
assuming this time, that $ H_{s,t}, H'_{s,t} $ are real (rather than
complex). Projections $ Q_E $ are still well-defined, and the spectral
measure $ \mu $ satisfying \eqref{13}, \eqref{14} emerges as
before. Alternatively we may consider the complex embedding $
H^{\prime\C}_{s,t} \subset H^\C_{s,t} $. Corresponding projections are
just $ Q_E^\C $ (since $ \( (H/H')_E \)^\C = (H^\C/H^{\prime\C})_E
$). It follows that the spectral measure for the complex embedding is
the complexification of $ \mu $, that is, $ A \mapsto \( \mu(A) \)^\C
$ (here $ \( \mu(A) \)^\C $ is the complexification of the projection
$ \mu(A) $). Indeed, \eqref{14} shows that the equality between the
two projection-valued Borel measures holds on sets of the form $ A =
\{ C : C \subset E \} $; therefore it holds for all sets. The
conclusion follows.

\begin{lemma}
The complexification $ A \mapsto \( \mu(A) \)^\C $ of the spectral
measure $ \mu $ of an embedding $ H'_{s,t} \subset H_{s,t} $ of real
\cpHss\ is equal to the spectral 
measure\index{spectral measure on $ \cC $, projection-valued!and
 complexification}
of the corresponding embedding
$ H^{\prime\C}_{s,t} \subset H^\C_{s,t} $ of complex \cpHss.
\end{lemma}

All said in sections 3--7 holds for real and complex Hilbert spaces
equally well.\index{real versus complex}\index{complex versus real}

\section{Continuous products of measure classes}

`Square roots of measures', introduced by Accardi \cite{Ac}, are
quite useful for constructing product systems, which was suggested by
A.~Vershik. For definitions and basic facts see
\cite[Sect.~3.3.5]{Ar}, \cite[Sect.~14.4]{Ar-b}. The following
definition is somewhat more restrictive than Arveson's, since I
restrict myself to measure classes generated by a single measure.

\begin{definition}
(a) A \emph{measure class}\index{measure class}
is a triple $ (X,\B,\M) $ consisting of a
standard\footnote{%
 Another restriction in addition to Arveson's.}
Borel space $ (X,\B) $ and a nonempty set $ \M $ of finite positive
measures on $ (X,\B) $ such that for some (therefore, every) $ \mu \in
\M $ and for every finite positive measure $ \nu $ on $ (X,\B) $,
\[
\nu \sim \mu \quad \text{if and only if} \quad \nu \in \M \, ,
\]
$ \nu \sim \mu $ denoting mutual absolute continuity. (Such $ \M $
will be called an equivalence class of
measures.)\index{equivalence class of measures}

(b) An \emph{isomorphism}\index{isomorphism!of measure classes}
of a measure class $ (X,\B,\M) $ to another
measure class $ (X',\B',\M') $ is a map $ \eta : X \to X' $ which is a
$ \bmod \, 0 $ isomorphism of measure spaces $ (X,\B,\mu) $ and $
(X',\B',\mu') $ for some $ \mu \in \M $, $ \mu' \in \M' $.
\end{definition}

As usual, all negligible\index{negligible set}
sets (that is, subsets of Borel sets of
measure $ 0 $ for some, therefore every, $ \mu \in \M $) may be added
to the \sif. Instead of requiring $ (X,\B) $ to be a standard Borel
space, we may require for some (therefore every) $ \mu \in \M $, that
$ (X,\B,\mu) $ is a Lebesgue-Rokhlin measure space; it means a measure
space, isomorphic $ \bmod\,0 $ to some interval $ (0,a) $ with
Lebesgue measure, or a discrete (finite or countable) measure space,
or a combination of both. Up to $ \bmod\,0 $ isomorphism, the two
approaches are equivalent.

\begin{quote}
Throughout, either by assumption or by construction, all measure
spaces will be
Lebesgue-Rokhlin\index{measure space: always Lebesgue-Rokhlin}
spaces;\footnote{%
 In particular, all negligible sets are measurable.}
also, all claims and constructions will be invariant under $ \bmod\,0
$ isomorphisms.
\end{quote}

Hilbert spaces $ L_2(\mu) $, $ L_2(\nu) $ for $ \mu,\nu \in \M $ may
be glued together via the unitary operator
\[
L_2(\mu) \ni f \mapsto \sqrt{ \frac\mu\nu } f \in L_2(\nu) \, ;
\]
here $ \frac\mu\nu $ is the Radon-Nikodym derivative (denoted also by $
\frac{d\mu}{d\nu} $). These spaces may be treated as `incarnations' of
a single Hilbert space $ L_2 (X,\B,\M)
$.\footnote{%
 More formally: an element of $ L_2 (X,\B,\M) $ is, by definition, a
 family $ \psi = (\psi_\mu)_{\mu\in\M} $ such that $ \psi_\mu \in L_2
 (X,\B,\mu) $ for every $ \mu \in \M $, and $ \psi_\nu = \sqrt{
 \frac\mu\nu } \psi_\mu $ for all $ \mu,\nu \in \M $.}%
\index{L2@$ L_2 (X,\B,\M) $, Hilbert space out of a measure class}
The general form of an
element of $ L_2 (X,\B,\M) $ is $ f \sqrt\mu $, where $ \mu \in \M $
and $ f \in L_2(X,\B,\mu) $, taking into account the relation
\[
f \sqrt\mu = \bigg( \sqrt{ \frac\mu\nu } f \bigg) \sqrt\nu \, .
\]

Any isomorphism of measure classes induces naturally a unitary
operator between the corresponding Hilbert spaces.

The product of two measure classes is defined naturally, and
\[
L_2 \( (X,\B,\M) \times (X',\B',\M') \) = L_2 (X,\B,\M) \ot L_2
(X',\B',\M') \, ;
\]
that is, we have a canonical unitary operator between these spaces,
namely, $ f \sqrt\mu \ot f' \sqrt{\mu'} \mapsto (f\ot f') \sqrt{ \mu
\ot \mu' } $, where $ (f\ot f') (\om,\om') = f(\om) f'(\om') $. I'll
write in short $ L_2 (\M) $ instead of $ L_2 (X,\B,\M) $; thus, $
L_2 ( \M \times \M' ) = L_2 ( \M ) \ot L_2 ( \M' ) $. Everyone knows
the similar fact for measure spaces, $ L_2 ( \mu \times \mu' ) = L_2 (
\mu ) \ot L_2 ( \mu' ) $.

Given isomorphisms of measure classes $ \eta : \M \to \M_1 $ (that is,
$ \eta : (X,\B,\M) \to (X_1,\B_1,\M_1) $) and $ \eta' : \M' \to \M'_1
$, the map $ \eta \times \eta' $ defined by $ ( \eta \times \eta' )
(\om,\om') = \( \eta(\om), \eta'(\om') \) $, is an isomorphism of
measure classes, $ \eta \times \eta' : \M \times \M' \to \M_1 \times
\M'_1 $.

\begin{definition}\label{3.2}
(a) A continuous product of measure classes (`\cpmc', for 
short)\index{CP MC@\cpmc\ = continuous product of measure classes}
is a family $ \( (X_{s,t}, \B_{s,t}, \M_{s,t} ) \)_{-\infty\le
s<t\le\infty} $ of measure classes such that
\[
(X_{r,t}, \B_{r,t}, \M_{r,t} ) = (X_{r,s}, \B_{r,s}, \M_{r,s} ) \times
(X_{s,t}, \B_{s,t}, \M_{s,t} )
\]
whenever $ -\infty\le r<s<t\le\infty $, and $ \M_{-\infty,\infty} $ is
non-degenerate (that is, not just the zero measure).

\begin{sloppypar}
(b) An
isomorphism\index{isomorphism!of \cpmcs}
of one \cpmc\ $ \( (X_{s,t}, \B_{s,t}, \M_{s,t} )
\)_{-\infty\le s<t\le\infty} $ to another \cpmc\ $ \( (X'_{s,t},
\B'_{s,t}, \M'_{s,t} ) \)_{-\infty\le s<t\le\infty} $ is a family $
(\eta_{s,t})_{-\infty\le s<t\le\infty} $ where each $ \eta_{s,t} $ is
an isomorphism of measure classes, $ (X_{s,t}, \B_{s,t}, \M_{s,t} ) $
to $ (X'_{s,t}, \B'_{s,t}, \M'_{s,t} ) $, and
\[
\eta_{r,t} = \eta_{r,s} \times \eta_{s,t}
\]
whenever $ -\infty\le r<s<t\le\infty $.
\end{sloppypar}
\end{definition}

I'll write simply $ (\M_{s,t})_{s<t} $ instead of $ \( (X_{s,t},
\B_{s,t}, \M_{s,t} ) \)_{s<t} $. Given a \cpmc\ $ (\M_{s,t})_{s<t} $, we
may construct the corresponding \cpHs\ $ (H_{s,t})_{s<t} $; just $
H_{s,t} = L_2 (\M_{s,t}) $.

Does every \cpHs\ (up to isomorphism) emerge from some \cpmc? I do not
know.

All product systems constructed in this paper emerge (via \cpHss) from
\cpmcs.

Another useful language for describing \cpmcs\ is based on sub-\sif
s. By a sub-\sif\ on a (Lebesgue-Rokhlin!) measure space $ (X,\B,\mu)
$, or a measure class space $ (X,\B,\M) $, I always mean a sub-\sif\
of $ \B $ that contains all negligible
sets.\index{sub-\sif, always completed}
The sub-\sif\
generated by the union $ \F_1 \cup \F_2 $ of sub-\sif s $ \F_1, \F_2 $
will be denoted by $ \F_1 \vee \F_2 $.

\begin{definition}
Sub-\sif s $ \F_1, \dots, \F_n $ on a measure class $ (X,\B,\M) $ are
\emph{independent,}\footnote{%
 They could also be called quasiindependent; see footnotes
 \ref{foo18}, \ref{foo24}.}%
\index{independent!sub-\sif s on a measure class}
if there exists a probability measure $ \mu \in \M
$ such that
\begin{equation}\label{16}
\mu ( A_1 \cap \dots \cap A_n ) = \mu(A_1) \dots \mu(A_n) \quad
\text{for all } A_1 \in \F_1, \dots, A_n \in \F_n \, .
\end{equation}
\end{definition}

For independent $ \F_1, \dots, \F_n $ the sub-\sif\ $ \F_1 \vee \dots
\vee \F_n $ generated by them will be denoted also by $ \F_1 \ot \dots
\ot \F_n $.

Given a product of two measure classes,
\[
(X,\B,\M) = (X_1,\B_1,\M_1) \times (X_2,\B_2,\M_2) \, ,
\]
we have two independent sub-\sif s $ \F_1, \F_2 $ such that $ \B =
\F_1 \ot \F_2 $; roughly,
\[
\F_1 = \{ A \times X_2 : A \in \B_1 \} \, , \quad
\F_2 = \{ X_1 \times B : B \in \B_2 \} \, ,
\]
though, all negligible sets must be added.

And conversely, every two independent sub-\sif s $ \F_1, \F_2 $ such
that $ \B = \F_1 \ot \F_2 $ emerge from a representation of $
(X,\B,\M) $ (up to isomorphism) as a product; in fact, $
(X_k,\B_k,\M_k) = (X,\B,\M)/\F_k $ is the quotient space.

\begin{remark}\label{3.4}
If $ \F_1, \F_2, \F_3 $ are sub-\sif s on a measure class $ (X,\B,\M)
$ such that
\begin{gather*}
\F_1 \text{ and } (\F_2\vee\F_3) \text{ are independent},\\
\F_2 \text{ and } \F_3 \text{ are independent},
\end{gather*}
then
\[
\F_1, \F_2, \F_3 \text{ are independent}.
\]
Indeed, we can split off the first factor, and then split the
rest. (It may become a difficulty, if you do not use the technique of
products.)
\end{remark}

\begin{sloppypar}
Given a \cpmc\ $ \( (X_{s,t}, \B_{s,t}, \M_{s,t} ) \)_{-\infty\le
s<t\le\infty} $, we have a family of sub-\sif s $
(\F_{s,t})_{-\infty\le s<t\le\infty} $ on the measure class $
(X,\B,\M) =
(X_{-\infty,\infty},\B_{-\infty,\infty},\M_{-\infty,\infty}) $
satisfying $ \F_{-\infty,\infty} = \B $ and
\begin{equation}\label{17}
\F_{r,t} = \F_{r,s} \ot \F_{s,t} \quad \text{whenever } -\infty\le
s<t\le\infty \, .
\end{equation}
And conversely, every such family of sub-\sif s on a measure class
emerges from a \cpmc.
\end{sloppypar}

Condition \eqref{17} is evidently equivalent to the combination of two
conditions,
\begin{gather}
\F_{r,t} = \F_{r,s} \vee \F_{s,t} \quad \text{whenever } -\infty\le
 r<s<t\le\infty \, , \label{18} \\
\F_{-\infty,s} \text{ and } \F_{s,\infty} \text{ are independent for
 all } s \in \R \, . \label{19}
\end{gather}

\begin{remark}\label{3.5a}
It is tempting to restrict ourselves to one-parameter families $
(\F_{-\infty,s})_{s\in\R} $ and $ (\F_{s,\infty})_{s\in\R} $, since $
\F_{s,t} = \F_{-\infty,s} \cap \F_{s,\infty} $. However, there is an
obstacle, illustrated by the following elementary example. Let
\begin{gather*}
X = \{-1,+1\}^3 = \{ (\tau_1,\tau_2,\tau_3) : \tau_1,\tau_2,\tau_3 =
 \pm1 \} \, , \displaybreak[0] \\
\B = \text{(all subsets of $ X $)} \, , \displaybreak[0] \\
\mu = \text{(the uniform probability distribution on $ X $)} \, ,
 \displaybreak[0] \\
\F_1 \text{ is generated by } \tau_1 \, , \displaybreak[0] \\
\F_{1,2} \text{ is generated by } \tau_1 \text{ and } (1-\tau_1)
 \tau_2 + (1+\tau_1) \tau_2 \tau_3 \, , \displaybreak[0] \\
\F_{2,3} \text{ is generated by } \tau_2 \text{ and } \tau_3 \, ,
\displaybreak[0] \\
\F_3 \text{ is generated by } \tau_3 \, .
\end{gather*}
Then $ \F_1 \subset \F_{1,2} $, $ \F_{2,3} \supset \F_3 $, and $ \F_1
\ot \F_{2,3} = \B = \F_{1,2} \ot \F_3 $. Nevertheless, the sub-\sif\ $
\F_2 = \F_{1,2} \cap \F_{2,3} $ is trivial, and so, $ \F_1 \vee \F_2
\vee \F_3 = \F_1 \vee \F_3 \ne \B $.
\end{remark}

The following definition is equivalent to Definition \ref{3.2}.

\begin{definition}\label{3.5}
\begin{sloppypar}
(a) A \cpmc\index{CP MC@\cpmc\ = continuous product of measure classes}
is $ \( (X,\B,\M), (\F_{s,t})_{-\infty\le s<t\le\infty} \)
$ where $ (X,\B,\M) $ is a measure class and $ \F_{s,t} $ are sub-\sif s
satisfying \eqref{18}, \eqref{19} and $ \F_{-\infty,\infty} = \B $.
\end{sloppypar}

(b) An
isomorphism\index{isomorphism!of \cpmcs}
of one \cpmc\ $ \( (X,\B,\M), (\F_{s,t})_{-\infty\le
s<t\le\infty} \) $ to another \cpmc\ $ \( (X',\B',\M'),
(\F'_{s,t})_{-\infty\le  s<t\le\infty} \) $ is an
isomorphism of measure classes, $ (X,\B,\M) $ to $ (X',\B',\M') $,
that sends each $ \F_{s,t} $ to $ \F'_{s,t} $.
\end{definition}

\begin{definition}
\begin{sloppypar}
A \emph{homogeneous continuous product of measure
classes}\index{HCP MC@\hcpmc\ = homogeneous continuous product of
 measure classes}\index{homogeneous!CPMC@\cpmc}
(`\hcpmc', for
short) is $ \( (X,\B,\M), (\F_{s,t})_{-\infty\le s<t\le\infty},
(\ti\theta^u)_{u\in\R} \) $, where $ \( (X,\B,\M), (\F_{s,t})_{s<t} \)
$ is a \cpmc, and $ (\ti\theta^u)_{u\in\R} \) $ is a one-parameter
group of automorphisms $ \ti\theta^u $ of the measure class $
(X,\B,\M) $ such that
\[
\ti\theta^u (\F_{s,t}) = \F_{s+u,t+u} \quad \text{whenever } -\infty \le
s < t \le \infty, u \in \R \, , 
\]
and the corresponding unitary group $ (\theta^u)_{u\in\R} $ on $ L_2
(X,\B,\M) $ is strongly continuous.
\end{sloppypar}
\end{definition}

I leave to the reader elaborating on the latter condition, and
defining isomorphisms of \hcpmcs.

Clearly, every \hcpmc\ leads to an \hcpHs, and (via Theorem \ref{1.5})
to a product system.\index{product system!out of \hcpmc}

\section{Continuous products of probability spaces}

By a probability space\index{probability space}
I always mean a (Lebesgue-Rokhlin!) measure
space $ (X,\B,\mu) $ such that $ \mu(X) = 1 $. By an
isomorphism\index{isomorphism!of probability spaces}
of probability spaces I mean a $ \bmod \, 0 $ isomorphism (measure
preserving). The product of two probability spaces is a probability
space, of course.

A continuous product of probability spaces (`\cpps', for
short)\index{CP PS@\cpps\ = continuous product of probability spaces}
is defined similarly to Def.~\ref{3.2}(a); that is,
\[
( X_{r,t}, \B_{r,t}, \mu_{r,t} ) = ( X_{r,s}, \B_{r,s}, \mu_{r,s} )
\times ( X_{s,t}, \B_{s,t}, \mu_{s,t} ) \, ,
\]
where $ \mu_{s,t} $ are probability measures. Definition \ref{3.2}(b)
is also adapted readily.

\begin{sloppypar}
An equivalent definition, similar to Definition~\ref{3.5}: a 
\cpps\index{CP PS@\cpps\ = continuous product of probability spaces}
is $ \( (X,\B,\mu), \linebreak[0]
(\F_{s,t})_{-\infty\le s<t\le\infty} \) $, where $
(X,\B,\mu) $ is a probability space, and $ \F_{s,t} $ are sub-\sif s
satisfying \eqref{17} (or equivalently, \eqref{18} and \eqref{19});
this time,
independence\index{independent!sub-\sif s on a probability space}
means \eqref{16} for the given $ \mu $ (rather
than some equivalent measure). Definition \ref{3.5}(b) is also adapted
readily.
\end{sloppypar}

Given a \cpps\ $ (\mu_{s,t})_{s<t} $, we may construct the corresponding
\cpHs\ $ (H_{s,t})_{s<t} $; just $ H_{s,t} = L_2 (\mu_{s,t})
$. Alternatively, we may construct first the corresponding \cpmc\ $
(\M_{s,t})_{s<t} $, where $ \M_{s,t} $ consists of all measures
equivalent to $ \mu_{s,t} $; and second, out of $ (\M_{s,t})_{s<t} $
we construct $ (H_{s,t})_{s<t} $, $ H_{s,t} = L_2 (\M_{s,t}) $. Both
ways lead to the same (up to isomorphism) \cpHs\ $ (H_{s,t})_{s<t}
$. Thus, the class of \cpHss\ corresponding to \cppss\ is a part of the
class of \cpHss\ corresponding to \cpmcs. It is a proper
subclass. Indeed, any \cpHs\ of the form $ (L_2(\mu_{s,t}))_{s<t} $ has
a decomposable vector, namely, a constant function. (In fact, such a
\cpHs\ can be of type $ I $ or $ II $.) In contrast, a \cpHs\ of the
form $ (L_2(\M_{s,t}))_{s<t} $ need not have a decomposable
vector. (In fact, such a \cpHs\ can be of type $ I $, $ II $ or $ III
$.)

The lemma below establishes a relation between \cpmcs, \cppss\ and
decomposable vectors. Note that for every measure class $ (X,\B,\M) $
and every vector $ \psi \in L_2(\M) $ we have the corresponding
positive measure $ |\psi|^2 $ on $ (X,\B) $, satisfying $ |\psi|^2 (X)
= \|\psi\|^2 $.\index{zz@$ "|\psi"|^2 $, measure}
Namely, if $ \psi = f \sqrt\mu $ then $ |\psi|^2 =
|f|^2 \mu $ (that is, $ \frac{ d|\psi|^2 }{ d\mu } = |f|^2 $); the
result does not depend on $ \mu \in \M $. Also, for $ \psi_1 \in
L_2(\M_1) $, $ \psi_2 \in L_2(\M_2) $ we have $ | \psi_1 \ot \psi_2
|^2 = | \psi_1 |^2 \times | \psi_2 |^2 $ (the product measure).

\begin{lemma}\label{4.1}
Let $ (\M_{s,t})_{s<t} $ be a \cpmc, and $ \psi \in L_2 (\M_{0,1}) $ a
decomposable vector, $ \|\psi\| = 1 $; $ \psi = \psi_{0,1} $, $
\psi_{r,t} = \psi_{r,s} \ot \psi_{s,t} $, $ \psi_{s,t} \in H_{s,t} =
L_2 (\M_{s,t}) $, $ \| \psi_{s,t} \| = 1 $. Then measures $ \mu_{s,t}
= |\psi_{s,t}|^2 $ form a \cpps\ $ (\mu_{s,t})_{s<t} $.
\end{lemma}

\begin{proof}
$ \mu_{r,t} = |\psi_{r,t}|^2 = | \psi_{r,s} \ot \psi_{s,t} |^2 =
|\psi_{r,s}|^2 \times |\psi_{s,t}|^2 = \mu_{r,s} \times \mu_{s,t} $.
\end{proof}

Every \cpps\ $ (\mu_{s,t})_{s<t} $, be it of type $ I $ or $ II $, can
arise as $ \mu_{s,t} = |\psi_{s,t}|^2 $. Indeed, we may take the
corresponding \cpmc\ $ (\M_{s,t})_{s<t} $ and the decomposable vector
$ \psi_{s,t} = \sqrt{ \mu_{s,t} } $.

Homogeneous \cppss\ are defined evidently.
\hcppss, also called `noises', are studied in probability theory, see
a survey \cite{Ts+}. Classical noises (white, Poisson and their
combinations) lead to \hcpHss\ (and product systems) of type $ I
$. Nonclassical noises, found by J.~Warren and S.~Watanabe, lead to
\hcpHss\ (and product systems) of type $ II_n $, $ n > 0 $. Some of
them have spectral sets of any given Hausdorff dimension between $ 0 $
and $ 1/2 $. Some are time-asymmetric. Black noises, found by
A.~Vershik, myself and S.~Watanabe, lead to \hcpHss\ (and product
systems) of type $ II_0 $. See \cite{Ts+} for references and
details. Thus, a part of the results of the present work can be
derived from the theory of \hcppss. However, \hcpmcs\ give simpler
proofs and stronger results about
\hcpHss.\index{product system!out of \hcpps}

\section{Random sets, and type $ II_0 $}

Let $ \cC $\index{C@$ \cC $, the space of closed sets}
be the set of all closed sets $ C \subset \R $ (not just $
C \subset [0,1] $ as in Sect.~2). We equip $ \cC $ with the Effros
\sif\ $ \B $;\index{Effros \sif}
it is generated by functions $ C \mapsto \min_{s\in C}
|s-t| $ indexed by $ t \in \R $,\footnote{%
 Once again, $ \min_{s\in\emptyset} |s-t| = 1 $ for all $ t $, by
 definition.}
and it turns $ \cC $ into a standard Borel space.

For any $ s,t $ such that $ -\infty\le s<t\le\infty $ we define the
sub-\sif\ $ \B_{s,t} \subset \B $ generated by the map $ C \mapsto C
\cap [s,t] $, that is, by functions $ C \mapsto \min_{r\in C\cap[s,t]}
|r-u| $, $ u \in \R $. Clearly, these sub-\sif s satisfy
\eqref{18}. Till now, no measure appears, and no sets are
negligible. The next definition introduces an equivalence class $ \M $
of measures on $ \cC $ (thus turning the Borel space $ (\cC,\B) $ into
a measure class $ (\cC,\B,\M) $), after which all \sif s must be
completed with negligible sets.

\begin{definition}\label{5.1}
An equivalence class $ \M $ of measures on $ (\cC,\B) $ is
\emph{decomposable,}\index{decomposable!equivalence class of measures
 on $ \cC $}
if

(a) sub-\sif s $ \B_{-\infty,s} $ and $ \B_{s,\infty} $ on the measure
class $ (\cC,\B,\M) $ are independent for all $ s \in \R $;

(b) the set $ \{ C \in \cC : C \ni t \} $ is negligible, for each $ t
\in \R $.
\end{definition}

Condition (a) is just \eqref{19}; we see that every decomposable
equivalence class of measures $ \M $ leads to a \cpmc\ $ \( (\cC,\B,\M), 
(\B_{s,t})_{-\infty\le s<t\le\infty} \) $. The set $ \{ C \in \cC : C
\ni t \} $, belonging to the intersection $ \B_{-\infty,t} \cap
\B_{t,\infty} $  of independent sub-\sif s, must be either negligible
or co-negligible (that is, of negligible complement); the second
possibility is of no interest for us, and we exclude it by (b).

Looking forward to a relation between $ \M $ and spectral measures, I
find it convenient to say `almost all spectral sets belong to $ A_1 $'
instead of `$ \cC \setminus A_1 $ is negligible (w.r.t.\ $ \M
$)'. Similarly to Sect.~2 (recall \eqref{15}), almost all spectral
sets do not contain a given point $ t $ (due to \ref{5.1}(b)),
therefore, almost all spectral sets are nowhere dense, of Lebesgue
measure $ 0 $.

\begin{lemma}\label{5.2}
Let $ \M $ be a decomposable eqvivalence class of measures on $ \cC $,
and $ -\infty < s < t < \infty $. Then the set $ \{ C \in \cC : C \cap
[s,t] = \emptyset \} $ is not negligible.
\end{lemma}

\begin{proof}
Due to \ref{5.1}(b), each $ r $ has a neighborhood $ (r-\eps,r+\eps) $
such that the set $ \{ C \in \cC : C \cap (r-\eps,r+\eps) = \emptyset
\} $ is not negligible. Due to compactness of $ [s,t] $ we may choose
a single $ \eps $ for all $ r \in [s,t] $. Non-negligible sets $ \{ C
\in \cC : C \cap [s,s+\eps] = \emptyset \} $ and $ \{ C \in \cC : C
\cap [s+\eps,s+2\eps] = \emptyset \} $ belong to independent \sif s $
\B_{-\infty,s+\eps} $ and $ \B_{s+\eps,\infty} $ respectively.
Therefore their intersection $ \{ C \in \cC : C \cap [s,s+2\eps] =
\emptyset \} $ is not negligible. Continuing this way, we see that the
set $ \{ C \in \cC : C \cap [s,s+n\eps] \cap [s,t] = \emptyset \} $ is
not negligible for any $ n $.
\end{proof}

The quotient measure class $ (\cC,\B,\M) / \B_{s,t} $ may be thought
of as the measure class $ (\cC_{s,t},\ti\B_{s,t},\M_{s,t}) $, where $
\cC_{s,t} $ is the space of all closed subsets of $ [s,t] $, $
\ti\B_{s,t} $ its Effros \sif, and $ \M_{s,t} $ consists of images of
measures $ \mu \in \M $ under the map $ C \mapsto C \cap [s,t]
$. Lemma \ref{5.2} shows that the empty set (or rather, $ \{ \emptyset
\} $) is an atom (non-negligible point) of $
(\cC_{s,t},\ti\B_{s,t},\M_{s,t}) $. Consider the
probability measure $ \mu_{s,t} $ concentrated at the
atom. Clearly, $
\mu_{r,s} \times \mu_{s,t} = \mu_{r,t} $. In terms of the \cpHs\ $
(H_{s,t})_{s<t} = \( L_2(\M_{s,t}) \)_{s<t} $ it means that vectors $
\psi_{s,t} = \sqrt{\mu_{s,t}} $ satisfy
\[
\psi_{r,s} \ot \psi_{s,t} = \psi_{r,t} \quad \text{whenever } -\infty
< r < s < t < \infty \, .
\]
According to Definition \ref{2.2d}(a), $ \psi_{0,1} $ is a
decomposable vector of the \cpHs \linebreak
$ (H_{s,t})_{0\le s<t\le1} $
over $ [0,1] $; let us call it `the standard decomposable
vector'.\index{standard decomposable vector}
We see
that each decomposable equivalence class $ \M $ of measures on $ \cC $
leads to a \cpHs\ of type $ I $ or $ II $.

The spectral measure $ \mu $ is defined (see Sect.~2) for the pair of
the \cpHs\ $ \(L_2(\M_{s,t})\)_{0\le s<t\le1} $ and its standard
decomposable vector; $ \mu $ is a projection-valued Borel measure on $
(\cC_{0,1}, \ti\B_{0,1}) $.

As every projection-valued measure, $ \mu $ has its equivalence class
$ \M_\mu $ of measures; it consists of positive Borel measures on $
\cC_{0,1} $ of the form $ A \mapsto \ip{ \mu(A)\psi }{ \psi } $, where
$ \psi $ runs over all vectors of $ L_2(\M_{0,1}) $ such that $ \ip{
\mu(A)\psi }{ \psi } > 0 $ whenever $ \mu(A) \ne 0 $.

\begin{lemma}
$ \M_\mu = \M $ for every decomposable eqvivalence class $ \M $ of
measures on $ \cC $.
\end{lemma}

\begin{proof}
For any elementary set $ E \subset (0,1) $, the subspace $ (H/H')_E =
Q_E H $ of $ H = H_{0,1} $ may be described as consisting of vectors $
\psi = f \sqrt\mu $ such that the measure $ |\psi|^2 = |f|^2 \cdot \mu
$ is concentrated on the set $ \{ C \in \cC_{0,1} : C \subset E \}
$. It means that $ Q_E $ is the operator of multiplication by the
indicator\index{indicator}\index{zz@$ \One $, indicator function of a set}\footnote{%
 The indicator $ \One_S $ of a set $ S $ is a function equal to $ 1 $
 on $ S $ and $ 0 $ on its complement.}
of the set. Thus,
for all $ \psi $, $ \ip{ Q_E \psi }{ \psi } = |\psi|^2 \( \{ C : C
\subset E \} \) $. Using \eqref{14} we see that $ \ip{ \mu(A) \psi }{
\psi } = |\psi|^2 (A) $
for $ A \subset \cC_{0,1} $ of the form $ \{ C : C \subset E \} $;
therefore it holds for all Borel sets $ A $. Choosing $ \psi $ such
that $ \ip{ \mu(A) \psi }{ \psi } > 0 $ whenever $ \mu(A) \ne 0 $, we
see that $ |\psi|^2 \in \M $. On the other hand, the measure $ A
\mapsto \ip{ \mu(A) \psi }{ \psi } $ belongs to $ \M_\mu $.
\end{proof}

\begin{remark}
The spectral measure of an arbitrary vector $ \psi \in H_{0,1} $ is
given, for every Borel set $ A \subset \cC_{0,1} $, by  $ \ip{ \mu(A)
\psi }{ \psi } = |\psi|^2 (A) $ (as is shown in the proof above).
\end{remark}

\begin{remark}
We see that every `random set' satisfying conditions \ref{5.1}(a,b)
can appear as the spectral set (of some pair of \cpHss). The converse
is also true, with some reservations. Namely, let $ \mu $ be the
spectral measure of a pair ($ H'_{s,t} \subset H_{s,t} $ for $ 0 \le s
< t \le 1 $), and $ \M_\mu $ the corresponding equivalence class of
measures on $ \cC_{0,1} $. Then $ \M_\mu $ satisfies condition
\ref{5.1}(b) for all $ t \in (0,1) \setminus D $, where $ D $ is an at
most countable set (introduced in Sect.~2). Condition \ref{5.1}(a) is
also satisfied for $ s \in (0,1) \setminus D $. This fact will not be
used, and I give only a hint toward the proof:
\[
\mu ( A' \times A'' ) = \mu' (A') \ot \mu'' (A'')
\]
for all Borel sets $ A' \subset \cC_{0,s} $, $ A'' \subset \cC_{s,1}
$; here $ \mu' $ is the spectral measure of the restriction to $ [0,s]
$ of the given pair of \cpHss, thus, each $ \mu' (A') $ is a projection
in $ H_{0,s} $. Similarly, each $ \mu'' (A'') $ is a projection in $
H_{s,1} $.
\end{remark}

A probability measure $ \mu $ on $ (\cC,\B) $ will be called
decomposable,\index{decomposable!probability measure on $ \cC $}
if

(a) sub-\sif s $ \B_{-\infty,s} $ and $ \B_{s,\infty} $ on the
probability space $ (\cC,\B,\mu) $ are independent for all $ s \in \R
$;

(b) [the same as \ref{5.1}(b)].

Such a measure $ \mu $ leads to a \cpps\ $ \( (\cC,\B,\mu), 
(\B_{s,t})_{-\infty\le s<t\le\infty} \) $.

\begin{lemma}
A probability measure $ \mu $ on $ (\cC,\B) $ is decomposable if and
only if it describes the Poisson point process corresponding to a
nonatomic locally finite measure $ \nu $ on $ \R $.
\end{lemma}

\begin{proof}
`If' is evident; `only if' will be proven. Consider the number of
points in a bounded interval, $ N = | C \cap [s,t] | $ (an integer or
infinity) as a random variable on $ (\cC,\B,\mu) $. By Lemma
\ref{5.2}, $ \Pr{ N=0 } > 0 $ (`$ \Pr{\dots} $' stands for probability). We
divide $ [s,t] $ into $ n $ equal intervals (`cells') and consider the
random number $ N_n $ of occupied cells. Clearly, $ N_n \uparrow N $
for $ n \to \infty $. Due to independence, $ \Pr{ N_n = 0 } =
(1-p_1^{(n)}) \dots (1-p_n^{(n)}) $, where $ p_k^{(n)} $ is the
probability of cell $ k $ being occupied. We know that $ \inf_n
(1-p_1^{(n)}) \dots (1-p_n^{(n)}) \ge \Pr{ N=0 } > 0 $; therefore $
\sup_n \Ex N_n = \sup_n \( p_1^{(n)} + \dots + p_n^{(n)} \) \le -\ln
\Pr{ N=0 } < \infty $. Also, $ \max_n \( p_1^{(n)}, \dots, p_n^{(n)}
\) \to 0 $ for $ n \to \infty $ (recall the compactness argument in
the proof of Lemma \ref{5.2}). By a well-known theorem of probability
theory, $ N $ has a Poisson distribution. It remains to use
independence.
\end{proof}

\begin{corollary}\label{5.6}
Every decomposable probability measure on $ (\cC,\B) $ is concentrated
on the set of locally finite sets $ C \in \cC $ \textup{(that is, $ C
$ has no accumulation points).}
\end{corollary}

A closed set $ C \in \cC $ is called
perfect,\index{perfect set}
if it has no isolated points. The only locally finite perfect set is
the empty set.

\begin{lemma}\label{5.8}
Let $ \M $ be a decomposable eqvivalence class of measures on $ \cC $
such that almost all sets $ C \in \cC $ are perfect. Then the
corresponding \cpHs\ over $ [0,1] $ is of type $ II_0 $
\textup{(}unless almost all sets $ C \in \cC $ are empty\textup{).}
\end{lemma}

\begin{proof}
Let $ \psi \in H_{0,1} $ be a decomposable vector. We have to prove
that it is, up to a coefficient, the standard decomposable vector;
that is, the measure $ |\psi|^2 $ is concentrated at the atom $ C =
\emptyset $. On one hand, the measure $ |\psi|^2 $, being absolutely
continuous w.r.t.\ $ \M_{0,1} $, is concentrated on (the set of all)
perfect sets.\footnote{%
 If $ C $ is a perfect set, then $ C \cap [0,1] $ is perfect, provided
 that $ 0 \notin C $, $ 1 \notin C $, which holds a.s.}
On the other hand, the measure $ |\psi|^2 $, being decomposable (by
Lemma \ref{4.1}), is concentrated on (the set of all) locally finite
sets (by Corollary \ref{5.6}, adapted to $ \cC_{0,1} $).
\end{proof}

\begin{corollary}\label{5.9}
\begin{sloppypar}
Let $ \M, \M' $ be decomposable equivalence classes of measures on $
\cC $, both such that almost all sets $ C \in \cC $ are
perfect. Assume that their images $ \M_{0,1}, \M'_{0,1} $ under the
map $ C \mapsto C \cap [0,1] $ are different equivalence classes of
measures on $ \cC_{0,1} $. Then the corresponding \cpHss\ $ \(
L_2(\M_{s,t}) \)_{0\le s<t\le1} $, $ \( L_2(\M'_{s,t}) \)_{0\le
s<t\le1} $ over $ [0,1] $ are non-isomorphic.
\end{sloppypar}
\end{corollary}

\section{Constructing random sets}

\begin{sloppypar}
The set of zeros of a Brownian motion is a good example of a random
set. More generally, we consider a stationary, sample-continuous,
Fellerian Markov process $ (\xi_t)_{-\infty<t<\infty} $ on a
metrizable compact topological space $ K $. In other words, we have a
strongly continuous semigroup $ (T_t)_{0\le t<\infty} $ of positive
unit-preserving operators $ T_t $ on
the space $ C(K) $ of continuous functions $ f : K \to \R $, the
conjugate semigroup $ (T^*_t)_{0\le t<\infty} $ on the space $ M(K) $
of measures on $ K $, and a probability measure $ \nu $ on the space $
C (\R,K) $ of continuous maps $ \xi : \R \to K $ (called `sample
paths') such that
\begin{multline}\label{193}
\int_{C (\R,K)} f_1 (\xi(t_1)) \dots f_n (\xi(t_n)) \, d\nu(\xi) =
 \\
= \int_{C (\R,K)} f_1 (\xi(t_1)) \dots f_{n-1} (\xi(t_{n-1})) \(
 T_{t_n-t_{n-1}} f_n \) (\xi(t_{n-1})) \, d\nu(\xi)
\end{multline}
whenever $ n \in \{2,3,\dots\} $, $ -\infty < t_1 < \dots < t_n <
\infty $ and $ f_1,\dots,f_n \in C(K) $. In
more probabilistic language, the conditional distribution of $
\xi(t_n) $ given $ \xi(t_1), \dots, \xi(t_{n-1}) $ is $
T^*_{t_n-t_{n-1}} \de_{\xi(t_{n-1})} $; here $ \de_a $ is the unit
mass at $ a $. Note that the distribution of $ \xi(0) $ (as well as
any $ \xi(t) $) is a probability measure $ \nu_0 $ on $ K $, invariant
in the sense that $ T^*_t \nu_0 = \nu_0 $ for all $ t $. The measure $
\nu $ is uniquely determined by $ \nu_0 $ and $ (T_t)_{0\le t<\infty}
$, and is stationary (that is, invariant under time
shifts).\footnote{%
 Existence of an invariant measure $ \nu_0 $ for a given semigroup
 follows easily from compactness; namely, $ \| (T^*_t-\One) \frac1s
 \int_0^s T^*_r \, dr \| \le \frac{2t}{s} \to 0 $ for $ s \to \infty
 $, therefore each weak$^*$-limit point (for $ s \to \infty $) of $
 \frac1s \int_0^s (T^*_r) (\nu_1) \, dr $ is an invariant measure,
 irrespective of the choice a probability measure $ \nu_1 $. However,
 we do not need it.}
\end{sloppypar}

Given a closed subset $ K_0 \subset K $ such that $ \nu_0 (K_0) = 0 $,
we consider the random set $ C = \{ t \in \R : \xi(t) \in K_0 \} $; it
is a measurable map from the probability space $ \( C(\R,K), \nu \)
$ to the Borel space $ \cC $ of closed subsets of $ \R $. The
distribution of the random set is a probability measure $ \mu $ on $
\cC $ that satisfies Condition \ref{5.1}(b), and is
stationary. However, Condition \ref{5.1}(a) needs special effort.

Assume that the random set is non-empty a.s., and consider the first
(after $ t=0 $) hit of $ K_0 $,
\[
\tau (\xi) = \min \{ t \in [0,\infty) : \xi(t) \in K_0 \} \, ;
\]
the conditional distribution of the pair $ (\tau, \xi(\tau)) $, given
$ \xi(0) $, is a measurable map from $ (K,\nu_0) $ to probability
measures on $ [0,\infty) \times K_0 $.

\newtheorem*{asme1}{The first condition of smearing}
\begin{asme1}\index{smearing condition}
There exists an equivalence class of measures on $ [0,\infty) \times
K_0 $ that contains almost all conditional distributions of $ (\tau,
\xi(\tau)) $ given $ \xi(0) $. (That is, for $ \nu_0 $-almost all $ a
\in K $ the conditional distribution of $ (\tau, \xi(\tau)) $ given $
\xi(0)=a $ belongs to the equivalence class.)
\end{asme1}

Note that $ \{0\} \times K_0 $ is negligible in $ [0,\infty) \times
K_0 $, since $ \xi(0) \notin K_0 $ a.s.

\begin{lemma}\label{6.1}
The first condition of smearing ensures that the random set satisfies
Condition \textup{\ref{5.1}(a).}
\end{lemma}

\begin{proof}
Due to stationarity, it is enough to consider $ \B_{-\infty,0} $ and $
\B_{0,\infty} $. The conditional
distribution of the future $ C \cap [0,\infty) $ of the random set $ C
$, given the past $ \xi|_{(-\infty,0]} $ of the process, is the
mixture of its conditional distributions given $ (\tau, \xi(\tau)) $,
over the conditional distribution of $ (\tau, \xi(\tau)) $ given $
\xi(0) $. Using the strong Markov property and the first condition of
smearing we see that $ \xi(0) $ influences the conditional
distribution but not its equivalence class.
\end{proof}

\newtheorem*{asme2}{The second condition of smearing}
\begin{asme2}\index{smearing condition}
For every $ t > 0 $, measures $ T^*_t \de_a $ (for all $ a \in K $)
are (pairwise) equivalent.
\end{asme2}

The first condition of smearing constrains $ \nu $ and $ K_0 $, while
the second condition of smearing constrains only $ \nu $, or rather,
the semigroup $ (T_t) $. Nevertheless:

\begin{lemma}\label{6.2}
The second condition of smearing implies the first condition of
smearing.
\end{lemma}

\begin{proof}
Assume the contrary: there exists a Borel set $ B \subset [0,\infty)
\times K_0 $ such that the conditional probability $ \cP{
(\tau,\xi(\tau)) \in B }{ \xi(0)=a } $ vanishes for some but not all $ a
\in K \setminus K_0 $. We may assume that $ B \subset (0,\infty)
\times K_0 $ and moreover, $ B \subset [t,\infty) \times K_0 $ for
some $ t > 0 $. By the Markov property,
\[
\cP{ (\tau,\xi(\tau)) \in B }{ \xi(0)=a } = \int_K \cP{
(\tau,\xi(\tau)) \in B }{ \xi(t)=b } \, d \( T^*_t \de_a \) (b) \, .
\]
The integral vanishes for some but not all $ a $, in contradiction to
equivalence of measures $ T^*_t \de_a $.
\end{proof}

\begin{figure}
\includegraphics{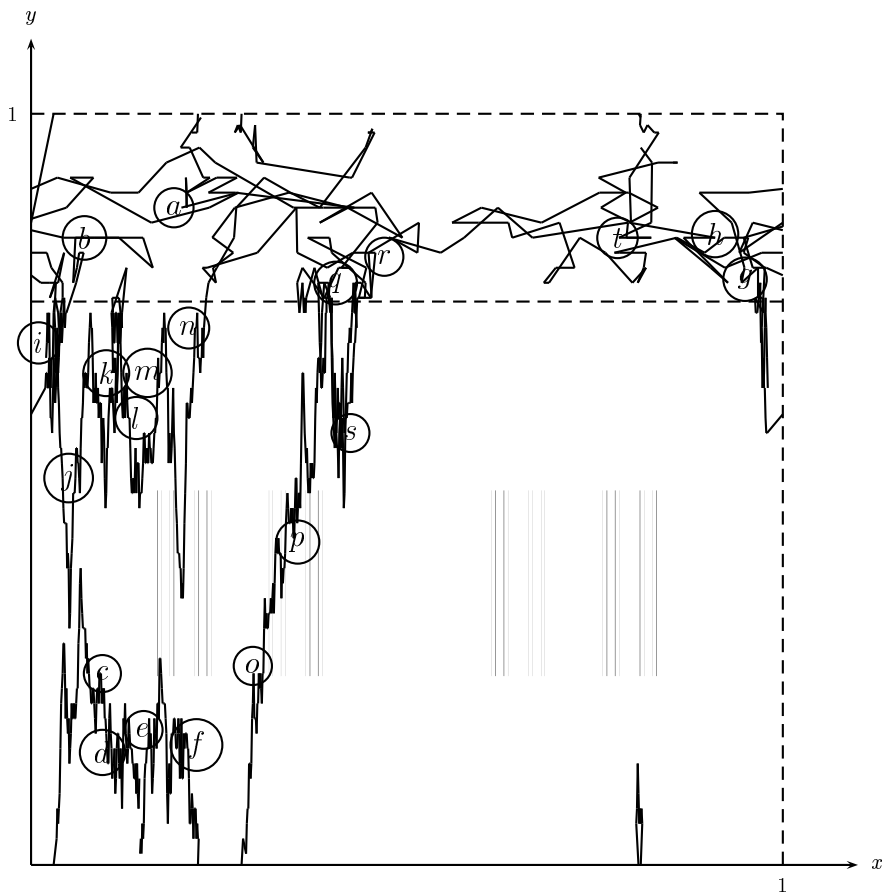}
\caption{A sample path of the diffusion process on the torus. Time
marks are shown by circled letters $ a,b,\dots,s,t $. The set $ K_0 $
is shown in gray.}
\label{fig1}
\end{figure}

Specifically, we construct a diffusion process $ \xi(\cdot) $ on the
two-dimensional torus $ K = \R^2/\Z^2 = \{ (x,y) : 0\le x\le1, 0\le
y\le1 \} $ with the identification: $ x=0 $ to $ x=1 $, and $ y=0 $ to
$ y=1 $. The process is governed by the differential operator
\begin{equation}\label{196}
L = \frac12 \frac{\pd^2}{\pd y^2} + f(y) \frac{\pd}{\pd x} + \frac12
g^2 (y) \frac{\pd^2}{\pd x^2} \, ,
\end{equation}
where $ f $ and $ g $ are $ C^\infty $ functions on the circle $ \R/\Z
$ such that
\begin{gather*}
f(y) > 0 \text{ and } f'(y) > 0 \text{ for } 0 \le y \le \frac34 \, ;
 \\
g(y) = 0 \text{ for } 0 \le y \le \frac34 \, ; \quad
g(y) > 0 \text{ for } \frac34 < y < 1 \, .
\end{gather*}
In terms of coordinates, $ \xi(t) = \( x(t), y(t) \) $, we see that $
y(\cdot) $ is just the standard Brownian motion (on the circle), and
the behavior of $ x(\cdot) $ depends on $ y(\cdot) $ dramatically:
\begin{equation}\label{198}
x'(t) = f(y(t)) \quad \text{when } 0 < y(t) < \frac34 \, .
\end{equation}
Here, $ x(\cdot) $ is continuously differentiable and monotone (on the
circle); in contrast, when $ \frac34 < y(t) < 1 $, $ x(\cdot) $ is as
irregular as the Brownian motion. See Fig.~\ref{fig1}.

The subset $ K_0 \subset K $ is
\[
K_0 = K_1 \times \Big[ \frac14, \frac12 \Big] \, ,
\]
where $ K_1 $ is a perfect subset of $ (0,1) $.

The operator $ L $ is not elliptic, since $ g^2(\cdot) $ is not strictly
positive. The degeneration impedes smearing, but gives us much better
control over the random set. Indeed, the process $ \xi(t) = \( x(t),
y(t) \) $ can hit $ K_0 $ only when $ y(t) \in [\frac14,\frac12] $;
here, the random set is determined by $ x(t) $ only, and $ x(t) $ is
an increasing, continuously differentiable function. Thus, the random
set is rather similar to the given set $ K_1 $.

On the other hand, it is essential that $ f $ is not constant (on the
relevant interval). Otherwise, if, say, $ f(\cdot) = 1 $, the (local)
metric structure of the set $ K_1 $ is reproduced by the random set
exactly, with no random deformation, which is inconsistent with
\ref{5.1}(a).

\begin{lemma}\label{6.3}
The semigroup generated by $ L $ satisfies the second condition of
smearing.
\end{lemma}

\begin{proof}
Let $ B $ be a Borel subset of the torus. We consider the function
\[
u(t,a) = \( T^*_t \de_a \) (B) \, .
\]
It is enough to prove for each $ t>0 $, that either $ u(t,a) = 0 $ for
all $ a \in K $, or $ u(t,a) > 0 $ for all $ a \in K $.

Assume for a while that $ u $ is a smooth function, of class $ C^2 $;
then $ u $ is a classical solution of the partial differential
equation $ \( \frac{\pd}{\pd t} - L \) u = 0 $. It is a second order
PDE with nonnegative characteristic form, and we may use the strong
maximum principle based on elliptic connectivity, see
\cite[Sect.~3.1]{OR}. Circles $ \{ (x,y,t) : y \in \R/\Z \} $,
parametrized by $ x \in \R/\Z $ and $ t>0 $, are lines of ellipticity
for the PDE $ \( \frac{\pd}{\pd t} - L \) u = 0 $ (see
\cite[p.~208]{OR} for the definition), due to the term $ \pd^2 / 2 \pd
y^2 $ of $ L $. Also, circles $ \{ (x,y,t) : x \in \R/\Z \} $ for $ y
\in \(\frac34,1\) $ and $ t>0 $ are lines of ellipticity, due to the
term $ g^2(y) \pd^2 / 2 \pd x^2 $. Thus, for each $ t>0 $, the torus $ K
\times \{t\} = \{ (x,y,t) : x,y \in \R/\Z \} $ is a set of elliptic
connectivity (see \cite[p.~209]{OR}). By \cite[Th.~3.1.2]{OR}, on such
a set either $ u=0 $ everywhere, or $ u>0 $ everywhere.

It remains to prove that $ u $ is of class $ C^2 $. (Now it is
essential that $ f(\cdot) \ne \const $.) The function $ u $ is a weak
solution of the PDE $ \( \frac{\pd}{\pd t} - L \) u = 0 $; it is
enough to check hypoellipticity of the operator $ \frac{\pd}{\pd t} -
L $. According to the well-known sufficient condition of H\"ormander
(see \cite[Th.~2.5.2]{OR}), it is enough to check that the system
\[
X_1 = \frac{\pd}{\pd y} \, , \quad X_2 = g(y) \frac{\pd}{\pd x} \, ,
\quad X_0 = f(y) \frac{\pd}{\pd x} + \frac{\pd}{\pd t}
\]
of first order differential operators (these represent $ L $ as $
\frac12 X_1^2 + \frac12 X_2^2 + X_0 $) has the full rank, equal to the
dimension $ 3 $ of $ K \times (0,\infty) $; see \cite[p.~167]{OR}. In
the domain $ \frac34 < y < 1 $ the full rank is evident, since vectors
\[
\begin{pmatrix} 0 \\ 1 \\ 0 \end{pmatrix} \, , \quad
\begin{pmatrix} g(y) \\ 0 \\ 0 \end{pmatrix} \, , \quad \begin{pmatrix}
f(y) \\ 0 \\ 1 \end{pmatrix} 
\]
are linearly independent. In the other domain, $ 0 \le y \le \frac34
$, $ X_2 $ vanishes, and we have only two operators. However, their
commutator
\[
[X_1,X_0] = \frac{\pd}{\pd y} \bigg( f(y) \frac{\pd}{\pd x} +
\frac{\pd}{\pd t} \bigg) - \bigg( f(y) \frac{\pd}{\pd x} +
\frac{\pd}{\pd t} \bigg) \frac{\pd}{\pd y} = f'(y) \frac{\pd}{\pd x}
\]
gives us the needed third vector.
\end{proof}

We have $ L \One = 0 $, but also $ L^* \One = 0 $ (since $ L^* $
differs from $ L $ only in the sign of $ f $),
which means that the Lebesgue measure $ \nu_0 $ on the torus $ K $ is
invariant under $ (T^*_t) $.\footnote{%
 Using hypoellipticity and the strong maximum principle, it is easy to
 see that the invariant measure is unique. However, we do not need
 it.}
Choosing once and for all functions $ f $ and $ g $, we hold fixed the
(distribution $ \nu $ of the) diffusion process $ \xi(\cdot) $. The
only parameter of our random set is the perfect set $ K_1 \subset
(0,1) $ of Lebesgue measure $ 0 $. Every such $ K_1 $ gives us a
measure $ \mu $ on $ \cC $, and its equivalence class $ \M $ is
decomposable; indeed, it satisfies Conditions \ref{5.1}(a,b).

Recall that $ K_0 = K_1 \times [\frac14,\frac12] $ is the union of
vertical segments. The set $ K_1 \times \{\frac14,\frac12\} $ of their
endpoints is polar, that is, never visited by the process $ \xi $,
under appropriate conditions on $ K_1 $, as will be shown by
constructing a `barrier function'.

\begin{lemma}\label{6.4}
There exists a function $ h : K \setminus \{ (\frac12,\frac12) \} \to
[0,\infty) $ such that

\textup{(a)} $ h $ is twice continuously differentiable on $ K \setminus \{
(\frac12,\frac12) \} $;

\textup{(b)} $ \displaystyle \sup_{ (x,y) \in K \setminus \{ (\frac12,\frac12)
\} } \sqrt{ |x-\tfrac12| } \, | (Lh) (x,y) | < \infty $;

\textup{(c)} $ \displaystyle \liminf_{x\to\frac12-} \sqrt{\tfrac12-x} \,
h(x,\tfrac12) > 0 $.
\end{lemma}

\begin{proof}
The function $ h $ will be constructed on a (punctured) neighborhood
of the point $ (\frac12,\frac12) $, assuming that $ f(\frac12) = 1 $
(just to simplify calculations; the general case follows by rescaling
$ x $ or $ y $). For convenience, we turn from $ x $ and $ y $ to new
coordinates $ z = y - \frac12 $, $ t = \frac12 - x $;\footnote{%
 The variable $ t $, local to the proof, has nothing in common with,
 say, $ t $ of $ \( \frac{\pd}{\pd t} - L \) $.}
now,
\[
L = \frac12 \frac{\pd^2}{\pd z^2} - \( 1 + \al z + O(z^2) \)
\frac{\pd}{\pd t} \, ,
\]
where $ \al = f'(1/2) $. We define
\[
h(z,t) = \frac1{\sqrt t} \exp \bigg( \! -\frac{z^2}{2t} + a \frac{z^3}t +
az \bigg) \quad \text{for } t > 0
\]
and $ h(z,t) = 0 $ for $ t \le 0 $; parameter $ a \in \R $ will be
chosen later.

Condition (c) becomes
\[
\liminf_{x\to0+} \sqrt t \, h(0,t) > 0 \, ,
\]
which clearly holds. Condition (b) becomes
\[
Lh (z,t) = O (1/\sqrt t) \quad \text{for } t > 0 \, ;
\]
this one remains to be checked (for some $ a $).

We have
\[
\frac{\pd^2}{\pd z^2} h(z,t) = \frac1{\sqrt t} \exp(\dots) \cdot
\bigg( \Big( -\frac z t + 3a \frac{z^2}t + a \Big)^2 + \Big( -\frac1t
+ 6a \frac z t \Big) \bigg) \, ;
\]
substituting $ z = r \sqrt t $,
\begin{multline*}
\frac{\pd^2}{\pd z^2} h(z,t) = \frac1{\sqrt t} \exp \bigg(
 -\frac{r^2}2 + ar^3 \sqrt t + ar\sqrt t \bigg) \cdot \bigg( \Big(
 -\frac r{\sqrt t} + 3a r^2 + a \Big)^2 - \frac1t + 6a \frac r{\sqrt
 t} \bigg) \\
= \frac1{\sqrt t} \exp(\dots) \cdot \bigg( \frac{r^2}t - 6a
 \frac{r^3}{\sqrt t} - 2a \frac r{\sqrt t} + O(r^4+1) - \frac1t + 6a
 \frac r{\sqrt t} \bigg) = \\
= \frac1{\sqrt t} \exp(\dots) \cdot \bigg( \frac{r^2-1}t + \frac{ 4ar
 - 6ar^3 }{ \sqrt t } + O(r^4+1) \bigg) \, .
\end{multline*}
Similarly,
\begin{multline*}
\frac{\pd}{\pd t} h(z,t) = \frac1{\sqrt t} \exp(\dots) \cdot \bigg(
 \! -\frac1{2t} + \frac{z^2}{2t^2} - a \frac{z^3}{t^2} \bigg) = \\
= \frac1{2\sqrt t} \exp(\dots) \cdot \bigg( \frac{r^2-1}t - 2a
\frac{r^3}{\sqrt t} \bigg) \, .
\end{multline*}
Taking into account that $ r\sqrt t = z = O(1) $ we get
\begin{multline*}
L h(z,t) = \frac1{2\sqrt t} \exp(\dots) \times \\
\times \bigg( \frac{r^2-1}t +
 \frac{ 4ar - 6ar^3 }{\sqrt t} + O(r^4+1) - \( 1 + \al r \sqrt t +
 O(r^2t) \) \Big( \frac{r^2-1}t - 2a \frac{r^3}{\sqrt t} \Big) \bigg) =
 \\
= \frac1{2\sqrt t} \exp(\dots) \cdot \bigg( \frac{ 4ar - 6ar^3 + 2ar^3
 - \al r(r^2-1) }{\sqrt t} + O(r^4+1) \bigg) = \\
= \frac1{2\sqrt t} \exp(\dots) \cdot \bigg( \! -\frac{ (4a+\al)r(r^2-1)
 }{\sqrt t} + O(r^4+1) \bigg) \, .
\end{multline*}
We choose the parameter
\[
a = - \frac\al4
\]
and get
\[
L h(z,t) = \frac1{2\sqrt t} \exp \bigg( \! -\frac{r^2}2 + a r^3 \sqrt t +
a r \sqrt t \bigg) \cdot O(r^4+1) \, .
\]
For $ |z| $ small enough (namely, $ |z| \le \frac1{4|a|} $),
\[
| L h(z,t) | \le \frac1{2\sqrt t} \exp \bigg( \! -\frac{r^2}4 \bigg)
\cdot O(r^4+1) = O \Big( \frac1{\sqrt t} \Big) \, .
\]
\end{proof}

From now on we restrict ourselves to sets $ K_1 $ such that
\begin{equation}\label{16a}
\Lambda_{1/2} (K_1) = 0 \, ,
\end{equation}
where $ \Lambda_{1/2} $ is the Hausdorff measure of dimension $ 1/2
$. The condition is satisfied by all sets $ K_1 $ of Hausdorff
dimension $ <1/2 $, and some sets of Hausdorff dimension $ =1/2 $.

\begin{lemma}\label{6.5}
Almost surely, the process $ \xi(\cdot) $ does not hit the set $ K_1
\times \{ \frac14, \frac12 \} $.
\end{lemma}

\begin{proof}
By \eqref{16a}, for every $ \eps $, the set $ K_1 $ can be covered by a
collection of intervals $ [x_k,x_k+\de_k] $ such that $ \sum_k
\sqrt{\de_k} \le \eps $. Consider such a function $ u_\eps : K \to
[0,\infty] $:
\[
u_\eps (x,y) = \sum_k \sqrt{\de_k} \, h ( x-x_k-2\de_k+\tfrac12, y )
\, ,
\]
where $ h $ is given by Lemma {6.4}. Condition \ref{6.4}(b) implies
\[
\int_0^1 \sup_y | L h(x,y) | \, dx \le \const \cdot \int
\frac{dx}{\sqrt{|x-\frac12|}} = C < \infty \, ,
\]
and therefore
\[
\int_0^1 \sup_y | L u_\eps (x,y) | \, dx \le C \sum_k \sqrt{\de_k} \le
C \eps \, .
\]
Condition \ref{6.4}(c) gives us $ c > 0 $ such that
\[
u_\eps (x,\tfrac12) > c \quad \text{for each } x \in K_1 \, .
\]
Also,
\[
| u_\eps (x,y) | \le M_y \sum_k \sqrt{\de_k} \le M_y \eps
\]
where $ M_y = \sup_x h(x,y) < \infty $ for $ y \ne \frac12 $.
Note that $ c, C, M_y $ do not depend on $ \eps $. We do so for each $
\eps $ of a sequence $ \eps_1, \eps_2, \dots $ such that $ \sum \eps_n
< \infty $ and consider
\[
u(x,y) = \sum_n u_{\eps_n} (x,y) \, .
\]
Then
\begin{gather}
\int_0^1 \sup_y | L u (x,y) | \, dx \le C \sum_n \eps_n < \infty \, ,
 \label{21} \\
u \Big( x, \frac12 \Big) = \infty \quad \text{for each } x \in K_1 \,
 , \label{22} \\
| u(x,y) | \le M_y \sum_n \eps_n < \infty \quad \text{for each } y \ne
 \frac12 \, . \label{23}
\end{gather}
The random process $ u(\xi(\cdot)) $ is a semimartingale (as long as
it is finite), whose drift
is $ (Lu) (\xi(\cdot)) $. When inside the region $ 0 < y < \frac34 $,
the drift is locally integrable (in time) due to \eqref{21},
\eqref{198}. Assume
that a sample path $ \xi(\cdot) $ hits $ K_1 \times \{\frac12\} $ at
some time $ t $ the first time. Then the martingale part of the
semimartingale at $ t- $ is bounded from below but unbounded from
above (due to \eqref{22}), which can happen only with probability $ 0
$ (by an evident generalization of \cite[Prop.\ V.1.8]{RY}). The same
holds for $ K_1 \times \{ \frac14 \} $.
\end{proof}

Given a sample path $ \xi(\cdot) $ of our diffusion process, we may
consider two closed sets on the time axis; one of them is our random
set $ C = \{ t : \xi(t) \in K_1 \times [\frac14,\frac12] \} $, the
other is $ D = \{ t : \xi(t) \in \R/\Z \times \{\frac14,\frac12\} \}
$. In terms of coordinates $ x(t), y(t) $ of the point $ \xi(t) $, $ D
= \{ t : y(t) = \frac14 \} \cup \{ t : y(t) = \frac12 \} $. By Lemma
6.5, $ C \cap D = \emptyset $. Both sets are unbounded. We get two
alternating sequences, $ \dots < s_{-1} < t_{-1} < s_0 < t_0 < s_1 <
t_1 < s_2 < \dots $ such that $ s_k \in C $, $ t_k \in D $, $ (s_k,t_k)
\cap D = \emptyset $, $ (t_k,s_{k+1}) \cap C = \emptyset $ for all $ k
\in \Z
$. Thus, $ C \subset \cup_k [s_k,t_k] $, and for all $ t \in \cup_k
[s_k,t_k] $ we have $ \frac14 \le y(t) \le \frac12 $; therefore $
x'(t) = f(y(t)) $ --- recall \eqref{198}. A
portion\index{portion of the random set}
$ C \cap [s_k,t_k] $ of the set $ C $,
equal to $ \{ t \in [s_k,t_k] : x(t) \in K_1 \} $, results from the
corresponding portion of $ K_1 $ by a smooth (of class $ C^1 $)
monotone transformation $ x(\cdot)|_{[s_k,t_k]} $. Of course, `the
corresponding portion of $ K_1 $' should be understood cyclically; $
K_1 \subset \R/\Z $, and $ x(\cdot) $ can make several revolutions
during $ [s_k,t_k] $. Note also that $ f(\frac14) \le x'(t) \le
f(\frac12) $ for $ t \in \cup_k [s_k,t_k] $. The numbers $ f(\frac14),
f(\frac12) $ are arbitrary, except for the constraint $ 0 < f(\frac14)
< f(\frac12) $. For every $ \eps > 0 $ the function $ f $ can be
chosen so that $ f(\frac14) = 1 $ and $ f(\frac12) = 1 + \eps $,
ensuring $ 1 \le x'(t) \le 1+\eps $ for $ t \in \cup_k [s_k,t_k] $.

The Hausdorff dimension of $ C \cap [s_k,t_k] $ is evidently equal to
the Hausdorff dimension of the corresponding portion of $ K_1 $. If
the local Hausdorff dimension of $ K_1 $ is equal to a given $ \al $
at every point of $ K_1 $, then the same holds for $ C $. Such a $ K_1
$ exists for each $ \al \in (0,\frac12) $; it leads to a decomposable
equivalence class of measures on $ \cC $ (Lemmas \ref{6.1}, \ref{6.2},
\ref{6.3}), a \cpmc\ (Sect.~5), and a \cpHs\ (Sect.~3) of type $ II_0
$ (Lemma \ref{5.8}). They all are invariant under time shifts; thus,
the \cpHs\ is in fact an \hcpHs\, and leads to a product system
(Theorem \ref{1.5}) of type $ II_0 $ (that is, having one-dimensional
spaces $ D(t) $ of decomposable vectors, see \cite[Remark
3.3.3]{Ar}). These product systems for different $ \al $ are
non-isomorphic (Corollary \ref{5.9} and Lemma \ref{1.6}), which proves
the following result.

\begin{theorem}\label{6.6}
There is a continuum of mutually non-isomorphic product systems of
type $ II_0 $.
\end{theorem}

It answers in the affirmative two questions of Arveson
\cite[p.~12]{Ar96}: (1) Are there uncountably many non-isomorphic
product systems? (2) (see also \cite[p.~167]{Ar99}) Is there a product
system of type $ II_0 \, $?

\section{Time reversal}

Recalling Sect.~1, a history $ \( (e^{itX})_{t\in\R}, \B(H_-) \ot \One \)
$ leads to an $ E_0 $-semigroup $ \al^+ $ and its product system. The
time-reversed history $ \( (e^{-itX})_{t\in\R}, \One \ot \B(H_+) \) $
leads to another $ E_0 $-semigroup $ \al^- $ and another product
system. These two product systems are called opposite (to each other);
see \cite[Sect.~3.5, 3.6]{Ar}.

For any \cpHs\ $ \( (H_{s,t})_{-\infty\le s<t\le\infty},
(W_{r,s,t})_{-\infty\le r<s<t\le\infty} \) $, its time-reversed
\cpHs\index{time-reversed!CPHS@\cpHs}
is, by definition, $ \( (\ti H_{s,t})_{-\infty\le s<t\le\infty},
\linebreak[0]
(\ti W_{r,s,t})_{-\infty\le r<s<t\le\infty} \) $, where $ \ti H_{s,t}
= H_{-t,-s} $ and $ \ti W_{r,s,t} ( x \ot y ) = W_{-t,-s,-r} ( y \ot x
) $ for $ x \in \ti H_{r,s} $, $ y \in \ti H_{s,t} $. For an \hcpHs,
we add $ \ti\theta^u_{s,t} = \theta^{-u}_{-t,-s}
$.\index{time-reversed!HCPHS@\hcpHs}

Each tail trivial \hcpHs\ leads to a history (Corollary
\ref{1.4}). Each \hcpHs\ leads to a product system (Theorem
\ref{1.5}). The time-reversed \hcpHs\ leads to the time-reversed
history (under tail triviality), and the opposite product system
(irrespective of tail triviality).

For a \cpmc\ $ \( (X,\B,\M), (\F_{s,t})_{-\infty\le s<t\le\infty} \) $
(recall Definition \ref{3.5}), its time-reversed
\cpmc\index{time-reversed!CPMC@\cpmc}
is, by
definition, $ \( (X,\B,\M), (\F_{-t,-s})_{-\infty\le s<t\le\infty} \)
$. It leads to the time-reversed \cpHs. The same for the homogeneous
case (\hcpmc\index{time-reversed!HCPMC@\hcpmc}
and \hcpHs), and for \cppss\ (homogeneous or not).

For a decomposable equivalence class of measures on $ (\cC,\B) $
(recall Definition \ref{5.1}), its time-reversed equivalence class of
measures on $ (\cC,\B) $ is, by definition, its image under the map $
\cC \ni C \mapsto \{ t : (-t) \in C \} \in \cC $. It leads to the
time-reversed \cpHs. (Homogeneity is respected, too.)

For the setup of the beginning of Sect.~6 (recall $ K $, $ (T_t) $,
$ \nu $, $ \nu_0 $), the time-reversed random process is described by the
image $ \ti\nu $ of $ \nu $ under the map $ C(\R,K) \ni \xi \mapsto
\ti\xi \in C(\R,K) $ given by $ \ti\xi(t) = \xi(-t) $. It corresponds
to the same measure $ \nu_0 $ on $ K $, but a different semigroup $
(\ti T_t) $,
\[
\int_K ( \ti T_t f) g \, d\nu_0 = \int_K f (T_t g) \, d\nu_0
\]
for all $ f,g \in C(K) $. Indeed, such $ (\ti T_t) $ satisfies a
time-reversed version of \eqref{193},
\begin{multline*}
\int_{C (\R,K)} f_1 (\xi(t_1)) \dots f_n (\xi(t_n)) \, d\nu(\xi) =
 \\
= \int_{C (\R,K)} \( \ti T_{t_2-t_1} f_1 \) (\xi(t_2)) f_2 (\xi(t_2))
\dots f_n (\xi(t_n)) \, d\nu(\xi)
\end{multline*}
for $ t_1 < \dots < t_n $.
(Use \eqref{193} for reducing $ n $ by $ 1 $ until $ n=2 $.)

For the semigroup $ (T_t) $ generated by $ L $ of \eqref{196}, the
time-reversed semigroup $ (\ti T_t) $ is generated by the operator
\[
L^* = \frac12 \frac{\pd^2}{\pd y^2} - f(y) \frac{\pd}{\pd x} + \frac12
g^2 (y) \frac{\pd^2}{\pd x^2} \, ,
\]
which shows that the time-reversed process $ \ti \xi (\cdot) $ may be
identified (in distribution) with the process $ \( 1-x(t), y(t) \) $,
where $ \(x(t),y(t)\) = \xi(t) $; here, time reversal appears to be
equivalent to $ x $ reversal.

If the $ x $ reversal leaves $ K_1 $ invariant, that is, $ x \in K_1 $
if and only if $ (1-x) \in K_1 $ for all $ x \in (0,1) $ (which is
evidently compatible with any Hausdorff dimension), then our
random set is symmetric (in distribution), that is, the corresponding
decomposable equivalence class of measures is equal to its time
reversal, which gives us a sharpening of Theorem \ref{6.6}.

\begin{theorem}\label{7.1}
There is a continuum of mutually nonisomorphic symmetric\footnote{%
 Each isomorphic to its opposite, in other words, anti-isomorphic to
 itself.}
product systems of type $ II_0 $.
\end{theorem}

In order to find asymmetric product systems, we need another invariant
of a perfect set; indeed, Hausdorff dimensions of $ K_1 $ and $ \ti
K_1 = \{ x : (1-x) \in K_1 \} $ are evidently equal.

Let us restrict ourselves to perfect sets $ K_1 $ such that each
connected component of the complement $ (\R/\Z) \setminus K_1 $ is an
interval of length of the form, say, $ 10^{-n} $, $ n = 1,2,\dots
$.\footnote{%
 Or, less demanding, of length that belongs to $ \cup_n [10^{-n}, 2
 \cdot 10^{-n} ] $.}
On the other hand, we choose $ f $ such that $ f(\frac14) = 1 $ and $
f(\frac12) \le 2 $. The (random) transformation, relating a portion of
the random set $ C $ to a portion of $ K_1 $, cannot change a length
too much (as explained before Theorem \ref{6.6}). For each component
of $ \R \setminus C $, the corresponding exponent $ n $ is determined
uniquely, unless the component is situated between two portions. Thus,
the combinatorial structure of $ K_1 $ is inherited by $ C $, though a
(locally) finite number of high hierarchical levels are spoiled by
switching between portions. By the combinatorial structure I mean data
like the following: between a pair of nearest intervals of length $
10^{-20} $
there are three intervals of length $ 10^{-21} $; before the first of
these three we see two intervals of length $ 10^{-22} $, while
between the first and second (of the three) we see only one interval
of length $ 10^{-22} $\dots, etc. An infinite array of data may be
encoded in the combinatorial structure of $ K_1 $, and most of the
data are preserved in $ C $. It gives us an alternative way to Theorem
\ref{6.6}, unrelated to the Hausdorff dimension. More important, it gives
us a lot of asymmetric random sets; indeed, the combinatorial
structure need not be symmetric at all. The idea may be implemented as
follows.

Let $ a = (a_k)_{k=1}^\infty $, $ b = (b_k)_{k=1}^\infty $ be two
infinite sequences of numbers $ 0 $ and $ 1 $. Introduce the set $
M_{a,b} $ of all numbers of the form $ \sum_{k=1}^n c_k 3^{-k} $,
where $ n \in \{1,2,\dots\} $, each $ c_k $ is either $ -1 $ or $ +1
$, and the combination $ c_{2k} = c_{2k+1} = -1 $ is forbidden for
each $ k $ such that $ a_k = 0 $; similarly, the combination $
c_{2k} = c_{2k+1} = +1 $ is forbidden whenever $ b_k = 0 $. For any $
x = \sum_{k=1}^n c_k 3^{-k} \in M_{a,b} $ we define $ p(x) =
10^{-n} $, and consider the closure $ K_{a,b} $ of the set
\[
\bigg\{ \sum_{x\in M_{a,b}\cap(-\infty,y)} p(x) : y \in \R \bigg\} \, .
\]

\begin{lemma}
\textup{(a)} For any sequences $ a, b $ the set $ K_{a,b} $ is a
perfect set of Lebesgue measure zero.

\textup{(b)} Let $ x \in K_{a,b} $, $ \eps > 0 $, and $ f :
(x-\eps,x+\eps) \to \R $ be an increasing function such that
\[
\frac15 \le \frac{ f(t) - f(s) }{ t - s } \le 5 \quad
\text{\textup{whenever} } x-\eps \le s < t \le x+\eps \, .
\]
Assume that $ a', b' $ is another pair of sequences such that for
every $ t \in (x-\eps,x+\eps) $,
\[
t \in K_{a,b} \quad \text{if and only if} \quad f(t) \in K_{a',b'} \,
.
\]
Then $ a_k = a'_k $ and $ b_k = b'_k $ for all $ k $ large enough.

\textup{(c)} Assume the same as in \textup{(b)}, but let $ f $ be a
\textup{decreasing} function, $ \frac15 \le \frac{ f(s) - f(t) }{ t -
s } \le 5 $. Then $ a_k = b'_k $ and $ b_k = a'_k $ for all $ k $
large enough.
\end{lemma}

The proof is left to the reader. The conclusion follows.

\begin{theorem}\label{7.2}
There is a continuum of mutually non-isomorphic asymmetric product
systems of type $ II_0 $.
\end{theorem}

It answers in the negative a question of Arveson \cite[p.~6]{Ar89},
namely, whether every product system is symmetric, or not.

\section{FHS space: logarithm of a Hilbert space}

\smallskip

\hfill\parbox{9cm}{%
\dots we have not yet achieved a satisfactory understanding of the
existence of ``logarithms'' for product systems.\\
\mbox{}\hfill W.~Arveson \cite[Sect.~1.6]{Ar-b}
}

\medskip

The classical understanding of the well-known Fock [-Cook]
exponentiation was widely believed to be satisfactory, but is now
changing dramatically. See \cite[Sect.~3.3.1]{Ar} for a concise
summary of the classical understanding; a Hilbert space $ e^H $ arises
from a given Hilbert space $ H $, and $ e^{H_1\oplus H_2} = e^{H_1}
\ot e^{H_2} $ (up to a canonical isomorphism), and the unitary group
of $ H $ acts on $ e^H $ unitarily, by $ U \mapsto \Gamma(U) $. Close
relations to Gaussian random variables and Weyl (canonical)
commutation relations are well-known, see \cite[Chap.~IV and
XIII]{Ja}.

\begin{sloppypar}
Let us start with a finite-dimensional Hilbert space $ H $ over
$ \R $ (that is, Euclidean space). It carries the $ n $-dimensional
Lebesgue measure ($ n = \dim H $), denoted by $ \mes_H $. We may
simply define
\begin{equation}\label{20a}
e^H = L_2 \( (H,\mes_H), \R \) \, ,
\end{equation}
the space of square integrable real-valued functions on $ H $; its
complexification is $ e^H_\C = L_2 \( (H,\mes_H), \C \) $. Why denote
the space by $ e^H \, $? Since $ e^{H_1\oplus H_2} = e^{H_1} \ot
e^{H_2} $ (just because $ \mes_{H_1\oplus H_2} = \mes_{H_1} \times
\mes_{H_2} $). The action $ U \mapsto \Gamma(U) $ is evident: $ \Gamma
(U) f(x) = f ( U^{-1} x ) $. Also, the
additive group of $ H $ acts on $ e^H $ (and $ e^H_\C $) by shifts, $
(V_x f)(y) = f(x+y) $. Note that $ \Gamma(U) V_x = V_{Ux} \Gamma(U) $,
which means that $ \Gamma $ and $ V $ together form a representation
of the group of motions of $ H $.\footnote{%
 Given that $ \Gamma (U) f(x) $ is $ f ( U^{-1} x ) $ rather than $
 f(Ux) $, it would be more consistent to define $ V_x f(y) $ as $
 f(y-x) $ rather than $ f(y+x) $. However, $ f(y+x) $ conforms to the
 tradition of Weyl relations.}
On the complex space $ e^H_\C $ there is another action of the
additive group of $ H $, by such multiplications: $ (U_x f)(y) =
e^{i\ip x y} f(y) $; these unitary operators satisfy the Weyl relations, $
U_x V_y = e^{-i\ip x y} V_y U_x $. There is $ f_0 \in e^H_\C $ such
that $ \ip{V_x f_0}{f_0} = \exp (-\frac18 \|x\|^2 ) $ and $ \ip{U_x
f_0}{f_0} = \exp (-\frac12 \|x\|^2 ) $ for all $ x $,\footnote{%
 Coefficients $ \frac18 $, $ \frac12 $ (rather than $ \frac14 $, $
 \frac14 $) may seem strange, but conform to probabilistic tradition.}
namely, $ f_0 (x) = (2\pi)^{-n/4} \exp (-\frac14 \|x\|^2 ) $ where $ n
= \dim H $; and it is unique up to a phase coefficient. The measure $
\ga_H = |f_0|^2 \cdot \mes_H $, that is, $ (2\pi)^{-n/2} \exp
(-\frac12 \|x\|^2 ) \, dx $, is the standard Gaussian measure on $ H
$. Note that $ \ga_{H_1 \oplus H_2} = \ga_{H_1} \times \ga_{H_2}
$. Thus, we may replace \eqref{20a} with
\begin{equation}\label{21a}
e^H = L_2 \( (H,\ga_H), \R \)
\end{equation}
and get $ e^{H_1\oplus H_2} = e^{H_1} \ot e^{H_2} $ again, as well as
the unitary representation $ \Gamma (U) f(x) = f (U^{-1} x) $. And no
wonder; the two spaces, $ L_2 \( (H,\mes_H), \R \) $ and $ L_2 \(
(H,\ga_H), \R \) $ may be identified by the \emph{canonical} unitary
map, $ L_2 \( (H,\ga_H), \R \) \ni f \mapsto f \cdot |f_0|^2 \in L_2
\( (H,\mes_H), \R \) $. Formulas \eqref{20a}, \eqref{21a} may be united
(recall Sect.~3):
\begin{equation}\label{22a}
e^H = L_2 \( (H,\M_H), \R \) \, ,
\end{equation}
where $ \M_H $ is the equivalence class of measures on $ H $ that
contains $ \mes_H $ and $ \ga_H $. The measure class $ (H,\M_H) $ will
be called the \emph{standard measure class of} $ H $.
\end{sloppypar}

Taking into account that $ \mes_H $ and $ \ga_H $ fail in
infinite dimension, we reformulate the framework. The first twist: we
treat $ \ga_H $ as a measure on the dual space $ H' $ (the space of
all linear maps $ H \to \R $) rather than $ H $. It is basically the
same, due to the evident canonical correspondence between $ H $ and $
H' $, $ x \mapsto \ip \cdot x $. The random element of $ H' $ is a
random linear function on $ H $. The second twist: we restrict the
random function to a basis of $ H $. Nothing is changed in finite
dimension, since a linear functional is uniquely determined by its
restriction to the basis. However, in infinite dimension we escape the
problem of (dis)continuity of the functional. That is the idea; its
implementation follows.

Let $ H $ be an \emph{infinite-dimensional} separable Hilbert space
over $ \R $, and $ E $ an orthonormal basis of $ H $ (treated as a
countable subset of $ H $ rather than a sequence). The space $ \R^E $
of all functions $ \xi : E \to \R $ is the product $ \prod_{e\in
E} \R $ of countably many copies of $ \R $. We equip (each copy of) $
\R $ with the standard Gaussian measure $ \ga^1 $ (just $
(2\pi)^{-1/2} e^{-u^2/2} \, du $) and consider the probability space
\[
(\R, \ga^1)^E = \prod_{e\in E} (\R, \ga^1) = ( \R^E, \ga^E ) \, .
\]
Of course, it is isomorphic to the space of sequences $ (\R, \ga^1)^\N
= ( \R^\N, \ga^\N ) $, $ \N = \{ 1,2,\dots \} $. An evident
isomorphism between $ ( \R^E, \ga^E ) $ and $ ( \R^{E'}, \ga^{E'} ) $
suggests itself for any two orthonormal bases $ E, E' \subset H
$ (provided that the bases are enumerated). However, the evident
isomorphism is irrelevant; what we need, is another, not so evident
isomorphism between these spaces, constructed below.

Given an orthonormal basis $ E $ and a vector $ x \in H $, we may
define $ \xi_x $\index{xi@$ \xi_x $, linear random variable}
for almost all $ \xi \in \R^E $ by
\[
\xi_x = \sum_{e\in E} \ip x e \xi_e \, ;
\]
the series (of orthogonal terms) converges in $ L_2 (\R^E, \ga^E)
$.\footnote{%
 In fact, it converges almost everywhere (provided that the basis is
 enumerated), but we do not need it.}
Thus, we extend $ \xi $ from $ E $ to $ E \cup \{x\} $.\footnote{%
 Note that we do not try to extend $ \xi $ to the whole $ H $.}
It is easy to see that $ \xi_x $ is a normal random variable, $ \xi_x
\sim N(0,\|x\|^2) $. Also, pairs $ (\xi_x,\xi_y) $ have
two-dimensional normal distributions, $ \Ex \xi_x \xi_y = \ip x y $.

Given the second orthonormal basis $ E' \subset H $, we extend $ \xi $
from $ E $ to $ E \cup E' $ and observe that the distribution of the
family $ (\xi_e)_{e\in E'} $ is equal to $ \ga^{E'} $. We get a
(non-evident) isomorphism between the two probability spaces, $ (\R^E,
\ga^E) $ and $ (\R^{E'}, \ga^{E'}) $.

Does it mean that we can glue together probability spaces $ (\R^E,
\ga^E) $ for all bases $ E \, $? On the level of individual points (of
probability spaces), it does not. Indeed, each isomorphism is defined
almost everywhere; the union of a continuum of negligible sets need
not be negligible.\footnote{%
 In fact, here it appears to be co-negligible.}
However, on the level of equivalence classes (of measurable sets or
measurable functions) the glueing is well-defined, and we get a single
Hilbert space
\begin{equation}\label{star}
e^H = L_2 ( \R^E, \ga^E ) \quad \text{for all } E \, ,
\end{equation}
where $ E $ runs over orthonormal bases of $ H $. In particular, for
every vector $ x \in H $, the random variable $ \xi_x \in e^H $ is
well-defined, $ \xi_x \sim N(0,\|x\|^2) $, $ \Ex \xi_x \xi_y = \ip x y
$, and the joint distributions of such random variables are normal. We
get the object well-known as `the isonormal process on $ H $' or `a
Gaussian Hilbert space indexed by $ H $', see \cite[Def.~1.18]{Ja}. A
shorter description is given there, at the expense of considerable
arbitrariness. The description proposed above is canonical; everything
is uniquely determined by $ H $ (not only up to isomorphism), as far
as negligible sets are neglected.

Till now, it is just another description of the classical $ e^H $ for
a real Hilbert space $ H $. Recall however $ (H,\M_H) $ in
\eqref{22a}. For the infinite dimension we may still rewrite
\eqref{star} as
\[
e^H = L_2 ( \R^E, \M^E ) \quad \text{for all } E \, ,
\]
where $ \M^E $ is the equivalence class of measures that contains $
\ga^E $. Doing so, we have no reason to restrict ourselves to
\emph{orthonormal} bases and \emph{measure preserving} maps!

For a finite-dimensional $ H $, each basis $ E $ (orthonormal or not)
can be used. The joint distribution of random variables $ \xi_e $, $ e
\in E $, is some non-degenerate normal distribution (Gaussian measure)
in $ \R^n $ (not just $ \ga^n $); anyway, it belongs to the standard
measure class $ \M_H $.

In infinite dimension we cannot use arbitrary bases of $ H $. Even
such a seemingly innocent case as $ \ip{e_k}{e_l} = 0 $ but, say, $
\ip{e_k}{e_k} = 2 $ (rather than $ 1 $), leads to a distribution
non-equivalent to $ \ga^\N $. Still, \emph{some} non-orthonormal bases
can be used. In other words, we can use bases orthonormal w.r.t.\
\emph{some} other norms on $ H $. The relevant equivalence relation
for norms was described in 1958 by J.~Feldman, J.~H\'ajek and
I.~Segal, see \cite[Th.~3.3.7]{Ar}, \cite[Th.~14.3.1]{Ar-b}. The
following definition is applicable both to real and complex spaces,
however, all FHS spaces will be assumed real, unless otherwise
stated.\index{real versus complex}\index{complex versus real}

\begin{definition}\label{8.1}
(a) An \emph{FHS-equivalence operator}\index{FHS!equivalence operator}
between separable Hilbert
spaces $ H_1, H_2 $ is a linear homeomorphism $ L : H_1 \to H_2 $ such
that $ \One - L^* L $ is a Hilbert-Schmidt operator on $ H_1 $.

(b) \emph{FHS-equivalent norms}\index{FHS!equivalent norms}
on a linear space $ E $ are norms $
\|\cdot\|_1, \|\cdot\|_2 $ on $ E $ such that $ (E,\|\cdot\|_1) $ and
$ (E,\|\cdot\|_2) $ are separable Hilbert spaces, and the identity map
$ E \to E $ is an FHS-equivalence operator between these Hilbert
spaces.

(c) An \emph{FHS space}\index{FHS!space}
is a pair  $ (G,\cN) $ of a linear space $ G $
(over $ \R $, unless it is stated to be complex) and a nonempty set $
\cN $ of norms on $ G $, called \emph{admissible
norms,}\index{admissible norm}
such that

every admissible norm turns $ G $ into a separable Hilbert space;

all admissible norms are mutually FHS-equivalent;

every norm FHS-equivalent to an admissible norm is also admissible.

(d) An \emph{isomorphism}\index{isomorphism!of FHS spaces}
of an FHS space $ (G_1,\cN_1) $ to another
FHS space $ (G_2,\cN_2) $ is a linear homeomorphism $ L : G_1 \to G_2
$ such that $ \| Lx \|_2 = \| x \|_1 $ for some admissible norms $ \|
\cdot \|_1 $ on $ G_1 $ and $ \| \cdot \|_2 $ on $ G_2 $.
\end{definition}

Every Hilbert space is also an FHS space, and every FHS space is also
a linear topological space. Often we say `an FHS space $ G $', leaving
$ \cN $ implicit.

A subset $ E $ of an FHS space $ G $ will be called a
(quasi)orthonormal\footnote{\label{foo18}
 The shorter term `orthonormal' is unambiguous, as far as $ G $ is
 only an FHS space (that is, no one norm is singled out).
 Nevertheless, to be on the safe side, I sometimes use the longer term
 `quasiorthonormal'.}
\index{orthonormal basis of FHS space}
\index{quasiorthonormal basis of FHS space}
basis of $ G $, if there exists an admissible norm $ \|\cdot\| $ on $
G $ such that $ E $ is an orthonormal basis of the Hilbert space $
(G,\|\cdot\|) $.

\begin{definition}
Let $ E, E' $ be (quasi)orthonormal bases of an FHS space $ G $. An
isomorphism between measure classes\footnote{%
 As before, $ \M^E $ is the equivalence class of measures, that
 contains the product measure $ \ga^E $.}
$ (\R^E,\M^E) $ and $ (\R^{E'},\M^{E'}) $ is called
\emph{canonical}, if for every $ x \in G $, the following two
measurable functions (or rather, their equivalence classes) are sent
to each other by the isomorphism:
\[
\R^E \ni \xi \mapsto \sum_{e\in E} \ip x e \xi_e \, , \quad \text{and}
\quad \R^{E'} \ni \xi \mapsto \sum_{e\in E'} \ip x e \xi_e \, ;
\]
here, numbers $ \ip x e $ are defined as the coefficients of the
expansion $ x = \sum_e \ip x e e $ (the sum is over $ E $ or $ E' $
respectively).
\end{definition}

\begin{lemma}
Let $ G $ be an FHS space.

\textup{(a)} For any two \textup{(}quasi\textup{)}orthonormal bases $
E, E' $ of $ G $, a canonical isomorphism between $ (\R^E, \M^E) $ and
$ (\R^{E'}, \M^{E'}) $ exists and is unique \textup{($ \bmod \, 0 $).}

\textup{(b)} For any three \textup{(}quasi\textup{)}orthonormal bases
$ E, E', E'' $ of $ G $, canonical isomorphisms make the diagram
\[
\xymatrix{
(\R^E, \M^E)
  \ar[rr]
  \ar[rd] & &
(\R^{E''}, \M^{E''}) \\
& (\R^{E'}, \M^{E'})
  \ar[ru]
}
\]
commutative. In addition, the canonical isomorphism $ (\R^E, \M^E)
\to (\R^E, \M^E) $ is the identity map.
\end{lemma}

The proof is left to the reader.

So, every FHS space $ G $ leads to a measure class $ (X_G, \M_G) $
that is canonically isomorphic $ \bmod \, 0 $ to each $ (\R^E, \M^E)
$.\footnote{%
 More formally: a measurable set ($ \bmod \, 0 $) in $ (X_G, \M_G) $
 is, by definition, a family $ A = (A_E) $ such that $ A_E $ is a
 measurable set ($ \bmod \, 0 $) in $ (\R^E, \M^E) $ for every
 orthonormal basis $ E $ of $ G $, and for any two such bases, $
 A_{E_1} $ and $ A_{E_2} $ are sent to each other ($ \bmod \, 0 $) by
 the canonical isomorphism between $ (\R^{E_1}, \M^{E_1}) $ and $
 (\R^{E_2}, \M^{E_2}) $. Measurable functions ($ \bmod \, 0 $) on $
 (X_G, \M_G) $ are defined similarly. Note that we do not try to define
 a point of $ X_G $ in a canonical way.}
We call $ (X_G, \M_G) $ the standard measure class of $ G
$.\footnote{%
 Recall that for finite dimension we took $ X_H = H $. For infinite
 dimension, $ X_G $ may be thought of as some extension of the
 topological dual $ G' $, and then $ G' $ is a negligible subset of $
 X_G $.}
Each vector
$ x \in G $ leads to a measurable function $ \xi_x : X_G \to \R $ (or
rather, equivalence class), and the map $ x \mapsto \xi_x $ is
linear. Every admissible norm $ \|\cdot\| $ on $ G $ leads to a
measure $ \ga_{\|\cdot\|} $ on $ X_G $, $ \ga_{\|\cdot\|} \in \M_G
$. Such measures will be called \emph{Gaussian
measures}\index{Gaussian measure!over FHS space}
of $ \M_G
$. With respect to the measure $ \ga_{\|\cdot\|} $ we have
\begin{equation}\label{8.35}
\xi_x \sim N(0,\|x\|^2) \quad \text{for all } x \in G \, .
\end{equation}
The direct sum of FHS spaces corresponds to the product of their
standard measure classes,
\begin{equation}\label{24.5}
( X_{G_1 \oplus G_2}, \M_{G_1 \oplus G_2} ) = ( X_{G_1}, \M_{G_1} )
\times ( X_{G_2}, \M_{G_2} )
\end{equation}
up to a canonical isomorphism.

\begin{definition}
The \emph{exponential}\index{exponential of FHS space}
of an FHS space $ G $ over $ \R $ is such a Hilbert space (over $ \R
$):
\[
e^G = L_2 \( (X_G,\M_G), \R \) \, ,
\]
where $ (X_G,\M_G) $ is the standard measure class of $ G $. Its
complexification, $ e^G_\C = L_2 \( (X_G,\M_G), \C \) $, is the
complex exponential of $ G $.
\end{definition}

Any isomorphism between FHS spaces, $ W : G_1 \to G_2 $, leads to the
corresponding isomorphism of measure classes $ X_1 \to X_2 $; namely,
$ \om_1 \mapsto \om_2 $ means $ \xi_x (\om_1) = \xi_{Wx} (\om_2)
$. Also, a non-homogeneous linear transformation may be introduced: $
\om_1 \mapsto \om_2 $ means $ \xi_x (\om_1) = \xi_{Wx} (\om_2) + \ip x
y $, where $ y \in G'_1 $ is a given linear functional on $ G_1 $. The
isomorphism of measure classes induces a unitary operator $
\Gamma(W,y) = U : e^{G_1} \to
e^{G_2} $\index{Gamma@$ \Gamma $, action on exponential space}
between Hilbert spaces;
namely, $ U (f_1 \sqrt{\mu_1}) = f_2 \sqrt{\mu_2} $ whenever $
f_1(\om_1) = f_2(\om_2) $ and $ \mu_1(d\om_1) = \mu_2(d\om_2) $. Note
that $ \Gamma (\One,y) \Gamma(W,0) = \Gamma (W,W^* y) = \Gamma (W,0)
\Gamma (\One,W^* y) $; of course, $ \ip x {W^* y} = \ip {Wx} y $.

\begin{lemma}\label{8.5}
\textup{(a)} $ e^{G_1 \oplus G_2} = e^{G_1} \ot e^{G_2} $ \textup{(}up
to canonical isomorphism\textup{),} for any FHS spaces $ G_1, G_2 $.

\textup{(b)} The group of automorphisms of an FHS space $ G $ acts on
$ e^G $ by unitary operators, $ W \mapsto \Gamma(W) $.

\textup{(c)} The additive group of the \textup{(topological)} dual
space $ G' $ of an FHS space $ G $ acts on $ e^G $ \textup{(}and $
e^G_\C $\textup{)} by unitary operators $ 
V_y $\index{Vy@$ V_y $, action on exponential space}
such that $ V_y \Gamma(W) = \Gamma(W) V_{W^* y} $; the additive
group of $ G $ acts on $ e^G_\C $ by unitary operators
$ U_x $;\index{Ux@$ U_x $, action on exponential space}
these
operators satisfy Weyl relations, $ U_x V_y = e^{-i\ip x y} V_y U_x
$. For each Gaussian measure $ \ga \in \M_G $ we have $ \ip{V_y
\sqrt\ga}{\sqrt\ga} = \exp(-\frac18 \|y\|^2) $ and $ \ip{U_x
\sqrt\ga}{\sqrt\ga} = \exp(-\frac12 \|x\|^2) $ for all $ x \in G $, $
y \in G' $; here $ \|\cdot\| $ is the admissible norm on $ G $
corresponding to $ \ga $ \textup{(}or its dual norm on $ G'
$\textup{).} Also, $ e^G $
is spanned by $ \{ U_x \sqrt\ga : x \in G \} $, as well as by $ \{ V_y
\sqrt\ga : y \in G' \} $.

\textup{(d)} Let $ G $ be an FHS space, $ (X,\B,\M) $ a measure
class, $
\ga\in\M $ a probability measure, and $ \ti\eta : G \to L_2(\ga) $ a
linear map, isometric from $ (G,\|\cdot\|) $ to $ L_2(\ga) $ for some
admissible norm $ \|\cdot\| $ on $ G $, and each $ \ti\eta (x) $ has
the normal distribution $ N(0,\|x\|^2) $ w.r.t.\ $ \ga $, and the
image $ \ti\eta (G) \subset L_2(\ga) $ generates the whole \sif\ $ \B
$. Then there exists one and only one isomorphism of measure classes $
\eta : (X_G,\M_G) \to (X,\M) $ inducing $ \ti\eta $ by $ \ti\eta(x)(\om)
= x (\eta^{-1}\om) $ for $ \om \in X $, $ x \in G $.
\end{lemma}

The proof is left to the reader.

The next result implies that the obstacle mentioned in Remark
\ref{3.5a} does not appear when dealing with Gaussian measures. Its
natural domain is linear topological spaces, but we need only FHS
spaces. See also Remark \ref{9.05}.

\begin{lemma}\label{8.6}
Let $ G $ be a linear topological space and $ G_1, G_{1,2}, G_{2,3},
G_3 $ its \textup{(closed linear)} subspaces such that $ G_1 \subset
G_{1,2} $, $ G_{2,3} \supset G_3 $, and
\[
G_{1,2} \oplus G_3 = G = G_1 \oplus G_{2,3} \, .
\]
Then the subspace $ G_2 = G_{1,2} \cap G_{2,3} $ satisfies
\[
G_1 \oplus G_2 \oplus G_3 = G \, , \quad
G_1 \oplus G_2 = G_{1,2} \, , \quad
G_2 \oplus G_3 = G_{2,3} \, .
\]
\textup{(Direct sums are treated in the topological sense.)}
\end{lemma}

\begin{proof}
The decomposition $ G = G_1 \oplus G_{2,3} $ determines a projection $
P_1 : G \to G $ such that $ P_1 G = G_1 $ and $ 1-P_1 = P_{2,3} $ is
also a projection, $ P_{2,3} G = G_{2,3} $. The same for $ P_3 $ and $
P_{1,2} = 1 - P_3 $. The inclusion $ G_1 \subset G_{1,2} $ gives $ P_1
= P_{1,2} P_1 = (1-P_3) P_1 $, that is, $ P_3 P_1 = 0 $; similarly, $
P_1 P_3 = 0 $. So, our projections commute with each
other. Introducing
\[
P_2 = P_{12} P_{23} = P_{23} P_{12} = (1-P_1) (1-P_3) = 1 - P_1 - P_3
\, ,
\]
we have $ P_2^2 = P_{1,2} P_{2,3} P_{1,2} P_{2,3} = P_{1,2}^2
P_{2,3}^2 = P_{1,2} P_{2,3} = P_2 $; that is, $ P_2 $ is also a
projection, and $ P_1 + P_2 + P_3 = 1 $. It follows that $ P_2 G = \(
(P_1+P_2) G \) \cap \( (P_2+P_3) G \) = G_{1,2} \cap G_{2,3} = G_2 $.
\end{proof}

\section{Continuous sums and off-white noises}

Roughly speaking, a continuous sum of FHS spaces is a family $
(G_{s,t})_{s<t} $ of FHS spaces $ G_{s,t} $ such that $ G_{r,s} \oplus
G_{s,t} = G_{r,t} $. One can turn the idea into a definition in the
style of \ref{1.1} (and \ref{3.2}) by introducing isomorphisms $
W_{r,s,t} : G_{r,s} \oplus G_{s,t} \to G_{r,t} $ of FHS
spaces. However, I prefer an equivalent definition in the style of
\ref{3.5}.

A (closed linear) subspace $ G_1 $ of an FHS space $ G $ inherits from
$ G $ the structure of an FHS space. Given two subspaces $ G_1, G_2
\subset G $, we may consider the map $ G_1 \oplus G_2 \ni x \oplus y
\mapsto x+y \in G $ and ask whether it is an isomorphism of FHS
spaces, or not. If it is, we say that $ G_1, G_2 $ are
(quasi)orthogonal\footnote{\label{foo24}
 Similarly to Footnote \ref{foo18}, the shorter term is unambiguous,
 but the longer term is safer.}%
\index{orthogonal FHS subspaces}
\index{quasiorthogonal FHS subspaces}
(see \cite[Sect.~3.3.3]{Ar}, \cite[Sect.~14.1.1]{Ar-b}) and write `$ G
= G_1 \oplus G_2 $ in the FHS sense'.\index{FHS!sense}
Do not confuse it with the relation `$ G = G_1 \oplus G_2 $ in the
topological sense';\index{topological sense}
the latter
means that the map $ G_1 \oplus G_2 \ni x \oplus y \mapsto x+y \in G $
is a linear homeomorphism, but not necessarily an isomorphism of FHS
spaces. Similarly, the relation `$ G = G_1 \oplus G_2 \oplus G_3 $ in
the FHS sense' means, by definition, that the map $ G_1 \oplus G_2
\oplus G_3 \ni x \oplus y \oplus z \mapsto x+y+z \in G $ is an
isomorphism of FHS spaces.

\begin{remark}\label{9.05}
Let $ G $ be an FHS space and $ G_1, G_{1,2}, G_{2,3}, G_3 $ its
subspaces such that $ G_1 \subset G_{1,2} $, $ G_{2,3} \supset G_3 $,
and
\[
G_{1,2} \oplus G_3 = G = G_1 \oplus G_{2,3} \quad \text{in the FHS
sense.}
\]
Then the subspace $ G_2 = G_{1,2} \cap G_{2,3} $ satisfies
\[
G_1 \oplus G_2 \oplus G_3 = G \, , \quad
G_1 \oplus G_2 = G_{1,2} \, , \quad
G_2 \oplus G_3 = G_{2,3} \quad \text{in the FHS sense.}
\]

Indeed, we apply Lemma \ref{8.6} and combine FHS isomorphisms $ G_1
\oplus G_{2,3} \to G $ and $ G_2 \oplus G_3 \to G_{2,3} $ into $ G_1
\oplus G_2 \oplus G_3 \to G_1 \oplus G_{2,3} \to G $.
\end{remark}

\begin{definition}
A \emph{continuous sum of FHS spaces}\index{continuous sum of FHS spaces}
is $ \( G, (G_{s,t})_{-\infty\le
s<t\le\infty} \) $ where $ G $ is an FHS space, and $ G_{s,t} \subset
G $ are subspaces such that
\[
G_{r,s} \oplus G_{s,t} = G_{r,t} \quad \text{whenever } -\infty \le r
< s < t \le \infty \, ,
\]
and $ G_{-\infty,\infty} = G $.
\end{definition}

\begin{corollary}
If $ \( G, (G_{s,t})_{s<t} \) $ is a continuous sum of FHS spaces,
then $ (e^{G_{s,t}})_{s<t} $ is a real \cpHs, and $
(e_\C^{G_{s,t}})_{s<t} $ is a \textup{(complex)} \cpHs.
\end{corollary}

(To be rigorous, canonical isomorphisms $ W_{r,s,t} : e^{G_{r,s}} \ot
e^{G_{s,t}} \to e^{G_{r,s}\oplus G_{s,t}} \to e^{G_{r,t}} $ should be
included.)

The transition from the continuous sum to the \cpHs\ is mediated by a
\cpmc. Namely, $ e^{G_{s,t}} = L_2 \( (X_{s,t}, \M_{s,t}), \R \) $ and
$ e_\C^{G_{s,t}} = L_2 \( (X_{s,t}, \M_{s,t}), \C \) $,
where $ (X_{s,t}, \M_{s,t}) $ is the standard measure class of $
G_{s,t} $; and $ (X_{r,s}, \M_{r,s}) \times (X_{s,t}, \M_{s,t}) =
(X_{r,t}, \M_{r,t}) $ by \eqref{24.5}. Alternatively, we may treat the
\cpmc\ according to \ref{3.5} (rather than \ref{3.2}); that is, on the
standard measure class $ (X,\B,\M) $ of $ G $ we have sub-\sif s $
\B_{s,t} $ such that $ \B_{r,s} \otimes \B_{s,t} = \B_{r,t} $. Namely,
$ \B_{s,t} $ is generated by measurable functions $ \xi_x $, $ x \in
G_{s,t} $. A similar transition from a \emph{homogeneous} continuous
sum of FHS spaces to a \hcpHs\ via a \hcpmc\ is discussed below.

\begin{definition}
\begin{sloppypar}
A \emph{homogeneous continuous sum of FHS
spaces}\index{homogeneous!continuous sum of FHS spaces}
is $ \( G,
(G_{s,t})_{-\infty\le s<t\le\infty}, (\ti\theta^u)_{u\in\R} \) $ where
$ \( G, (G_{s,t})_{-\infty\le s<t\le\infty} \) $ is a continuous sum
of FHS spaces, and $ (\ti\theta^u)_{u\in\R} $ is a one-parameter group
of automorphisms $ \ti\theta^u $ of the FHS space $ G $ such that
\[
\ti\theta^u (G_{s,t}) = G_{s+u,t+u} \quad \text{whenever } -\infty \le
s < t \le \infty, u \in \R \, , 
\]
and the corresponding unitary group $ \theta^u = \Gamma (\ti\theta^u)
$ on $ e^G $ is strongly continuous.
\end{sloppypar}
\end{definition}

\begin{sloppypar}
No need to elaborate on the latter condition, since all we need is
such a simple sufficient condition:
\begin{equation}\label{24}
\parbox{10.5cm}{%
there exists an admissible norm $ \|\cdot\| $ on $ G $ such that
$ (\ti\theta^u)_{u\in\R} $ is a strongly continuous \emph{unitary}
group on the Hilbert space $ (G,\|\cdot\|) $.}
\end{equation}
Indeed, it implies that $ \theta^u \sqrt{\ga} = \sqrt{\ga} $ (here $
\ga = \ga_{\|\cdot\|} $) and $ \theta^u V_y \sqrt{\ga} =
V_{(\ti\theta^{-u})^* y} \theta^u \sqrt{\ga} \to V_y \sqrt{\ga} $ when
$ u \to 0 $, for each $ y \in G' $; here $ \theta^u = \Gamma
(\ti\theta^u) $. Taking into account that the set $ \{ V_y \sqrt{\ga}
: y \in G' \} $ spans $ e^G $ we conclude that $ \theta^u \to \One $
strongly on $ e^G $ when $ u \to 0 $.
\end{sloppypar}

Let $ \( G, (G_{s,t})_{s<t}, (\ti\theta^u)_{u\in\R} \) $ be a
homogeneous continuous sum of FHS spaces, and $ (e^{G_{s,t}})_{s<t} $
the corresponding real \cpHs. The restriction $ \ti\theta^u_{s,t} $ of $
\theta^u $ to $ G_{s,t} $, being an isomorphism $ \ti\theta^u_{s,t} :
G_{s,t} \to G_{s+u,t+u} $ of FHS spaces, leads to the corresponding
unitary operator $ \theta^u_{s,t} : e^{G_{s,t}} \to e^{G_{s+u,t+u}}
$. Clearly, $ \ti\theta^u_{r,s} \oplus \ti\theta^u_{s,t} =
\ti\theta^u_{r,t} $, which ensures \ref{1.2}(a1). Condition
\ref{1.2}(a2) follows from the relation $ \ti\theta^v \ti\theta^u =
\ti\theta^{u+v} $. Condition \ref{1.2}(a3) is the strong continuity.

\begin{corollary}\label{9.4}
If $ \( G, (G_{s,t})_{s<t}, (\ti\theta^u)_{u\in\R} \) $ is a
homogeneous continuous sum of FHS spaces, then $ \(
(e^{G_{s,t}})_{s<t}, (\theta^u_{s,t})_{s<t,u\in\R} \) $ is a real
\hcpHs, and its complexification is a \textup{(complex)} \hcpHs.
\end{corollary}

In combination with Theorem \ref{1.5} it means that each homogeneous
continuous sum of FHS spaces leads to a product
system.\index{product system!out of FHS spaces}

The simplest nontrivial example of a homogeneous continuous sum of FHS
spaces is given by $ G_{s,t} = L_2 (s,t) \subset L_2 (\R) $. These are
Hilbert spaces, and $ G_{r,s} \oplus G_{s,t} = G_{r,t} $ in the
Hilbert-space sense, which is much better than needed. The
corresponding Gaussian measure describes the white
noise;\index{noise!white}
the \hcpmc\ is in fact an \hcpps, and the product system is of type $ I $.

\begin{sloppypar}
In order to get more interesting product systems, we modify the metric
of $ L_2(\R) $, keeping it shift-invariant. The idea is to define $
\| f \|_C^2 = \iint f(s) f(t) { C(s-t) } \, ds \, dt $; however,
technically, it is better to use the Fourier transform, $ \hat f(\la)
= (2\pi)^{-1/2} \int f(t) e^{-i\la t} \, dt
$,\index{Fourier transform, normalized to be unitary}
and define
\[
\| f \|_\nu^2 = \int | \hat f(\la) |^2 \, \nu(d\la) \, .
\index{zz@$ "\"|\cdot"\"|_\nu $, shift-invariant norm}
\]
The idea is implemented as follows.
\end{sloppypar}

Let $ \nu $\index{nu@$ \nu $, spectral measure on $ [0,\infty) $}
be a positive $ \sigma $-finite measure on $ [0,\infty) $ 
such that $ \int (1+\la^2)^{-m} \, \nu(d\la) < \infty $ for $ m $
large enough. We exclude the trivial case $ \nu = 0 $, and consider
the Hilbert space $ L_2(\nu) $, and its subspace $ G $ spanned by
functions $ \hat f $ when $ f $ runs over all compactly supported
functions $ f : \R \to \R $ of class $ C^{2m} $ (in fact, $ G $ is the
whole $ L_2(\nu)
$). Requiring $ f $ to vanish on $ (t,\infty) $ we get a subspace $
G_{-\infty,t} $ of $ G $. Similarly, $ G_{t,\infty} $ is spanned by $
\hat f $ for $ f $ vanishing on $ (-\infty,t) $. Despite the complex
nature of $ L_2(\nu) $, we downgrade $ G $ to a Hilbert space over $
\R $ and further, to a real FHS space. Only the case $ m=1 $ is
important for us (see \cite{Ts02} for the general case).

\begin{theorem}\label{9.5} \cite[Th.~3.2]{Ts02}.
The following two conditions on $ \nu $ are equivalent:

\textup{(a)} $ G = G_{-\infty,0} \oplus G_{0,\infty} $ in the FHS
sense, and $ \int (1+\la^2)^{-1} \, \nu(d\la) < \infty $;

\textup{(b)} $ \nu(d\la) = e^{\phi(\la)} \, d\la $ for some $ \phi :
[0,\infty) \to \R $ such that
\begin{equation}\label{25a}
\int_0^\infty \int_0^\infty \frac{ | \phi (\la_1) - \phi (\la_2) |^2
}{ | \la_1 - \la_2 |^2 } \, d\la_1 \, d\la_2 < \infty \, .
\end{equation}
\end{theorem}

\begin{lemma}\label{9.6}
If $ \nu $ satisfies the \textup{(equivalent)} conditions of Theorem
\textup{\ref{9.5},} then subspaces $ G_{s,t} = G_{-\infty,t} \cap
G_{s,\infty} $ of $ G $ form a continuous sum of FHS spaces.
\end{lemma}

\begin{proof}
First, $ G_{r,s} \oplus G_{s,t} = G_{r,t} $ in the topological sense,
which is proven similarly to Lemma \ref{8.6}. Second, $ G_{r,s} $ and $
G_{s,t} $ are quasi-orthogonal, since $ G_{-\infty,s} $ and $
G_{s,\infty} $ are.
\end{proof}

Unitary operators $ \ti\theta^u : L_2(\nu) \to L_2(\nu) $ defined by $
\ti\theta^u g(\la) = e^{-iu\la} g(\la) $ satisfy $ \ti\theta^u
(G_{s,t}) = G_{s+u,t+u} $ and \eqref{24}. Therefore $ \( G,
(G_{s,t})_{s<t}, (\ti\theta^u)_{u\in\R} \) $ is a homogeneous
continuous sum of FHS spaces, and by Corollary \ref{9.4}, Hilbert
spaces
\[
H_{s,t} = e^{G_{s,t}}
\]
form a real \hcpHs\ (which leads to a product system).

The theory of random processes can contribute to our understanding of the
spaces $ G $ and $ H = e^G $. It is well-known that every positive $
\sigma $-finite measure $ \nu $ on $ [0,\infty) $ such that $ \int
(1+\la^2)^{-m} \, \nu(d\la) < \infty $ (for $ m $ large enough) is the
spectral measure of a stationary Gaussian generalized random process $
\xi(\cdot) $. In other words, there exists one and only one
shift-invariant Gaussian
measure\index{Gaussian measure!on generalized functions}
$ \ga_\nu $ on the space $ \Om $ of
(tempered, Schwartz; real-valued) distributions\footnote{%
 Sorry, `a distribution in the space of distributions' may be
 confusing. A `probability distribution' is just a probability measure
 (intended to describe a random element of the corresponding space). In
 contrast, a generalized function, also called `distribution', is a
 more singular (than a measure) object over $ \R $, generally not
 positive; for example, a derivative $ \delta^{(n)} $ of Dirac's
 delta-function.}
(generalized functions) over $ \R $ such that
\[
\int_\Om \bigg( \int_\R f(t) \xi(t) \, dt \bigg)^2 \, \ga_\nu (d\xi) =
\| \hat f \|^2_{L_2(\nu)}
\]
for every compactly supported $ f $ of class $ C^{2m} $. More
probabilistically,
\begin{equation}\label{27}
\int_\R f(t) \xi(t) \, dt \sim N \( 0, \| \hat f \|^2_{L_2(\nu)} \) \,
.
\end{equation}
Lemma \ref{8.5}(d) gives us a canonical isomorphism between two
measure classes, the standard measure class $ (X_G,\M_G) $ of the FHS
space $ G $ and $ (\Om,\M_\nu) $, where $ \M_\nu $ is the equivalence
class containing the measure $ \ga_\nu $. We treat them as identical:
\[
(X_G,\M_G) = (\Om,\M_\nu) \, ; \quad H = L_2 \( (\Om,\M_\nu), \R \) =
L_2 \( (\Om,\ga_\nu), \R \) \, ;
\]
the random variable $ \int f(t) \xi(t) \, dt $ is nothing but $ \xi_x
$ (as in \eqref{8.35}) for $ x = f $. However, smooth test functions $
f $ represent \emph{continuous} linear functionals on $ \Om $, while
points $ x $ of $ G $ represent \emph{measurable} linear
functionals\index{measurable linear functional}
on $ \Om $. The corresponding $ \hat f $ are dense in $ G $, and $
\hat f $ for $ f $ that vanish on $
(-\infty,0) $ are dense in $ G_{0,\infty} $.

If $ \nu $ satisfies equivalent conditions of Theorem
\ref{9.5}, then the random process $ \xi(\cdot) $ is a so-called
off-white noise.\index{noise!off-white}
Its past (the restriction of $ \xi(\cdot) $ to $
(-\infty,0) $) and future (the restriction of $ \xi(\cdot) $ to $
(0,\infty) $) are stochastically dependent (unless the noise is
white), but not too much; the product of marginal distributions is
equivalent to the given distribution.

From now on we restrict ourselves to measures $ \nu $ that satisfy
equivalent conditions of Theorem \ref{9.5}.

The function $ e^{\phi(\la)} = \nu(d\la) / d\la $ will be called the
\emph{spectral
density.}\index{spectral density (or measure) on $ [0,\infty) $!of
 off-white noise}
Condition \ref{9.5}(b) ensures that $
\nu(d\la) / d\la > 0 $ almost everywhere;\footnote{%
 In fact, $ \int \( \frac{ \nu(d\la) }{ d\la } \)^p \, \frac{ d\la }{
 1+\la^2 } < \infty $ for all $ p \in (-\infty,\infty) $, but we do
 not need it.}
thus $ \nu $ is equivalent to the Lebesgue measure.

A sufficient condition for \eqref{25a} is available (see
\cite[Prop.~3.6(b)]{Ts02}):
\begin{equation}\label{28}
\text{$ \phi $ is continuously differentiable, and} \quad
\int_0^\infty \bigg| \frac{d}{d\la} \phi (\la) \bigg|^2 \la \,
 d\la < \infty \, .
\end{equation}
In particular, the sufficient condition is satisfied by any strictly
positive smooth function $ \la \mapsto e^{\phi(\la)} = \nu(d\la) /
d\la $ such that for $ \la $ large enough, one of the following
equalities holds:
\begin{gather}
\frac{\nu(d\la)}{d\la} = \ln^\al \la \, , \quad -\infty < \al <
 \infty \, ; \label{293} \\
\frac{\nu(d\la)}{d\la} = \exp ( - \ln^\al \la ) \, , \quad 0 < \al <
 \frac12 \, ; \label{294}
\end{gather}
see \cite[Examples 3.11, 3.12]{Ts02}.

\begin{sloppypar}
The case of $ \ln^\al \la $ with $ \al > 0 $ shows that the spectral
density of an off-white noise need not be bounded.
From now on we assume boundedness (in addition to
\ref{9.5}(b)):\footnote{%
 After an appropriate change of $ \nu(d\la) / d\la $ on a negligible
 set, of course.}
\begin{equation}\label{29}
\sup_{\la\in[0,\infty)} \frac{\nu(d\la)}{d\la} < \infty \, .
\end{equation}
(Thus, $ \ln^\al \la $ fits only for $ -\infty < \al \le 0 $.)
It follows that $ \| \hat f \|_{L_2(\nu)} \le \const \cdot \| \hat f
\|_{L_2(0,\infty)} = \const \cdot \| f \|_{L_2(\R)} $. We have a
continuous operator $ \iota : L_2 (\R) \to G $,\footnote{%
 By $ L_2(\R) $ I mean the space of real-valued (not complex-valued)
 functions on $ \R $.}
\begin{align*}
G &\text{ is the closure of } \iota ( L_2(\R) ) \, ; \\
G_{-\infty,t} &\text{ is the closure of } \iota ( L_2(-\infty,t) ) \,
 ; \\
G_{s,\infty} &\text{ is the closure of } \iota ( L_2(s,\infty) ) \,
 .
\end{align*}
In contrast to $ C^{2m} $, the space $ L_2(\R) $ is closed under
multiplication by the indicator $ \One_{(s,t)} $ of an interval $
(s,t) \subset \R $. Projections $ Q_{-\infty,t} $ and $ Q_{t,\infty} =
\One - Q_{-\infty,t} $ on $ G $, corresponding to the decomposition $
G = G_{-\infty,t} \oplus G_{t,\infty} $, satisfy
\begin{align*}
Q_{-\infty,t} \iota(f) = \iota ( f \cdot \One_{(-\infty,t)} ) \, ; \\
Q_{s,\infty} \iota(f) = \iota ( f \cdot \One_{(s,\infty)} ) \, .
\end{align*}
Therefore the projection $ Q_{s,t} = Q_{-\infty,t} Q_{s,\infty} =
Q_{s,\infty} Q_{-\infty,t} $ satisfies
\[
Q_{s,t} \iota(f) = \iota ( f \cdot \One_{(s,t)} ) \, ,
\]
which implies that
\[
G_{s,t} \text{ is the closure of } \iota ( L_2(s,t) ) \, .
\]
The random variable 
$ \xi_x $\index{xi@$ \xi_x $, linear random variable}
(recall \eqref{8.35}) corresponding to $ x
= \iota(f) $ where $ f \in L_2 (\R) $ will be \emph{denoted} by
\[
\int_{-\infty}^\infty f(t) \xi(t) \, dt \, .
\index{zz@$ \int f(t) \xi(t) \, dt $, linear random variable}
\]
Now, \eqref{27} holds for all $ f \in L_2 (\R) $. When $ f $ is good
enough, the new meaning of $ \int f(t) \xi(t) \, dt $ conforms $
\ga_\nu $-almost everywhere to the meaning used in \eqref{27}. Of
course, $ \int_s^t f(u) \xi(u) \, du $ stands for $
\int_{-\infty}^\infty \( f(u) \One_{(s,t)} (u) \) \xi(u) \, du $. Note
that $ \int_s^t f(u) \xi(u) \, du $ is measurable w.r.t.\ the
sub-\sif\ $ \B_{s,t} $ corresponding to $ G_{s,t} $ (these sub-\sif s
form a \cpmc). Also, $ \B_{s,t} $ is generated by $ \{ \int_s^t f(u)
\xi(u) \, du : f \in L_2 (s,t) \} $.
\end{sloppypar}

\section{Type $ III $}

\smallskip

\hfill\parbox{9cm}{%
We still do not know a ``natural'' construction of a product system
without any units\dots\ We believe that there should be a natural way
of constructing such product systems, and we offer that as a basic
unsolved problem.\\
\mbox{}\hfill W.~Arveson \cite[Introduction]{Ar97}
}

\medskip

We consider the homogeneous continuous sum of FHS spaces constructed
in Sect.~9 for a spectral density that satisfies \ref{9.5}(b),
\eqref{29} and
\begin{equation}\label{26}
\frac{ \nu(d\la) }{ d\la } \to 0 \quad \text{for } \la \to \infty
\end{equation}
(the white noise is now excluded). It will be proven that the
corresponding product system is of type $ III $. Recall that $ \| f
\|_\nu = \| \hat f \|_{L_2(\nu)} $ for $ f \in L_2(\R) $.

\begin{lemma}\label{10.1}
Let $ f \in L_2 (\R) $. Define $ f_n \in L_2(\R) $ by
\[
f_n (t) = f(t) \cdot \sgn \sin \pi n t = \begin{cases}
 f(t) &\text{\textup{for} $ t \in \cup_{k\in\Z} ( \frac{2k}n,
  \frac{2k+1}n ) $},\\
 -f(t) &\text{\textup{otherwise}}.
\end{cases}
\]
Then $ \| f_n \|_\nu \to 0 $ for $ n \to \infty $.
\end{lemma}

\begin{proof}
Taking into account that $ \|\cdot\|_\nu \le \const \cdot \|\cdot\| $
(here $ \|\cdot\| $ is the usual norm of $ L_2(\R) $), we may restrict
ourselves to $ f $ of a dense subset of $ L_2(\R) $. Thus, we assume
that $ \hat f $ has a compact support.

The Fourier transform of the periodic function $ t \mapsto \sgn \sin
\pi n t $ is concentrated on the (shifted) lattice $ \pi n (1+2\Z) =
\{ \dots, -3\pi n, -\pi n, \pi n, 3\pi n, \dots \} $. Let $ n $ be
large enough, such that the interval $ [-\pi n, \pi n] $ contains the
support of $ \hat f $. Then $ \hat f_n $ consists of non-overlapping
shifted copies of $ \hat f $, multiplied by Fourier coefficients of
the periodic function. We have
\[
\| f_n \|^2_\nu = \int | \hat f_n (\la) |^2 \nu(d\la) \, d\la \le
\const \cdot \| \hat f \|^2 \cdot \sup_{\la \in \supp \hat f_n} \frac{
\nu(d\la) }{ d\la } \to 0 \, ;
\]
indeed, for large $ n $, the support of $ \hat f_n $ is distant from
the origin.
\end{proof}

The following result holds for real and complex \cpHss\ equally well.

\begin{lemma}\label{10.2}
The \cpHs\ over $ [0,1] $, corresponding to the considered continuous
sum of FHS spaces, is of type $ III $.
\end{lemma}

\begin{proof}
The \cpHs\ arises from a \cpmc\ formed by sub-\sif s $
(\B_{s,t})_{0\le s<t\le1} $ on the measure class $ (\Om,\M_\nu) $. We
choose $ f \in L_2(0,1) $ such that $ \| f \|_\nu \ne 0 $ (in fact,
every $ f \ne 0 $ fits, since $ \nu $ is equivalent to the Lebesgue
measure). With respect to the Gaussian measure $ \ga_\nu $, the random
variable $ \int_0^1 f(u) \xi(u) \, du $ has a non-degenerate normal
distribution. Therefore $ \int_0^1 f(u) \xi(u) \, du $ has an
absolutely continuous distribution for every probability measure
of $ \M_\nu $, and moreover, for every probability measure
absolutely continuous w.r.t.\ $ \M_\nu $.

Assume that the \cpHs\ is not of type $ III $, and let $ \psi \in
H_{0,1} $ be a decomposable vector, $ \| \psi \| = 1 $. The
probability measure $ \mu = |\psi|^2 $ on $ (\Om,\B_{0,1}) $ is
absolutely continuous w.r.t.\ $ \M_\nu $, and makes sub-\sif s $ \B_{r,s}
$ and $ \B_{s,t} $ independent whenever $ 0 \le r < s < t \le 1 $ (see
Lemma \ref{4.1}). Therefore it makes independent the following two
random variables (for each $ n $ separately):
\[
\int_0^1 ( f(u) - f_n(u) ) \xi(u) \, du \quad \text{and} \quad
\int_0^1 ( f(u) + f_n(u) ) \xi(u) \, du \, ,
\]
where $ f_n $ are as in Lemma \ref{10.1}. (Indeed, it makes
the $ n $ sub-\sif s $ \B_{k/n,(k+1)/n} $ independent.)

By Lemma \ref{10.1} and \eqref{27}, measurable functions $ \int_0^1
f_n(u) \xi(u) \,
du $ converge to $ 0 $ (for $ n \to \infty $) in $ L_2 (\ga_\nu) $,
therefore, in measure w.r.t.\ $ \ga_\nu $, and therefore, also in
measure w.r.t.\ $ \mu $. So, w.r.t.\ $ \mu $ the random variable
\[
\int_0^1 f(u) \xi(u) \, du = \lim_n \int_0^1 ( f(u) - f_n(u) ) \xi(u)
\, du = \lim_n \int_0^1 ( f(u) + f_n(u) ) \xi(u) \, du
\]
is independent of itself! It means that the random variable is
degenerate (constant), which is impossible, since its distribution is
absolutely continuous.
\end{proof}

\section{The invariant via the logarithm}

Infinitesimality, as defined by Definition \ref{2.1}, will be used as
an invariant for discriminating \cpHss\ of type $ III $ that emerge by
exponentiation from continuous sums of FHS spaces. To this end,
infinitesimality w.r.t.\ such a \cpHs\ will be translated into the
language of FHS spaces.

\begin{definition}\label{11.1}
Let $ X $ be a metrizable topological space and $ X_1, X_2, \dotsc
\subset X $ closed subsets. We define 
$ \liminf_{n\to\infty} X_n $\index{liminf@$ \liminf $ for subspaces}
as
the set of limits of all convergent sequences $ x_1, x_2, \dots $ such
that $ x_1 \in X_1, x_2 \in X_2, \dots $
\end{definition}

The set $ \liminf X_n $ is always closed. If $ X $ is a linear
topological space and $ X_n $ are linear subspaces, then $ \liminf X_n
$ is a linear subspace. The following definition is formulated for
linear topological spaces, but we need it only for FHS spaces.

\begin{definition}
Let $ G $ be a linear topological space and $ G_1, G_2, \dots $ its
(closed linear) subspaces. We say that 
$ \limsup G_n = \{0\} $,\index{limsup@$ \limsup $ for subspaces}
if $ f
(x_n) \to 0 $ for every bounded sequence $ (x_n) $ such that $ x_n \in
G_n $, and every (continuous) linear functional $ f $ on $ G $.
\end{definition}

Let $ G $ be an FHS space and $ E_n, F_n $ its subspaces such that $ G
= E_n \oplus F_n $ (in the FHS sense) for each $ n $. We have the
tensor product of Hilbert spaces $ e^{E_n} \ot e^{F_n} = e^G $, and
two factors $ \A_n = \One \ot \B(e^{F_n}) $, $ \A'_n = \B(e^{E_n}) \ot
\One $ in the operator algebra $ \A = \B (e^G) $. These are real, but
may be complexified. Real-valued and complex-valued versions of
Condition (a) of the following theorem are equivalent, due to (an
evident generalization of) Lemma \ref{2.7}. We concentrate on the real
case until Lemma
\ref{11.7}.\index{real versus complex}\index{complex versus real}

\begin{theorem}\label{11.3}
The following two conditions are equivalent.

\textup{(a)} $ \sup_{A\in\A_n,\|A\|\le1} | \tr (AR) | \to 0 $ for $ n
\to \infty $ for each trace-class operator $ R \in \B(e^G) $
satisfying $ \tr (R) = 0 $;

\textup{(b)} $ \liminf E_n = G $ and $ \limsup F_n = \{0\} $.
\end{theorem}

The proof will be given after several lemmas and definitions.

\begin{sloppypar}
Let $ G $ be an FHS space and $ (X,\B,\M) $ its standard measure
class. Admissible norms on $ G $ are in a natural one-one
correspondence with Gaussian measures of $ \M $, $ \|\cdot\|
\leftrightarrow \ga_{\|\cdot\|} $. Moreover, the correspondence is a
homeomorphism, if both sets are equipped with appropriate
topologies.\index{convergence!for admissible norms}\index{convergence!for
Gaussian measures}
For any two admissible norms $ \|\cdot\|_1 $, $
\|\cdot\|_2 $ there exists a basis $ (e_k) $ of $ G $, orthogonal for
both norms, and eigenvalues $ \la_k = \| e_k \|_2 / \| e_k \|_1 $
satisfy $ \sum (\la_k^2-1)^2 < \infty $ (recall \ref{8.1}(a); $
1-\la_k^2 $ are eigenvalues of $ \One - L^* L $), or equivalently, $
\sum (\la_k-1)^2 < \infty $. A simple (basically, one-dimensional)
calculation gives
\begin{equation}\label{32}
\ip{ \sqrt{\ga_1} }{ \sqrt{\ga_2} } = \prod_k \bigg( \frac{
\la_k^{-1/2} + \la_k^{1/2} }{ 2 } \bigg)^{-1/2}
\end{equation}
where $ \ga_1, \ga_2 $ correspond to $ \|\cdot\|_1, \|\cdot\|_2 $. Of
course, convergence of the product is equivalent to convergence of the
series $ \sum_k (\la_k-1)^2 $. We may use, say, $ \arccos \ip{
\sqrt{\ga_1} }{ \sqrt{\ga_2} \, } $ as the distance between $ \ga_1 $ and
$ \ga_2 $, and $ \arccos \prod \( \frac{ \la_k^{-1/2} + \la_k^{1/2} }{
2 } \)^{-1/2} $ as the distance between $ \|\cdot\|_1 $ and $
\|\cdot\|_2 $,\footnote{%
 You can probably write a simpler equivalent metric.}
which makes the correspondence isometric. The distance
defines a topology on the set of all admissible norms. (No other
topology will be used on this set.) The distance $ \arccos \ip{
\sqrt{\ga_1} }{ \sqrt{\ga_2} \, } $ on the set of Gaussian measures of $
\M $ is equivalent to the variation distance $ \| \ga_1 - \ga_2 \| =
\int | \ga_1(d\om) - \ga_2(d\om) | $. The same holds for arbitrary
(not just Gaussian) probability measures of $ \M $; namely, $
2 - \sqrt{ 4 - \| \mu - \nu \|^2 } \le \| \sqrt\mu - \sqrt\nu \|^2 \le
\| \mu - \nu \| $. The variation
distance defines a topology on the set of all probability measures of
$ \M $. (No other topology will be used on this set.)
\end{sloppypar}

\begin{definition}\label{11.4}
(a) Let $ G $ be an FHS space, $ E_n, F_n \subset G $ subspaces, and $
\liminf E_n = G $. We say that $ F_n $ is asymptotically
orthogonal\index{asymptotically orthogonal}
to
$ E_n $, if there exists a convergent
sequence of admissible norms $ \|\cdot\|_n $ such that for each $ n $,
$ F_n $ is orthogonal to $ E_n $ w.r.t.\ $ \|\cdot\|_n $.

(b) Let $ (X,\B,\M) $ be a measure class, $ \cE_n, \F_n \subset
\B $ sub-\sif s, and $ \liminf \cE_n = \B $. We say that $ \F_n $ is
asymptotically
independent\index{asymptotically independent}
of $ \cE_n $, if there exists a
convergent sequence of probability measures $ \mu_n \in \M $ such that
for each $ n $, $ \F_n $ is independent of $ \cE_n $ w.r.t.\ $ \mu_n $.
\end{definition}

Here $ \liminf $ of \sif s is their $ \liminf $ (according to
\ref{11.1}) as closed subsets of the space of all measurable sets ($
\bmod \, 0 $) in $ (X,\B,\M) $. The topology of the space is defined
by the distance $ \mu ( A \bigtriangleup B ) $ between $ A $ and $ B
$; the metric depends on $ \mu \in \M $, but the topology does not.

\begin{remark}
In \ref{11.4}(a), `convergent' means, converging to \emph{some}
admissible norm $ \|\cdot\| $. It can be shown that the set of all
such norms $ \|\cdot\| = \lim_n \|\cdot\|_n $ is either the empty set,
or the whole space of admissible norms. In the latter case $ F_n $ is
asymptotically orthogonal to $ E_n $, in the former case it is
not. Similarly, under the conditions of \ref{11.4}(b), the set of all
possible measures $ \mu = \lim_n \mu_n $ is either the empty set, or
the whole $ \M $. However, we do not need it.
\end{remark}

\begin{theorem}\label{11.5}
Let $ G $ be an FHS space, $ E_n, F_n \subset G $ subspaces \textup{($
n = 1,2,\dotsc $)} such that $ \liminf E_n = G $, $ \limsup F_n =
\{0\} $, and $ G = E_n \oplus F_n $ \textup{(}in the FHS
sense\textup{)} for all $ n $. Then there exist subspaces $ C_n, D_n
\subset G $ such that $ E_n = C_n \oplus D_n $ \textup{(}in the FHS
sense\textup{)}, $ \liminf C_n = G $, and $ D_n \oplus F_n $ is
asymptotically orthogonal to $ C_n $.
\end{theorem}

The proof is postponed to the next section.

If $ (X,\B,\M) $ is a measure class and sub-\sif s $ \cE,\F \subset \B
$ satisfy $ \B = \cE \ot \F $, then $ (X,\B,\M) $ decomposes into the
product of two measure classes. These may be thought of as quotient
spaces $ (X,\B,\M) / \cE $ and $ (X,\B,\M) / \F $, or as $
(X,\cE,\M|_\cE) $ and $ (X,\F,\M|_\F) $ where $ \M|_\cE $ is the
equivalence class of measures on $ (X,\cE) $ containing restrictions $
\mu|_\cE $ of measures $ \mu \in \M $. Some $ \mu \in \M $ are product
measures satisfying the equality $ \mu = (\mu|_\cE) \ot (\mu|_\F)
$; others satisfy the equivalence $ \mu \sim (\mu|_\cE) \ot (\mu|_\F)
$. Also the Hilbert space decomposes into the tensor product, $ L_2
(X,\B,\M) = L_2(X,\cE,\M|_\cE) \ot L_2(X,\F,\M|_\F) $. Suppressing all
but \sif s we may write it shorter as $ L_2(\B) = L_2(\cE) \ot L_2(\F)
$. An element of $ L_2(\cE) $ may be written as $ f \sqrt{\mu|_\cE} $
where $ f $ is a $ \cE $-measurable function on $ X $, and $ \mu \in
\M $.

\begin{lemma}\label{11.6}
Let $ (X,\B,\M) $ be a measure class, $ \cE_n,\F_n \subset \B $
sub-\sif s \textup{($ n = 1,2,\dots $)}, $ \B = \cE_n \ot \F_n $ for
each $ n $, $ \liminf \cE_n = \B $, and $ \F_n $ asymptotically
independent of $ \cE_n $. Then there exists a probability measure $ \mu
\in \M $ such that
\[
\liminf_{n\to\infty} \( L_2 (\cE_n) \ot \sqrt{ \mu |_{\F_n} } \, \) =
L_2 (\B) \, .
\]
That is, for every $ \psi \in L_2 (\B) $ there exist $ \xi_n \in L_2
(\cE_n) $ such that $ \| \psi - \xi_n \ot \sqrt{ \mu |_{\F_n} } \|
\to 0 $ when $ n \to \infty $.
\end{lemma}

\begin{proof}
We take probability measures $ \mu, \mu_n \in \M $ such that $ \cE_n,
\F_n $ are $ \mu_n $-independent and $ \mu_n \to \mu $. We consider an
arbitrary finite set $ K \subset \R $, an arbitrary $ \B $-measurable
function $ f : X \to K $, and the corresponding vector $ \psi = f
\sqrt \mu \in L_2 (X,\B,\M) $. Using the fact that $ \liminf \cE_n = \B
$, we construct $ \cE_n $-measurable functions $ f_n : X \to K $ such
that $ f_n \to f $ in measure; then $ \| (f-f_n) \sqrt\mu \| \to 0 $.
Independence of $ \cE_n, \F_n $ w.r.t.\ $ \mu_n $ means that $
\sqrt{ \mu_n } = \sqrt{ \mu_n |_{\cE_n} } \ot \sqrt{ \mu_n |_{\F_n}
} $. Taking $ \xi_n = f_n \sqrt{ \mu_n |_{\cE_n} } \in L_2 (\cE_n) $
we get
\begin{multline*}
\| \psi - \xi_n \ot \sqrt{ \mu |_{\F_n} } \| = \| f \sqrt\mu - f_n
 \sqrt{ \mu_n |_{\cE_n} } \ot \sqrt{ \mu |_{\F_n} } \| \le \\
\le \| f \sqrt\mu - f_n \sqrt\mu \| + \| f_n \sqrt\mu - f_n
 \sqrt{\mu_n} \| + \| f_n \sqrt{\mu_n} - f_n  \sqrt{ \mu_n |_{\cE_n} }
 \ot \sqrt{ \mu |_{\F_n} } \| \le \\
\le \| (f-f_n) \sqrt\mu \| + \| f_n \|_\infty \| \sqrt\mu -
 \sqrt {\mu_n} \| + \| f_n  \sqrt{ \mu_n |_{\cE_n} } \ot \( \sqrt{
 \mu_n |_{\F_n} } - \sqrt{ \mu |_{\F_n} } \) \| \to 0 \, ;
\end{multline*}
indeed, $ \| f_n \|_\infty $ is bounded, $ \sqrt{\mu_n} \to
\sqrt\mu $, $ \| f_n  \sqrt{ \mu_n |_{\cE_n} } \| \le \| f_n
\|_\infty $, and $ \| \sqrt{ \mu_n |_{\F_n} } - \sqrt{ \mu |_{\F_n} }
\| \le \sqrt{ \| \mu_n |_{\F_n} - \mu |_{\F_n} \| } \le \sqrt{ \|
\mu_n - \mu \| } \to 0 $.
It remains to note that such vectors $ \psi $ (for all $ f $ and $ K
$) are dense in $ L_2 (X,\B,\M) $.
\end{proof}

If $ E $ is a subspace of an FHS space $ G $, then its standard
measure class is (canonically isomorphic to) a quotient
space of the standard measure class $ (X_G,\B_G,\M_G) $ of $ G $; in
other words, $ \B_E $ may be treated as a sub-\sif\ of $ \B_G $,
generated by $ \xi_x $ for $ x \in E $. Thus,
\[
E \subset G \quad \text{implies} \quad \B_E \subset \B_G \, .
\]
If $ E_n \subset G $ and $ \liminf E_n = G $, then $ \liminf \B_{E_n}
= \B_G $. Indeed, every $ \xi_x $ for $ x \in G $, being equal to $
\lim \xi_{x_n} $, $ x_n \in E_n $, is measurable w.r.t.\ the \sif\ $
\liminf \B_{E_n} $.

If $ E_n, F_n \subset G $, $ \liminf E_n = G $, and $ F_n $ is
asymptotically orthogonal to $ E_n $, then $ \B_{F_n} $ is
asymptotically independent of $ \B_{E_n} $. Indeed, we have $
\|\cdot\|_n \to \|\cdot\| $ and $ E_n \perp F_n $ w.r.t.\ $
\|\cdot\|_n  $; thus, Gaussian measures $ \ga, \ga_n \in \M $
corresponding to $ \|\cdot\|, \|\cdot\|_n $ satisfy $ \ga_n \to \ga $
and $ \B_{E_n} $, $ \B_{F_n} $ are independent w.r.t.\ $ \ga_n $.

\begin{proof}[Proof of theorem \ref{11.3}, part 1: (b) implies (a)]
Every trace-class operator $ R \in \B(e^G) $ satisfying $ \tr(R) = 0 $
is an (infinite) linear combination of two-dimensional
operators\footnote{%
 As usual, $ | \psi \rangle \langle \psi | $ stands for the
 one-dimensional projection $ x \mapsto \ip x \psi \psi $.}
$ R =
| \psi_1 \rangle \langle \psi_1 | - | \psi_2 \rangle \langle \psi_2 |
$, $ \psi_1, \psi_2 \in e^G $, $ \| \psi_1 \| = \| \psi_2 \| = 1 $. It
is enough to prove that
\[
\sup_{A\in\A_n,\|A\|\le1} | \ip{A\psi_1}{\psi_1} -
\ip{A\psi_2}{\psi_2} | \to 0 \quad \text{for } n \to \infty \, .
\]
Theorem \ref{11.5} (to be proven later) gives us subspaces $ C_n, D_n
\subset G $ such that $ E_n = C_n \oplus D_n $ (in the FHS sense),
$ \liminf C_n = G $, and $ D_n \oplus F_n $ is asymptotically
orthogonal to $ C_n $. The corresponding sub-\sif s satisfy $ \liminf
\B_{C_n} = \B_G $, and $ \B_{D_n\oplus F_n} $ is asymptotically
independent of $ \B_{C_n} $; also, $ \B_{C_n} \ot \B_{D_n\oplus F_n} =
\B_G $ (since $ C_n \oplus D_n \oplus F_n = G $ in the FHS sense).
Lemma \ref{11.6} gives us a representation
\[
\psi = \lim_{n\to\infty} ( \xi_n \ot \chi_n )
\]
for an arbitrary unit vector $ \psi \in e^G = L_2 (\B_G) $; here $
\xi_n \in e^{C_n} = L_2 (\B_{C_n}) $, $ \| \xi_n \| = 1 $, $ \chi_n \in
e^{D_n\oplus F_n} = L_2 (\B_{D_n\oplus F_n}) $, $ \| \chi_n \| = 1 $,
and these $ \chi_n $ (unlike $ \xi_n $) do not depend on $ \psi $.

Let $ A \in \A_n = \One \ot \B (e^{F_n}) \subset \B (e^G) $, $ \| A \|
\le 1 $. We have $ e^G = e^{C_n} \ot e^{D_n} \ot e^{F_n} = e^{C_n} \ot
e^{D_n\oplus F_n} $, $ A \in \One \ot \B (e^{D_n\oplus F_n}) $; $ A =
\One \ot \ti A $, $ \ti A \in \B (e^{D_n\oplus F_n}) $, $ \| \ti A \|
\le 1 $. Thus, $ \ip{ A(\xi_n \ot \chi_n) }{ \xi_n \ot \chi_n } = \ip{
\ti A \chi_n }{ \chi_n } $, which does not depend on $ \psi $. It
remains to note that $ | \ip{ A\psi }{ \psi } - \ip{ A(\xi_n \ot
\chi_n) }{ \xi_n \ot \chi_n } | \le \| A \| \cdot \| | \psi \rangle
\langle \psi | - | \xi_n \ot \chi_n \rangle \langle \xi_n \ot \chi_n |
\| \le 2 \| \psi - \xi_n \ot \chi_n \| \to 0 $.
\end{proof}

The first part is proven for the real case, but the second part will
be proven for the complex case (that is, complex Hilbert spaces out of
real FHS spaces); as was said, the two cases are
equivalent.\index{real versus complex}\index{complex versus real}

Let $ G = E \oplus F $ be the direct sum of two FHS spaces, $
\|\cdot\| $ an admissible norm on $ G $ (note that $ E, F $ need not
be orthogonal w.r.t.\ $ \|\cdot\| $), and $ \ga \in \M_G $ the
Gaussian measure on $ X_G $ corresponding to $ \|\cdot\| $. Recall
operators $ U_x, V_y $ ($ x \in G, y \in G' $) on the Hilbert space $
e_\C^G $ and their properties stated in \ref{8.5}(c). We treat $ E,F $
as subspaces of $ G $, and their duals $ E', F' $ as subspaces of $ G'
$ (according to $ G = E \oplus F $, irrespective of $ \|\cdot\| $). In
the algebra $ \B (e_\C^G) = \B ( e_\C^E \ot e_\C^F ) $ we introduce a
factor $ \A = \One \ot \B (e_\C^F) $. (Note that the vector $ \sqrt\ga
\in e_\C^E \ot e_\C^F $ need not be a product vector.)

\begin{lemma}\label{11.7}
There exists an increasing function $ M : [0,\infty) \to [0,\infty) $
such that $ M(r) > 0 $ for all $ r > 0 $, and for all $ G, E, F, \A,
\|\cdot\|, \ga $ as above,
\begin{align*}
\sup_{A\in\A,\|A\|\le1} | \ip{ A\sqrt\ga }{ \sqrt\ga } - \ip{ A V_y
 \sqrt\ga }{ V_y \sqrt\ga } | &\ge M \( \dist (y,E') \) \, , \\
\sup_{A\in\A,\|A\|\le1} | \ip{ A\sqrt\ga }{ \sqrt\ga } - \ip{ A U_x
 \sqrt\ga }{ U_x \sqrt\ga } | &\ge M \( \dist (x,E) \)
\end{align*}
for all $ x \in G $, $ y \in G' $; here the distance $ \dist(x,E) =
\inf_{z\in E} \| x - z \| $ is taken w.r.t.\ the norm $ \|\cdot\| $,
and $ \dist(y,E') $ --- w.r.t.\ the dual norm on $ G' $.
\end{lemma}

\begin{proof}
The operator $ A = U_x $ for $ x \in F $ belongs to $ \A $, and $ \| A
\| \le 1 $. We have
\[
\ip{ A\sqrt\ga }{ \sqrt\ga } = \exp \bigg( \! -\frac12 \|x\|^2 \bigg) \,
;
\]
\begin{multline*}
\ip{ A V_y \sqrt\ga }{ V_y \sqrt\ga } = \ip{ U_x V_y \sqrt\ga }{ V_y
 \sqrt\ga } = e^{-i\ip x y} \ip{ V_y U_x \sqrt\ga }{ V_y \sqrt\ga } = \\
= e^{-i\ip x y} \ip{ U_x \sqrt\ga }{ \sqrt\ga } = \exp \bigg( \! -i\ip
 x y - \frac12 \|x\|^2 \bigg) \, ;
\end{multline*}
\begin{multline*}
| \ip{ A\sqrt\ga }{ \sqrt\ga } - \ip{ A V_y \sqrt\ga }{ V_y \sqrt\ga }
 | = | 1 - e^{-i\ip x y} | \exp \bigg( \! -\frac12 \|x\|^2 \bigg) = \\
= 2 \Big| \sin \frac{ \ip x y }{ 2 } \Big| \exp \bigg( \! -\frac12
\|x\|^2 \bigg) \, .
\end{multline*}
Taking into account that $ \dist (y,E') = \sup_{x\in F,\|x\|\le1} |
\ip x y | $ we get the first inequality for
\[
M(r) = \sup_{u\ge0} 2\sin \Big( \frac12 ru \Big) \exp \Big( -\frac12
u^2 \Big) = \sup_{\phi\ge0} \bigg( \exp \Big( -\frac{\phi^2}{2r^2}
\Big) \cdot 2 \Big| \sin \frac\phi2 \Big| \bigg) \, .
\]
The proof of the second inequality is quite similar. Only $
U_x, V_y $ are interchanged, and $ \exp \( -\frac18 \|y\|^2 \) $
appears instead of $ \exp \( -\frac12 \|x\|^2 \) $, which leads to $
M(2r) $ instead of $ M(r) $.
\end{proof}

\begin{proof}[Proof of theorem \ref{11.3}, part 2: (a) implies (b)]
We choose a Gaussian measure $ \ga \in \M_G $ and apply (a) to the
operator $ R = | \sqrt\ga \rangle \langle \sqrt\ga | - | V_y \sqrt\ga
\rangle \langle V_y \sqrt\ga | $ for an arbitrary $ y \in G' $, and
also to $ R = | \sqrt\ga \rangle \langle \sqrt\ga | - | U_x \sqrt\ga
\rangle \langle U_x \sqrt\ga | $ for an arbitrary $ x \in G $:
\begin{align*}
\sup_{A\in\A,\|A\|\le1} | \ip{ A\sqrt\ga }{ \sqrt\ga } - \ip{ A V_y
 \sqrt\ga }{ V_y \sqrt\ga } | &\to 0 \, , \\
\sup_{A\in\A,\|A\|\le1} | \ip{ A\sqrt\ga }{ \sqrt\ga } - \ip{ A U_x
 \sqrt\ga }{ U_x \sqrt\ga } | &\to 0
\end{align*}
for $ n \to 0 $. Lemma \ref{11.7} gives
\[
\dist (y,E'_n) \to 0 \, , \quad \dist (x,E_n) \to 0 \, .
\]
The latter shows that $ x \in \liminf E_n $, which means that $
\liminf E_n = G $. The former shows that
\[
\sup_{z\in F_n,\|z\|\le1} | \ip y z | \to 0 \quad \text{for } n \to 0
\, ,
\]
which means that $ \limsup F_n = \{0\} $.
\end{proof}

So, Theorem \ref{11.3} follows from Theorem \ref{11.5}.

\section{Ensuring asymptotic orthogonality}

The sole goal of this section is to prove Theorem \ref{11.5}.

\begin{lemma}\label{12.1}
Let $ G $ be an FHS space, and $ F_n \subset G $ subspaces. Then the
following conditions are equivalent.

\textup{(a)} $ \limsup F_n = \{0\} $;

\textup{(b)} for every finite-dimensional subspaces $ E_1, E_2, \dotsc
$ such
that $ \liminf E_n = G $ there exist integers $ k_1 \le k_2 \le \dotsc
$ such that $ k_n \to \infty $ and $ F_n $ is asymptotically
orthogonal to $ E_{k_n} $.
\end{lemma}

\begin{proof}
The implication ``(b) $\Rightarrow$ (a)'' is simpler, and will not be
used; I leave it to the reader. Assume (a).
We choose an admissible norm on $ G $, thus turning $ G $ into a
Hilbert space. Condition (a) gives
\[
\forall x \in G \quad \sup_{f\in F_n, \|f\|\le1} \ip f x
\xrightarrow[n\to\infty]{} 0 \, ,
\]
that is,
\[
\forall x \in G \quad \angle (x,F_n) \xrightarrow[n\to\infty]{}
\frac\pi2 \, ,
\]
where the angle is defined by $ \cos \angle (x,F_n) = \sup \{ \ip x f
: f\in F_n, \|f\|\le1 \} $. Therefore
\[
\forall E \quad \angle (E,F_n) \xrightarrow[n\to\infty]{} \frac\pi2 \, ,
\]
where $ E $ runs over finite-dimensional subspaces, and $ \cos \angle
(E,F_n) = \sup \{ \ip e f : e \in E, f \in F_n, \|e\| \le 1, \|f\| \le
1 \} $. The following lemma completes the proof, provided that $ k_n $
tends to $ \infty $ slowly enough. Namely, in terms of $ \de
(\cdot,\cdot) $ introduced there, it suffices that $ \de \( \angle
(E_{k_n},F_n), \dim E_{k_n} \) \xrightarrow[n\to\infty]{} 0 $.
\end{proof}

\begin{lemma}\label{12.2}
Let $ G $ be an FHS-space, $ E, F \subset G $ subspaces, $ \dim (E) <
\infty $, $ E \cap F = \{ 0 \} $. Then for every admissible norm $ \|
\cdot \|_1 $ there exists an admissible norm $ \|\cdot\|_2 $ such that
$ E,F $ are orthogonal w.r.t.\ $ \|\cdot\|_2 $, and
\[
\dist ( \|\cdot\|_1, \|\cdot\|_2 ) \le \de \( \angle(E,F), \dim E \)
\]
for some function $ \de : [0,\frac\pi2] \times \{0,1,2,\dots\} \to
(0,\infty) $ such that for every $ n $, $ \de(\al,n) \to 0 $ for $
\al \to \frac\pi2 $.\footnote{%
 Of course, $ \de $ does not depend on $ G,E,F $.}
\end{lemma}

\begin{proof}
We equip $ G $ with the norm $ \|\cdot\|_1 $, thus turning $ G $ into
a Hilbert space, and consider orthogonal projections $ Q_E, Q_F $ onto
$ E,F $ respectively. Introduce subspaces $ E \cap
F^\perp $, $ E^\perp \cap F $, $ E^\perp \cap F^\perp $ (here $
E^\perp $ is the orthogonal complement of $ E $); the subspaces are
orthogonal to each other, and invariant under both $ Q_E $ and $ Q_F
$. Therefore
\[
G = G_0 \oplus (E\cap F^\perp) \oplus (E^\perp\cap F)
\oplus (E^\perp\cap F^\perp) \, , \quad \text{(orthogonal sum)}
\]
where $ G_0 $ is another subspace invariant under $ Q_E, Q_F $ (since
these operators are Hermitian). Introduce $ E_0 = E \cap G_0 $, $ F_0
= F \cap G_0 $, then $ Q_E h_0 = Q_{E_0} h_0 $ for all $ h_0 \in G_0 $
(since $ Q_E $ commutes with $ Q_{G_0} $), and $ Q_F h_0 = Q_{F_0} h_0
$. We may get rid of $ G_0^\perp $ by letting
\[
\| h_0 + h_1 \|_2^2 = \| h_0 \|_2^2 + \| h_1 \|_1^2 \quad \text{for
all $ h_0 \in G_0 $, $ h_1 \in G_0^\perp $.}
\]
In other words, we'll construct $ \|\cdot\|_2 $ on $ G_0 $ while
preserving both the given norm on $ G_0^\perp $ and the orthogonality
of $ G_0, G_0^\perp $. Now we forget about $ G_0^\perp $, assuming
that $ G = G_0 $, $ E = E_0 $, $ F = F_0 $.

So, we have $ E \cap F = \{ 0 \} $, $ E \cap F^\perp = \{ 0 \} $, $
E^\perp \cap F = \{ 0 \} $, $ E^\perp \cap F^\perp = \{ 0 \} $. The
latter implies $ \dim (F^\perp) \le \codim (E^\perp) = \dim E
$. Similarly, $ \dim F \le \dim E $. Therefore $ G $ is
finite-dimensional, $ \dim G \le 2 \dim E $.

Both $ Q_E $ and $ Q_F $ commute with the Hermitian operator $ C =
\frac12 (2Q_E-1) (2Q_F-1) + \frac12 (2Q_F-1) (2Q_E-1) $. The spectrum
of $ C $ consists of some numbers $ \cos 2\phi_k $ of multiplicity 2
(though, some $ \phi_k $ may coincide), and $ 0 < \phi_k < \frac\pi2 $
(the case $ \phi_k = 0 $ is excluded by $ E \cap F = \{ 0 \} $; the
case $ \phi_k = \pi/2 $ is excluded by $  E \cap F^\perp = \{ 0 \} $,
$ E^\perp \cap F = \{ 0 \} $, $ E^\perp \cap F^\perp = \{ 0 \}
$). Accordingly, $ G $ decomposes into the (orthogonal) direct sum
of planes, $ G = G_1 \oplus \dots \oplus G_d $, $ \dim G_k = 2 $,
invariant under $ Q_E, Q_F $. Subspaces $ E_k = E \cap G_k $, $ F_k =
F \cap G_k $ are two lines on the plane $ G_k $, and $ \angle
(E_k,F_k) = \phi_k $; $ k=1,\dots,d $; $ d \le \dim E $. Clearly,
\[
\angle (E,F) = \min ( \phi_1, \dots, \phi_d ) \, .
\]
We construct $ \|\cdot\|_2 $ on each $ G_k $ separately, while
preserving their orthogonality. Elementary 2-dimensional geometry
shows that the corresponding numbers $ \la'_k, \la''_k $ (mentioned in
\eqref{32}; two numbers
for each plane) are, in the optimal case,
\[
\la'_k = \( \tan \frac{\phi_k}2 \)^{-1/2} \, , \quad \la''_k = \( \tan
\frac{\phi_k}2 \)^{1/2} \, .
\]
The corresponding right-hand side of \eqref{32} is
\[
\cos \dist ( \|\cdot\|_1, \|\cdot\|_2 ) = \prod_{k=1}^d \bigg( \frac{
\tan^{-1/4} \frac{\phi_k}2 + \tan^{1/4} \frac{\phi_k}2 }{ 2 }
\bigg)^{-1} \, ,
\]
which gives the needed result for
\[
\de(\al,n) = \arccos \bigg( \frac{ \tan^{-1/4} \frac\al2 +
\tan^{1/4} \frac\al2 }{ 2 } \bigg)^{-n} \, .
\]
\end{proof}

\begin{proof}[Proof of theorem \ref{11.5}]
\begin{sloppypar}
We choose an admissible norm on $ G $, thus turning $ G $ into a
Hilbert space over $ \R $. Let $ L \subset G $ be a finite-dimensional
subspace, $
L \ne \{0\} $. For any given $ n $ consider the pair $ L, E_n $. Its
geometry may be described (similarly to the proof of Lemma \ref{12.2})
via angles $ \phi_1^{(n)}, \dots, \phi_{d_n}^{(n)} \in [0,\frac\pi2)
$, $ d_n \le \dim L $. This time, zero angles are allowed, since $ L
\cap E_n $ need not be $ \{0\} $. It may happen that $ d_n < \dim L $,
since $ L \cap E_n^\perp $ need not be $ \{0\} $. However,
\[
\sup_{x\in L, x\ne0} \angle (x,E_n) = \al_n \to 0 \quad \text{for } n
\to \infty \, ,
\]
since $ \liminf E_n = G $; for large $ n $ we have $ \al_n < \pi/2 $
which implies $ d_n = d =
\dim L $ and $ \max ( \phi_1^{(n)}, \dots, \phi_d^{(n)} ) = \al_n
$. We may send $ L $ into $ E_n $ rotating it by $ \phi_1^{(n)},
\dots, \phi_d^{(n)} $. In other words, there is a rotation $ U_n : G
\to G $ such that
\[
U_n (L) \subset E_n \quad \text{and} \quad \| U_n - 1 \| \le 2 \sin
\frac{ \al_n }{ 2 } \xrightarrow[n\to\infty]{} 0 \, .
\]
\end{sloppypar}

We choose subspaces $ L_k \subset G $ such that $ \dim L_k = k $ and $
\liminf L_k = G $.\footnote{%
 Of course, one may take $ L_1 \subset L_2 \subset \dotsb $}
Introduce
\[
\al_{k,n} = \sup_{x\in L_k, x\ne0} \angle (x,E_n) \, ,
\]
then $ \al_{k,n} \xrightarrow[n\to\infty]{} 0 $ for each $ k $. On the
other hand, introduce
\[
\be_{k,n} = \frac\pi2 - \angle (L_k,F_n) \, .
\]
Similarly to the proof of Lemma \ref{12.1} we have $ \be_{k,n}
\xrightarrow[n\to\infty]{} 0 $ for each $ k $; therefore\footnote{%
 Recall that $ \de(\cdot,\cdot) $ is introduced in Lemma \ref{12.2}.}
$ \de ( \frac\pi2 - \be_{k_n,n}, k_n ) \xrightarrow[n\to\infty]{} 0 $
if $ k_n $ tends to $ \infty $ slowly enough. However, we choose $ k_1
\le k_2 \le \dotsb $, $ k_n \to \infty $ so as to satisfy a stronger
condition:
\[
\de \Big( \frac\pi2 - \al_{k_n,n} - \be_{k_n,n}, \, k_n \Big)
\xrightarrow[n\to\infty]{} 0 \, .
\]
We take
\[
C_n = U_n ( L_{k_n} ) \, ,
\]
where rotations $ U_n $ satisfy $ U_n (L_{k_n}) \subset E_n $ and $ \|
U_n - 1 \| \le 2 \sin (\frac12 \al_{{k_n},n} ) \to 0 $. Then $ \liminf
C_n = G $, and
\[
\frac\pi2 - \angle (C_n,F_n) \le \al_{k_n,n} + \be_{k_n,n} \, ;
\]
due to Lemma \ref{12.2}, there exist admissible norms $ \|\cdot\|_n
\to \|\cdot\| $ such that $ F_n $ is orthogonal to $ C_n $ w.r.t.\ $
\|\cdot\|_n $. Consider the
orthogonal complement $ M_n $ of $ C_n $ w.r.t.\ $ \|\cdot\|_n $;
clearly, $ M_n $ is asymptotically orthogonal to $ C_n $. We have $
F_n \subset M_n $ and $ G = C_n \oplus M_n $ (in the FHS-sense). On
the other hand, $ C_n \subset E_n $ and $ G = E_n \oplus F_n
$. Applying Remark \ref{9.05} (for each $ n $ separately) to $ G_1 =
C_n $, $ G_{1,2} = E_n $, $ G_{2,3} = M_n $, $ G_3 = F_n $ we see that
the subspace
\[
D_n = E_n \cap M_n
\]
satisfies $ C_n \oplus D_n \oplus F_n = G $ and $ C_n \oplus D_n = E_n
$ (and also $ D_n \oplus F_n = M_n $). Asymptotical orthogonality of $
D_n \oplus F_n $ to $ C_n $ follows from their orthogonality w.r.t.\ $
\|\cdot\|_n $.
\end{proof}

\section{Calculating the invariant}

\smallskip

\hfill\parbox{9cm}{%
It is conjectured that there are uncountably many type $ III $ $ E_0
$-semigroups which are not cocycle conjugate. (So far we have only one
example!)\\
\mbox{}\hfill W.~Arveson \cite[p.~167]{Ar99}
}

\medskip

As was said in Sect.~2, we consider sequences $ (E_n) $ of elementary
sets
\begin{equation}\label{33}
E_n = \bigcup_{k=0}^{n-1} \bigg( \frac1n \Big( k + \frac{1-\eps_n}2
\Big), \, \frac1n \Big( k + \frac{1+\eps_n}2 \Big) \bigg) \, ;
\end{equation}
infinitesimality of $ (E_n) $ depends on $ (\eps_n) $ and the spectral
measure $ \nu $ (assumed to satisfy \ref{9.5}(b) and
\eqref{29}). Theorem \ref{11.3} states that $ (E_n) $ is infinitesimal
if and only if
\begin{gather}
\liminf G_{(0,1)\setminus E_n} = G_{(0,1)} \, , \label{34} \\
\limsup G_{E_n} = \{ 0 \} \, ; \label{35}
\end{gather}
here $ G_{E_n} $ is the closure in $ L_2 (\nu) $ of the set of Fourier
transforms $ \hat f $ of functions $ f \in L_2(E_n) $ (or rather, $ f
\in L_2(\R) $ vanishing outside $ E_n $).

\begin{lemma}\label{13.1}
If $ \eps_n \to 0 $ then \eqref{34} holds.
\end{lemma}

\begin{proof}
For every $ f \in L_2(0,1) $,
\[
\| f - f \cdot \One_{E_n} \|_\nu \le \const \cdot \, \| f - f \cdot
\One_{E_n} \|_{L_2(0,1)} \to 0 \, .
\]
\end{proof}

\begin{remark}
If $ \nu $ satisfies \eqref{26} then \eqref{34} holds already for the
constant sequence $ \eps_n = 1/2 $, due to (a slight modification of)
Lemma \ref{10.1}.
\end{remark}

The question is when \eqref{35} holds. Here is a sufficient
condition for arbitrary elementary sets $ E_n \subset (0,1) $, not
just of the form \eqref{33}.

\begin{lemma}\label{13.3}
Existence of a function $ \be : (0,\infty) \to (0,\infty) $ and
numbers $ \de_n \to 0 $ satisfying conditions \textup{(a)---(c)} below
is sufficient for \eqref{35}:

\textup{(a)} $ \be $ decreases, but $ \la \mapsto \la^2 \be(\la) $
increases;

\textup{(b)} $ \nu(d\la) \ge \be(\la) \, d\la $;

\textup{(c)} $\displaystyle \frac{ \mes \( (E_n)_{+\de_n} \) }{
\be(\de_n^{-1}) } \to 0 $ for $ n \to \infty $, where $ (E_n)_{+\de_n}
$ stands for the $ \de_n $-neighborhood of $ E_n $, and \textup{`$ \mes $'} for
the Lebesgue measure.
\end{lemma}

\begin{proof}
Let $ f_n \in L_2(E_n) $, $ \| f_n \|_\nu \le 1 $. We have to prove
that $ \hat f_n \to 0 $ weakly in $ L_2(\nu) $, that is, $ \int \hat
f_n(\la) \overline{ \hat g(\la) }
\, d\la \to 0 $ for every $ \hat g $ such that $ \int | \hat g(\la)
|^2 \frac{d\la}{\nu(d\la)} \, d\la < \infty $.\footnote{%
 For now $ \hat g $ is not assumed to be a Fourier transform of some $
 g $, but it will happen soon.}
We may restrict ourselves to a dense set of such $ \hat g $; thus we
assume that $ \hat g $ has a compact support. It follows that $ \int
| \hat g(\la) |^2 \, d\la < \infty $ (since $ \be $ is bounded away
from $ 0 $ on the support of $ \hat g $) and $ \int | \hat g(\la) | \,
d\la < \infty $. Therefore $ \hat g $ is the Fourier transform of a
function $ g \in L_2(\R) \cap L_\infty (\R) $.

For every $ \de \in (0,\infty) $,
\[
\frac{ \nu(d\la) }{ d\la } \ge \be(\de^{-1}) \bigg( \frac{ \sin \de\la
}{ \de\la } \bigg)^2 \quad \text{for all } \la \in (0,\infty) \, ;
\]
indeed, for $ \la \in (0,\de^{-1}) $ we have
\[
\bigg( \frac{ \de\la }{ \sin \de\la } \bigg)^2 \frac{ \nu(d\la) }{
d\la } \ge \frac{ \nu(d\la) }{ d\la } \ge \be(\la) \ge \be(\de^{-1})
\, ,
\]
and for $ \la \in (\de^{-1},\infty) $ we have
\[
\bigg( \frac{ \de\la }{ \sin \de\la } \bigg)^2 \frac{ \nu(d\la) }{
d\la } \ge (\de\la)^2 \frac{ \nu(d\la) }{ d\la } \ge \de^2 \la^2
\be(\la) \ge \de^2 \de^{-2} \be(\de^{-1}) \, .
\]
We define functions $ h_n $ by
\[
\hat h_n (\la) = \hat f_n (\la) \frac{ \sin \de_n \la }{ \de_n \la }
\, ;
\]
then $ h_n \in L_2 \( (E_n)_{+\de_n} \) $, since\footnote{%
 A constant appears, since we normalize Fourier transform so as to be
 unitary.}
\[
h_n(t) = \frac{ \const }{ 2\de_n } \int_{t-\de_n}^{t+\de_n} f_n(s) \,
ds \, .
\]
We have
\begin{multline*}
\int h_n^2 (t) \, dt = \int |\hat h_n(\la)|^2 \, d\la = \int |\hat
 f_n(\la)|^2 \bigg( \frac{ \sin \de_n \la }{ \de_n \la } \bigg)^2 \,
 d\la \le \\
\le \frac2{ \be(\de_n^{-1}) } \int_0^\infty |\hat f_n(\la)|^2 \,
 \nu(d\la) \, ;
\end{multline*}
\[
\be(\de_n^{-1}) \| h_n \|^2_{L_2} \le 2 \| f_n \|_\nu^2 \le 2 \, ; \quad
\sqrt{ \be(\de_n^{-1}) } \| h_n \|_{L_2} \le \sqrt 2 \, .
\]
However,
\begin{multline*}
\bigg| \int \hat h_n(\la) \overline{ \hat g(\la) } \, d\la \bigg| = \bigg| \int
 h_n(t) g(t) \, dt \bigg| = \\
= \bigg| \int h_n(t) \One_{(E_n)_{+\de_n}} (t) g(t) \, dt \bigg| \le
 \| h_n \|_{L_2} \| \One_{(E_n)_{+\de_n}} g \|_{L_2} \le \\
\le \frac{\sqrt 2 }{ \sqrt{ \be(\de_n^{-1}) } } \sqrt{ \mes
 (E_n)_{+\de_n} } \, \| g \|_\infty \to 0 \, .
\end{multline*}
It remains to note that
\begin{multline*}
\bigg| \int \hat f_n(\la) \overline{ \hat g(\la) } \, d\la - \int \hat h_n(\la)
\overline{ \hat g(\la) } \, d\la \bigg| \le \int | \hat f_n(\la) |
\bigg( 1 - \frac{
  \sin \de_n \la }{ \de_n \la } \bigg) | \hat g(\la) | \, d\la
  \le \\
\Big( \int |\hat f_n(\la)|^2 \, \nu(d\la) \Big)^{1/2} \Big( \int
 |\hat g(\la)|^2 \frac{ d\la }{ \nu(d\la) } \, d\la \Big)^{1/2}
 \max_{\la\in\supp\hat g} \bigg( 1 - \frac{ \sin \de_n \la }{
 \de_n \la } \bigg) \to 0 \, .
\end{multline*}
\end{proof}

Now we need a necessary condition for \eqref{35}. Sets $ E_n $ are
again assumed to be of the form \eqref{33}. We introduce
functions
\[
f_n = \frac1{\eps_n} \One_{E_n} \, , \quad f = \One_{(0,1)} \, .
\]
If $ \| f_n - f \|_\nu \to 0 $ then \eqref{35} is evidently
violated. Lemma \ref{10.1} gives us $ \| f_n - f \|_\nu \to 0 $ for
the constant sequence $ \eps_n = 1/2 $ (assuming \eqref{26}); now it
will be proven for
$ \eps_n $ tending to $ 0 $ but not too fast.

Recall that a slowly
varying\index{slowly varying}
function is a function $ \al : (0,\infty)
\to (0,\infty) $ such that
\[
\al (\la) = (1+o(1)) \exp \bigg( \int_1^\la \frac{ \eta(\la_1) }{
\la_1 } \, d\la_1 \bigg) \quad \text{for } \la \to \infty
\]
for some function $ \eta(\cdot) $ satisfying $ \eta(\la) \to 0 $ for $
\la \to \infty $; see \cite[Sect.~VIII.9]{Flr}. Observe that functions
\eqref{293}, \eqref{294} are slowly varying.

\begin{lemma}\label{13.4}
Let $ \al : (0,\infty) \to (0,\infty) $ be a slowly varying function
such that $ \nu(d\la) \le \al(\la) \, d\la $, and $ \eps_n \to 0 $
satisfy
\[
\frac1{\eps_n} \al \bigg( \frac{n}{\eps_n} \bigg) \to 0 \quad
\text{for } n \to \infty \, .
\]
Then $ \| f_n - f \|_\nu \to 0 $ for $ n \to \infty $.
\end{lemma}

\begin{proof}
An elementary calculation gives
\[
| \hat f_n (\la) | = \frac{\const}{\eps_n} \left| \frac{ \displaystyle
\sin \frac{\la\eps_n}n \sin \frac\la2 }{ \displaystyle \la \sin
\frac{\la}{2n} } \right| \, .
\]
We have to prove that $ \int | \hat f_n (\la) - \hat f(\la) |^2 \,
\nu(d\la) \to 0 $. First of all, $ \hat f_n (\la) \to \hat f(\la) $
uniformly on $ \la \in (0,M) $ for every $ M < \infty $. Therefore it
is enough to prove that $ \int_M^\infty | \hat f_n (\la) |^2 \,
\nu(d\la) \to 0 $ for $ M \to \infty $, $ n \to \infty $, $ \pi n > M $. We have
\begin{multline*}
\int_M^\infty | \hat f_n (\la) |^2 \al(\la) \, d\la \le
\frac{\const}{\eps_n^2} \int_M^\infty \al(\la) \frac{ \sin^2
 \frac{\la\eps_n}n }{ \la^2 } \frac{ \sin^2 \frac\la2 }{ \sin^2
 \frac\la{2n} } \, d\la = \\
= \frac{\const}{\eps_n^2} \bigg( \int_M^{\pi n} \dotsc \, d\la +
 \sum_{k=1}^\infty \int_{(2k-1)\pi n}^{(2k+1)\pi n} \dotsc \, d\la
 \bigg) \, .
\end{multline*}
The function $ \frac{ \sin^2 \frac\la2 }{ \sin^2 \frac\la{2n} } =
\big| \sum_{k=0}^{n-1} \exp (ik\la/n) \big|^2 $ has period $ 2\pi n $,
and its mean value (over the period) is equal to $ n $. On the other
hand, $ \frac1{\la^2} \sin^2 \frac{\la\eps_n}n \le \min \(
\frac{\eps_n^2}{n^2}, \frac1{\la^2} \) $. Thus, for $ n \to \infty $,
\begin{multline*}
\sum_{k=1}^\infty \int_{(2k-1)\pi n}^{(2k+1)\pi n} \al(\la) \frac{
 \sin^2 \frac{\la\eps_n}n }{ \la^2 } \frac{ \sin^2 \frac\la2 }{ \sin^2
 \frac\la{2n} } \, d\la \le \\
\le \sum_{k=1}^\infty \bigg( \max_{\la\in[(2k-1)\pi n,(2k+1)\pi n]}
 \al(\la) \bigg) \min \bigg( \frac{\eps_n^2}{n^2}, \frac1{(2k-1)^2
 \pi^2 n^2} \bigg) \cdot 2\pi n \cdot n \le \\
\le (1+o(1)) n \int_{\pi n}^\infty \al(\la) \min \Big(
 \frac{\eps_n^2}{n^2}, \frac1{\la^2} \Big) \, d\la = \\
= (1+o(1)) n \bigg( \frac{\eps_n^2}{n^2} \int_{\pi n}^{n/\eps_n}
 \al(\la) \, d\la + \int_{n/\eps_n}^\infty \al(\la) \la^{-2} \, d\la
 \bigg) \, .
\end{multline*}
However, $ \int_{\pi n}^{n/\eps_n} \al(\la) \, d\la = (1+o(1))
\frac{n}{\eps_n} \al \( \frac{n}{\eps_n} \) $; see \cite[Th.~1 of
Sect.~VIII.9]{Flr}. Similarly, $ \int_{n/\eps_n}^\infty \al(\la)
\la^{-2} \, d\la = (1+o(1)) \al \( \frac{n}{\eps_n} \)
\int_{n/\eps_n}^\infty \la^{-2} \, d\la = (1+o(1)) \frac{\eps_n}{n}
\al \( \frac{n}{\eps_n} \) $. So,
\[
\frac1{\eps_n^2} \sum_{k=1}^\infty \int_{(2k-1)\pi n}^{(2k+1)\pi n}
\dotsc \, d\la \le (1+o(1)) \frac{n}{\eps_n^2} \cdot 2 \frac{\eps_n}{n}
\al \Big( \frac{n}{\eps_n} \Big) = (1+o(1)) \frac{2}{\eps_n} \al \Big(
\frac{n}{\eps_n} \Big) \to 0 \, .
\]
It remains to prove that
\[
\frac1{\eps_n^2} \int_M^{\pi n} \al(\la) \frac{ \sin^2
 \frac{\la\eps_n}n }{ \la^2 } \frac{ \sin^2 \frac\la2 }{ \sin^2
 \frac\la{2n} } \, d\la \to 0 \quad \text{for } M \to \infty, \, n \to
 \infty, \, \pi n > M \, .
\]
We have $ \sin \frac\la{2n} \ge
\frac2\pi \frac\la{2n} $, $ \sin^2 \frac\la2 \le 1 $, $ \sin^2
\frac{\la\eps_n}n \le \( \frac{\la\eps_n}n \)^2 $, $ \al(\la) \le \|
\al \|_\infty $; thus
\begin{multline*}
\frac1{\eps_n^2} \int_M^{\pi n} \dotsc \, d\la \le \frac1{\eps_n^2} \|
 \al \|_\infty \Big( \frac\pi2 \Big)^2 \int_M^\infty \frac{ \(
 \frac{\la\eps_n}n \)^2 }{ \la^2 \( \frac\la{2n} \)^2 } \, d\la
 = \\
= \frac1{\eps_n^2} \| \al \|_\infty \Big( \frac\pi2 \Big)^2
 (2\eps_n)^2 \int_M^\infty \frac{d\la}{\la^2} = \pi^2 \| \al \|_\infty
 \frac1M \, .
\end{multline*}
\end{proof}

Having a sufficient condition (given by Lemma \ref{13.3}) and a
necessary condition (given by Lemma \ref{13.4}), we want to compare
them. If $ E_n $ in Lemma \ref{13.3} is of the form \eqref{33} then $
\mes \( (E_n)_{+\de_n} \) = \eps_n + 2n\de_n $ (unless it exceeds $ 1
$), and \ref{13.3}(c) becomes
\begin{equation}\label{38}
\frac{ \eps_n + 2n\de_n }{ \be(\de_n^{-1}) } \to 0 \quad \text{for
some } (\de_n)_{n=1}^\infty \, .
\end{equation}
If we try $ \de_n = \eps_n / n $, the condition becomes $ \eps_n /
\be(n/\eps_n) \to 0 $, that is,
\begin{equation}\label{39}
\frac1{\eps_n} \be \bigg( \frac n {\eps_n} \bigg) \to \infty \, ,
\end{equation}
which implies the necessary condition, $ (1/\eps_n) \al(n/\eps_n) \not
\to 0 $, provided that $ \al(\cdot) \ge \be(\cdot) $.

Moreover, \eqref{39} is equivalent to \eqref{38}, if $ \be $ decreases,
but the function $ \la \mapsto \la \be(\la) $ increases (which is a
bit more than \ref{13.3}(a)). Indeed, then $ \( 1 +
\frac{\la_1}{\la_2} \) \frac{ \be(\la_1) }{ \be(\la_2) } \ge 1 $ for
all $ \la_1, \la_2 $ (check two cases, $ \la_1 \le \la_2 $ and $ \la_1
\ge \la_2 $); taking $ \la_1 = n/\eps_n $ and $ \la_2 = 1/\de_n $ we
get
\[
\frac{ \eps_n + n \de_n }{ \be(1/\de_n) } \cdot \frac{ \be (n/\eps_n)
}{ \eps_n } \ge 1 \, ;
\]
thus \eqref{38} implies \eqref{39}.

Given $ n_1 < n_2 < \dots $, we may consider the sequence $
(E_{n_k})_{k=1}^\infty $ (where $ E_n $ are given by \eqref{33}) and
ask whether it is infinitesimal, or not. Both Lemmas \ref{13.3} and
\ref{13.4} (and their proofs) still hold when $ n $ is substituted
with $ n_k $, and `$ n \to \infty $' with `$ k \to \infty $'.

\begin{lemma}\label{13.5}
\begin{sloppypar}
Assume that $ \nu_1, \nu_2 $ are two measures, each satisfying
\textup{\ref{9.5}(b)} and \eqref{29}; $ \al, \be : (0,\infty) \to
(0,\infty) $; $ \al $ is a slowly varying function; $ \be $ decreases,
but $ \la \mapsto \la^2 \be(\la) $ increases; $ \nu_1 (d\la) \le
\al(\la) \, d\la $; $ \nu_2 (d\la) \ge \be(\la) \, d\la $; and $
\limsup_{\la\to\infty} \( \be(\la) / \al(\la) \) = \infty $. Then
there exists a sequence of elementary sets that satisfies \eqref{34},
\eqref{35} w.r.t.\ $ \nu_2 $ but violates \eqref{35} w.r.t.\ $ \nu_1
$.
\end{sloppypar}
\end{lemma}

\begin{proof}
\begin{sloppypar}
We take $ \la_k \uparrow \infty $ such that $ \be(\la_k) / \al(\la_k)
\to \infty $. We note that $ \la \al(\la) \to \infty $; in particular,
$ \la_k \al(\la_k) \to \infty $. However, $ \la_k \be(\la_k) \gg \la_k
\al(\la_k) $ (that is, $ \( \la_k \be(\la_k) \) / \( \la_k \al(\la_k)
\) \to \infty $), and we can choose integers $ n_k \uparrow \infty $
such that
\[
\la_k \al(\la_k) \ll n_k \ll \la_k \be(\la_k) \, .
\]
Taking $ \eps_k = n_k / \la_k $ we have $ \eps_k \to 0 $ and
\[
\frac1{\eps_k} \al \bigg( \frac{ n_k }{ \eps_k } \bigg) \to 0 \, ,
\quad \frac1{\eps_k} \be \bigg( \frac{ n_k }{ \eps_k } \bigg) \to
\infty \, .
\]
The latter gives
\[
\frac{ \eps_k + 2n_k \de_k }{ \be(\de_k^{-1}) } \to 0 \quad \text{when
} k \to \infty
\]
for some $ \de_k $ (namely, for $ \de_k = \eps_k / n_k $, recall
\eqref{38}, \eqref{39}). Defining (recall \eqref{33})
\[
E_k = \bigcup_{l=0}^{n_k-1} \bigg( \frac1{n_k} \Big( l +
\frac{1-\eps_k}2 \Big), \, \frac1{n_k} \Big( l + \frac{1+\eps_k}2
\Big) \bigg) \, ,
\]
we have $ \mes E_k = \eps_k $ and $ \mes \( (E_k)_{+\de_k} \) = \eps_k
+ 2n_k \de_k $. Lemma \ref{13.3} (generalized to subsequences) states
that $ (E_k)_{k=1}^\infty $ satisfies \eqref{35} w.r.t.\ $ \nu_2 $. On
the other hand, Lemma \ref{13.4} (generalized to subsequences) shows
that $ (E_k)_{k=1}^\infty $ does not satisfy \eqref{35} w.r.t.\ $
\nu_1 $. It remains to note that \eqref{34} is ensured (w.r.t.\ $
\nu_2 $) by Lemma \ref{13.1}.
\end{sloppypar}
\end{proof}

The following result holds for both cases, real and
complex.\index{real versus complex}\index{complex versus real}

\begin{theorem}\label{13.6}
There is a continuum of mutually non-isomorphic product systems of
type $ III $.
\end{theorem}

\begin{proof}
For each $ \al \in (0,\infty) $ we choose a positive $ \sigma $-finite
measure $ \nu_\al $ on $ [0,\infty) $ such that the function $
\nu_\al (d\la) / d\la $ is strictly positive, continuously differentiable,
and for $ \la $ large enough, it is equal to $ \ln^{-\al} \la $ (recall
\eqref{293}). Lemma \ref{9.6} gives us the
corresponding continuous sum of FHS spaces. As was said in Sect.~9,
the continuous sum is homogeneous; Corollary \ref{9.4} converts it
into an \hcpHs. By Lemma \ref{10.2}, it is of type $ III $. The
corresponding (by Theorem \ref{1.5}) product system is also of type $
III $ (see \cite[Remark 3.3.3]{Ar}). Functions $ \la \mapsto
\ln^{-\al} \la $ are slowly varying, so, Lemma \ref{13.5} and Theorem
\ref{11.3} show that these \hcpHss\ are mutually non-isomorphic. By
Lemma \ref{1.6}, the corresponding product systems are mutually
non-isomorphic.
\end{proof}

\bibliographystyle{amsalpha}

\printindex

\bigskip
\filbreak
{
\small
\begin{sc}
\parindent=0pt\baselineskip=12pt
\parbox{2.3in}{
Boris Tsirelson\\
School of Mathematics\\
Tel Aviv University\\
Tel Aviv 69978, Israel
\smallskip
\emailwww{tsirel@tau.ac.il}
{www.tau.ac.il/\textasciitilde tsirel/}
}
\end{sc}
}
\filbreak

\end{document}